\magnification = \magstep 1
\catcode`\@=11
\font\tensmc=cmcsc10      
\def\smc{\tensmc}

\def\hcorrection#1{\advance\hoffset by #1 }
\def\vcorrection#1{\advance\voffset by #1 }
\def\wlog#1{}
\newif\iftitle@
\outer\def\title{\title@true\vglue 24\p@ plus 12\p@ minus 12\p@
   \bgroup\let\\=\cr\tabskip\centering
   \halign to \hsize\bgroup\tenbf\hfill\ignorespaces##\unskip\hfill\cr}
\def\endtitle{\cr\egroup\egroup\vglue 18\p@ plus 12\p@ minus 6\p@}
\outer\def\author{\iftitle@\vglue -18\p@ plus -12\p@ minus -6\p@\fi\vglue
    12\p@ plus 6\p@ minus 3\p@\bgroup\let\\=\cr\tabskip\centering
    \halign to \hsize\bgroup\smc\hfill\ignorespaces##\unskip\hfill\cr}
\def\endauthor{\cr\egroup\egroup\vglue 18\p@ plus 12\p@ minus 6\p@}
\outer\def\heading{\bigbreak\bgroup\let\\=\cr\tabskip\centering
    \halign to \hsize\bgroup\smc\hfill\ignorespaces##\unskip\hfill\cr}
\def\endheading{\cr\egroup\egroup\nobreak\medskip}

\outer\def\proclaim#1{\medbreak\noindent\smc\ignorespaces
    #1\unskip.\enspace\sl\ignorespaces}
\outer\def\endproclaim{\par\ifdim\lastskip<\medskipamount\removelastskip
  \penalty 55 \fi\medskip\rm}
\outer\def\demo#1{\par\ifdim\lastskip<\smallskipamount\removelastskip
    \smallskip\fi\noindent{\smc\ignorespaces#1\unskip:\enspace}\rm
      \ignorespaces}
\outer\def\enddemo{\par\smallskip}
\newcount\footmarkcount@
\footmarkcount@=1
\def\makefootnote@#1#2{\insert\footins{\interlinepenalty=100
  \splittopskip=\ht\strutbox \splitmaxdepth=\dp\strutbox
  \floatingpenalty=\@MM
  \leftskip=\z@\rightskip=\z@\spaceskip=\z@\xspaceskip=\z@
  \noindent{#1}\footstrut\rm\ignorespaces #2\strut}}
\def\footnote{\let\@sf=\empty\ifhmode\edef\@sf{\spacefactor
   =\the\spacefactor}\/\fi\futurelet\next\footnote@}
\def\footnote@{\ifx"\next\let\next\footnote@@\else
    \let\next\footnote@@@\fi\next}
\def\footnote@@"#1"#2{#1\@sf\relax\makefootnote@{#1}{#2}}
\def\footnote@@@#1{$^{\number\footmarkcount@}$\makefootnote@
   {$^{\number\footmarkcount@}$}{#1}\global\advance\footmarkcount@ by 1 }

\hyphenation{man-u-script man-u-scripts ap-pen-dix ap-pen-di-ces}
\hyphenation{data-base data-bases}
\ifx\amstexloaded@\relax\catcode`\@=13
  \endinput\else\let\amstexloaded@=\relax\fi
\newlinechar=`\^^J
\def\eat@#1{}
\def\Space@.{\futurelet\Space@\relax}
\Space@. %
\newhelp\athelp@
{Only certain combinations beginning with @ make sense to me.^^J
Perhaps you wanted \string\@\space for a printed @?^^J
I've ignored the character or group after @.}
\def\futureletnextat@{\futurelet\next\at@}
{\catcode`\@=\active
\lccode`\Z=`\@ \lowercase
{\gdef@{\expandafter\csname futureletnextatZ\endcsname}
\expandafter\gdef\csname atZ\endcsname
   {\ifcat\noexpand\next a\def\next{\csname atZZ\endcsname}\else
   \ifcat\noexpand\next0\def\next{\csname atZZ\endcsname}\else
    \def\next{\csname atZZZ\endcsname}\fi\fi\next}
\expandafter\gdef\csname atZZ\endcsname#1{\expandafter
   \ifx\csname #1Zat\endcsname\relax\def\next
     {\errhelp\expandafter=\csname athelpZ\endcsname
      \errmessage{Invalid use of \string@}}\else
       \def\next{\csname #1Zat\endcsname}\fi\next}
\expandafter\gdef\csname atZZZ\endcsname#1{\errhelp
    \expandafter=\csname athelpZ\endcsname
      \errmessage{Invalid use of \string@}}}}
\def\atdef@#1{\expandafter\def\csname #1@at\endcsname}
\newhelp\defahelp@{If you typed \string\define\space cs instead of
\string\define\string\cs\space^^J
I've substituted an inaccessible control sequence so that your^^J
definition will be completed without mixing me up too badly.^^J
If you typed \string\define{\string\cs} the inaccessible control sequence^^J
was defined to be \string\cs, and the rest of your^^J
definition appears as input.}
\newhelp\defbhelp@{I've ignored your definition, because it might^^J
conflict with other uses that are important to me.}
\def\define{\futurelet\next\define@}
\def\define@{\ifcat\noexpand\next\relax
  \def\next{\define@@}%
  \else\errhelp=\defahelp@
  \errmessage{\string\define\space must be followed by a control
     sequence}\def\next{\def\garbage@}\fi\next}
\def\undefined@{}
\def\preloaded@{}
\def\define@@#1{\ifx#1\relax\errhelp=\defbhelp@
   \errmessage{\string#1\space is already defined}\def\next{\def\garbage@}%
   \else\expandafter\ifx\csname\expandafter\eat@\string
         #1@\endcsname\undefined@\errhelp=\defbhelp@
   \errmessage{\string#1\space can't be defined}\def\next{\def\garbage@}%
   \else\expandafter\ifx\csname\expandafter\eat@\string#1\endcsname\relax
     \def\next{\def#1}\else\errhelp=\defbhelp@
     \errmessage{\string#1\space is already defined}\def\next{\def\garbage@}%
      \fi\fi\fi\next}
\def\famzero{\fam\z@}

\def\lim{\mathop{\famzero lim}}

\def\min{\mathop{\famzero min}}

\def\textfont@#1#2{\def#1{\relax\ifmmode
    \errmessage{Use \string#1\space only in text}\else#2\fi}}
\textfont@\rm\tenrm
\textfont@\it\tenit
\textfont@\sl\tensl
\textfont@\bf\tenbf
\textfont@\smc\tensmc
\let\ic@=\/
\def\/{\unskip\ic@}
\def\textfonti{\the\textfont1 }
\def\t#1#2{{\edef\next{\the\font}\textfonti\accent"7F \next#1#2}}
\let\B=\=
\let\D=\.
\def~{\unskip\nobreak\ \ignorespaces}
{\catcode`\@=\active
\gdef\@{\char'100 }}
\atdef@-{\leavevmode\futurelet\next\athyph@}
\def\athyph@{\ifx\next-\let\next=\athyph@@
  \else\let\next=\athyph@@@\fi\next}
\def\athyph@@@{\hbox{-}}
\def\athyph@@#1{\futurelet\next\athyph@@@@}
\def\athyph@@@@{\if\next-\def\next##1{\hbox{---}}\else
    \def\next{\hbox{--}}\fi\next}
\def\.{.\spacefactor=\@m}
\atdef@.{\null.}
\atdef@,{\null,}
\atdef@;{\null;}
\atdef@:{\null:}
\atdef@?{\null?}
\atdef@!{\null!}
\def\srdr@{\thinspace}
\def\drsr@{\kern.02778em}
\def\sldl@{\kern.02778em}
\def\dlsl@{\thinspace}
\atdef@"{\unskip\futurelet\next\atqq@}
\def\atqq@{\ifx\next\Space@\def\next. {\atqq@@}\else
         \def\next.{\atqq@@}\fi\next.}
\def\atqq@@{\futurelet\next\atqq@@@}
\def\atqq@@@{\ifx\next`\def\next`{\atqql@}\else\def\next'{\atqqr@}\fi\next}
\def\atqql@{\futurelet\next\atqql@@}
\def\atqql@@{\ifx\next`\def\next`{\sldl@``}\else\def\next{\dlsl@`}\fi\next}
\def\atqqr@{\futurelet\next\atqqr@@}
\def\atqqr@@{\ifx\next'\def\next'{\srdr@''}\else\def\next{\drsr@'}\fi\next}

\def\textfontii{\the\textfont2 }
\def\{{\relax\ifmmode\lbrace\else
    {\textfontii f}\spacefactor=\@m\fi}
\def\}{\relax\ifmmode\rbrace\else
    \let\@sf=\empty\ifhmode\edef\@sf{\spacefactor=\the\spacefactor}\fi
      {\textfontii g}\@sf\relax\fi}
\def\nonhmodeerr@#1{\errmessage
     {\string#1\space allowed only within text}}
\def\linebreak{\relax\ifhmode\unskip\break\else
    \nonhmodeerr@\linebreak\fi}
\def\allowlinebreak{\relax
   \ifhmode\allowbreak\else\nonhmodeerr@\allowlinebreak\fi}
\newskip\saveskip@
\def\nolinebreak{\relax\ifhmode\saveskip@=\lastskip\unskip
  \nobreak\ifdim\saveskip@>\z@\hskip\saveskip@\fi
   \else\nonhmodeerr@\nolinebreak\fi}
\def\newline{\relax\ifhmode\null\hfil\break
    \else\nonhmodeerr@\newline\fi}
\def\nonmathaerr@#1{\errmessage
     {\string#1\space is not allowed in display math mode}}
\def\nonmathberr@#1{\errmessage{\string#1\space is allowed only in math mode}}
\def\mathbreak{\relax\ifmmode\ifinner\break\else
   \nonmathaerr@\mathbreak\fi\else\nonmathberr@\mathbreak\fi}
\def\nomathbreak{\relax\ifmmode\ifinner\nobreak\else
    \nonmathaerr@\nomathbreak\fi\else\nonmathberr@\nomathbreak\fi}
\def\allowmathbreak{\relax\ifmmode\ifinner\allowbreak\else
     \nonmathaerr@\allowmathbreak\fi\else\nonmathberr@\allowmathbreak\fi}
\def\pagebreak{\relax\ifmmode
   \ifinner\errmessage{\string\pagebreak\space
     not allowed in non-display math mode}\else\postdisplaypenalty-\@M\fi
   \else\ifvmode\penalty-\@M\else\edef\spacefactor@
       {\spacefactor=\the\spacefactor}\vadjust{\penalty-\@M}\spacefactor@
        \relax\fi\fi}
\def\nopagebreak{\relax\ifmmode
     \ifinner\errmessage{\string\nopagebreak\space
    not allowed in non-display math mode}\else\postdisplaypenalty\@M\fi
    \else\ifvmode\nobreak\else\edef\spacefactor@
        {\spacefactor=\the\spacefactor}\vadjust{\penalty\@M}\spacefactor@
         \relax\fi\fi}
\def\newpage{\relax\ifvmode\vfill\penalty-\@M\else\nonvmodeerr@\newpage\fi}
\def\nonvmodeerr@#1{\errmessage
    {\string#1\space is allowed only between paragraphs}}
\def\smallpagebreak{\relax\ifvmode\smallbreak
      \else\nonvmodeerr@\smallpagebreak\fi}
\def\medpagebreak{\relax\ifvmode\medbreak
       \else\nonvmodeerr@\medpagebreak\fi}
\def\bigpagebreak{\relax\ifvmode\bigbreak
      \else\nonvmodeerr@\bigpagebreak\fi}
\newdimen\captionwidth@
\captionwidth@=\hsize
\advance\captionwidth@ by -1.5in
\def\caption#1{}
\def\topspace#1{\gdef\thespace@{#1}\ifvmode\def\next
    {\futurelet\next\topspace@}\else\def\next{\nonvmodeerr@\topspace}\fi\next}
\def\topspace@{\ifx\next\Space@\def\next. {\futurelet\next\topspace@@}\else
     \def\next.{\futurelet\next\topspace@@}\fi\next.}
\def\topspace@@{\ifx\next\caption\let\next\topspace@@@\else
    \let\next\topspace@@@@\fi\next}
 \def\topspace@@@@{\topinsert\vbox to
       \thespace@{}\endinsert}
\def\topspace@@@\caption#1{\topinsert\vbox to
    \thespace@{}\nobreak
      \smallskip
    \setbox\z@=\hbox{\noindent\ignorespaces#1\unskip}%
   \ifdim\wd\z@>\captionwidth@
   \centerline{\vbox{\hsize=\captionwidth@\noindent\ignorespaces#1\unskip}}%
   \else\centerline{\box\z@}\fi\endinsert}
\def\midspace#1{\gdef\thespace@{#1}\ifvmode\def\next
    {\futurelet\next\midspace@}\else\def\next{\nonvmodeerr@\midspace}\fi\next}
\def\midspace@{\ifx\next\Space@\def\next. {\futurelet\next\midspace@@}\else
     \def\next.{\futurelet\next\midspace@@}\fi\next.}
\def\midspace@@{\ifx\next\caption\let\next\midspace@@@\else
    \let\next\midspace@@@@\fi\next}
 \def\midspace@@@@{\midinsert\vbox to
       \thespace@{}\endinsert}
\def\midspace@@@\caption#1{\midinsert\vbox to
    \thespace@{}\nobreak
      \smallskip
      \setbox\z@=\hbox{\noindent\ignorespaces#1\unskip}%
      \ifdim\wd\z@>\captionwidth@
    \centerline{\vbox{\hsize=\captionwidth@\noindent\ignorespaces#1\unskip}}%
    \else\centerline{\box\z@}\fi\endinsert}
\mathchardef\prime@="0230
\def\prime{{{}\prime@{}}}
\def\prim@s{\prime@\futurelet\next\pr@m@s}

\def\,{\relax\ifmmode\mskip\thinmuskip\else\thinspace\fi}
\def\!{\relax\ifmmode\mskip-\thinmuskip\else\negthinspace\fi}
\def\frac#1#2{{#1\over#2}}

\def\:{\nobreak\hskip.1111em{:}\hskip.3333em plus .0555em\relax}
\def\intic@{\mathchoice{\hskip5\p@}{\hskip4\p@}{\hskip4\p@}{\hskip4\p@}}
\def\negintic@
 {\mathchoice{\hskip-5\p@}{\hskip-4\p@}{\hskip-4\p@}{\hskip-4\p@}}
\def\intkern@{\mathchoice{\!\!\!}{\!\!}{\!\!}{\!\!}}
\def\intdots@{\mathchoice{\cdots}{{\cdotp}\mkern1.5mu
    {\cdotp}\mkern1.5mu{\cdotp}}{{\cdotp}\mkern1mu{\cdotp}\mkern1mu
      {\cdotp}}{{\cdotp}\mkern1mu{\cdotp}\mkern1mu{\cdotp}}}
\newcount\intno@
\def\iint{\intno@=\tw@\futurelet\next\ints@}
\def\iiint{\intno@=\thr@@\futurelet\next\ints@}
\def\iiiint{\intno@=4 \futurelet\next\ints@}
\def\idotsint{\intno@=\z@\futurelet\next\ints@}
\def\ints@{\findlimits@\ints@@}
\newif\iflimtoken@
\newif\iflimits@
\def\findlimits@{\limtoken@false\limits@false\ifx\next\limits
 \limtoken@true\limits@true\else\ifx\next\nolimits\limtoken@true\limits@false
    \fi\fi}
\def\multintlimits@{\intop\ifnum\intno@=\z@\intdots@
  \else\intkern@\fi
    \ifnum\intno@>\tw@\intop\intkern@\fi
     \ifnum\intno@>\thr@@\intop\intkern@\fi\intop}
\def\multint@{\int\ifnum\intno@=\z@\intdots@\else\intkern@\fi
   \ifnum\intno@>\tw@\int\intkern@\fi
    \ifnum\intno@>\thr@@\int\intkern@\fi\int}
\def\ints@@{\iflimtoken@\def\ints@@@{\iflimits@
   \negintic@\mathop{\intic@\multintlimits@}\limits\else
    \multint@\nolimits\fi\eat@}\else
     \def\ints@@@{\multint@\nolimits}\fi\ints@@@}
\def\Sb{_\bgroup\vspace@
        \baselineskip=\fontdimen10 \scriptfont\tw@
        \advance\baselineskip by \fontdimen12 \scriptfont\tw@
        \lineskip=\thr@@\fontdimen8 \scriptfont\thr@@
        \lineskiplimit=\thr@@\fontdimen8 \scriptfont\thr@@
        \Let@\vbox\bgroup\halign\bgroup \hfil$\scriptstyle
            {##}$\hfil\cr}
\def\endSb{\crcr\egroup\egroup\egroup}
\def\Sp{^\bgroup\vspace@
        \baselineskip=\fontdimen10 \scriptfont\tw@
        \advance\baselineskip by \fontdimen12 \scriptfont\tw@
        \lineskip=\thr@@\fontdimen8 \scriptfont\thr@@
        \lineskiplimit=\thr@@\fontdimen8 \scriptfont\thr@@
        \Let@\vbox\bgroup\halign\bgroup \hfil$\scriptstyle
            {##}$\hfil\cr}
\def\endSp{\crcr\egroup\egroup\egroup}
\def\Let@{\relax\iffalse{\fi\let\\=\cr\iffalse}\fi}
\def\vspace@{\def\vspace##1{\noalign{\vskip##1 }}}
\def\aligned{\,\vcenter\bgroup\vspace@\Let@\openup\jot\m@th\ialign
  \bgroup \strut\hfil$\displaystyle{##}$&$\displaystyle{{}##}$\hfil\crcr}
\def\endaligned{\crcr\egroup\egroup}
\def\matrix{\,\vcenter\bgroup\Let@\vspace@
    \normalbaselines
  \m@th\ialign\bgroup\hfil$##$\hfil&&\quad\hfil$##$\hfil\crcr
    \mathstrut\crcr\noalign{\kern-\baselineskip}}
\def\endmatrix{\crcr\mathstrut\crcr\noalign{\kern-\baselineskip}\egroup
                \egroup\,}
\newtoks\hashtoks@
\hashtoks@={#}
\def\format{\crcr\egroup\iffalse{\fi\ifnum`}=0 \fi\format@}
\def\format@#1\\{\def\preamble@{#1}%
  \def\c{\hfil$\the\hashtoks@$\hfil}%
  \def\r{\hfil$\the\hashtoks@$}%
  \def\l{$\the\hashtoks@$\hfil}%
  \setbox\z@=\hbox{\xdef\Preamble@{\preamble@}}\ifnum`{=0 \fi\iffalse}\fi
   \ialign\bgroup\span\Preamble@\crcr}
 
\let\hdots=\ldots
\def\cases{\left\{\,\vcenter\bgroup\vspace@
     \normalbaselines\openup\jot\m@th
       \Let@\ialign\bgroup$##$\hfil&\quad$##$\hfil\crcr
      \mathstrut\crcr\noalign{\kern-\baselineskip}}

\newif\iftagsleft@
\tagsleft@true
\def\TagsOnRight{\global\tagsleft@false}
\def\tag#1$${\iftagsleft@\leqno\else\eqno\fi
 \hbox{\def\pagebreak{\global\postdisplaypenalty-\@M}%
 \def\nopagebreak{\global\postdisplaypenalty\@M}\rm(#1\unskip)}%
  $$\postdisplaypenalty\z@\ignorespaces}
\interdisplaylinepenalty=\@M
\def\allowdisplaybreak@{\def\allowdisplaybreak{\noalign{\allowbreak}}}
\def\displaybreak@{\def\displaybreak{\noalign{\break}}}
\def\align#1\endalign{\def\tag{&}\vspace@\allowdisplaybreak@\displaybreak@
  \iftagsleft@\lalign@#1\endalign\else
   \ralign@#1\endalign\fi}
\def\ralign@#1\endalign{\displ@y\Let@\tabskip\centering\halign to\displaywidth
     {\hfil$\displaystyle{##}$\tabskip=\z@&$\displaystyle{{}##}$\hfil
       \tabskip=\centering&\llap{\hbox{(\rm##\unskip)}}\tabskip\z@\crcr
             #1\crcr}}
\def\lalign@
 #1\endalign{\displ@y\Let@\tabskip\centering\halign to \displaywidth
   {\hfil$\displaystyle{##}$\tabskip=\z@&$\displaystyle{{}##}$\hfil
   \tabskip=\centering&\kern-\displaywidth
        \rlap{\hbox{(\rm##\unskip)}}\tabskip=\displaywidth\crcr
               #1\crcr}}
\def\overrightarrow{\mathpalette\overrightarrow@}
\def\overrightarrow@#1#2{\vbox{\ialign{$##$\cr
    #1{-}\mkern-6mu\cleaders\hbox{$#1\mkern-2mu{-}\mkern-2mu$}\hfill
     \mkern-6mu{\to}\cr
     \noalign{\kern -1\p@\nointerlineskip}
     \hfil#1#2\hfil\cr}}}
\def\overleftarrow{\mathpalette\overleftarrow@}
\def\overleftarrow@#1#2{\vbox{\ialign{$##$\cr
     #1{\leftarrow}\mkern-6mu\cleaders\hbox{$#1\mkern-2mu{-}\mkern-2mu$}\hfill
      \mkern-6mu{-}\cr
     \noalign{\kern -1\p@\nointerlineskip}
     \hfil#1#2\hfil\cr}}}
\def\overleftrightarrow{\mathpalette\overleftrightarrow@}
\def\overleftrightarrow@#1#2{\vbox{\ialign{$##$\cr
     #1{\leftarrow}\mkern-6mu\cleaders\hbox{$#1\mkern-2mu{-}\mkern-2mu$}\hfill
       \mkern-6mu{\to}\cr
    \noalign{\kern -1\p@\nointerlineskip}
      \hfil#1#2\hfil\cr}}}
\def\underrightarrow{\mathpalette\underrightarrow@}
\def\underrightarrow@#1#2{\vtop{\ialign{$##$\cr
    \hfil#1#2\hfil\cr
     \noalign{\kern -1\p@\nointerlineskip}
    #1{-}\mkern-6mu\cleaders\hbox{$#1\mkern-2mu{-}\mkern-2mu$}\hfill
     \mkern-6mu{\to}\cr}}}
\def\underleftarrow{\mathpalette\underleftarrow@}
\def\underleftarrow@#1#2{\vtop{\ialign{$##$\cr
     \hfil#1#2\hfil\cr
     \noalign{\kern -1\p@\nointerlineskip}
     #1{\leftarrow}\mkern-6mu\cleaders\hbox{$#1\mkern-2mu{-}\mkern-2mu$}\hfill
      \mkern-6mu{-}\cr}}}
\def\underleftrightarrow{\mathpalette\underleftrightarrow@}
\def\underleftrightarrow@#1#2{\vtop{\ialign{$##$\cr
      \hfil#1#2\hfil\cr
    \noalign{\kern -1\p@\nointerlineskip}
     #1{\leftarrow}\mkern-6mu\cleaders\hbox{$#1\mkern-2mu{-}\mkern-2mu$}\hfill
       \mkern-6mu{\to}\cr}}}
\def\sqrt#1{\radical"270370 {#1}}
\def\dots{\relax\ifmmode\let\next=\ldots\else\let\next=\tdots@\fi\next}
\def\tdots@{\unskip\ \tdots@@}
\def\tdots@@{\futurelet\next\tdots@@@}
\def\tdots@@@{$\mathinner{\ldotp\ldotp\ldotp}\,
   \ifx\next,$\else
   \ifx\next.\,$\else
   \ifx\next;\,$\else
   \ifx\next:\,$\else
   \ifx\next?\,$\else
   \ifx\next!\,$\else
   $ \fi\fi\fi\fi\fi\fi}
\def\text{\relax\ifmmode\let\next=\text@\else\let\next=\text@@\fi\next}
\def\text@@#1{\hbox{#1}}
\def\text@#1{\mathchoice
 {\hbox{\everymath{\displaystyle}\def\textfonti{\the\textfont1 }%
    \def\textfontii{\the\textfont2 }\textdef@@ T#1}}
 {\hbox{\everymath{\textstyle}\def\textfonti{\the\textfont1 }%
    \def\textfontii{\the\textfont2 }\textdef@@ T#1}}
 {\hbox{\everymath{\scriptstyle}\def\textfonti{\the\scriptfont1 }%
   \def\textfontii{\the\scriptfont2 }\textdef@@ S\rm#1}}
 {\hbox{\everymath{\scriptscriptstyle}\def\textfonti{\the\scriptscriptfont1 }%
   \def\textfontii{\the\scriptscriptfont2 }\textdef@@ s\rm#1}}}
\def\textdef@@#1{\textdef@#1\rm \textdef@#1\bf
   \textdef@#1\sl \textdef@#1\it}

\def\textdef@#1#2{\def\next{\csname\expandafter\eat@\string#2fam\endcsname}%
\if S#1\edef#2{\the\scriptfont\next\relax}%
 \else\if s#1\edef#2{\the\scriptscriptfont\next\relax}%
 \else\edef#2{\the\textfont\next\relax}\fi\fi}
\scriptfont\itfam=\tenit \scriptscriptfont\itfam=\tenit
\scriptfont\slfam=\tensl \scriptscriptfont\slfam=\tensl
\mathcode`\0="0030
\mathcode`\1="0031
\mathcode`\2="0032
\mathcode`\3="0033
\mathcode`\4="0034
\mathcode`\5="0035
\mathcode`\6="0036
\mathcode`\7="0037
\mathcode`\8="0038
\mathcode`\9="0039
\def\Cal{\relax\ifmmode\let\next=\Cal@\else
     \def\next{\errmessage{Use \string\Cal\space only in math mode}}\fi\next}
\def\Cal@#1{{\fam2 #1}}
\def\bold{\relax\ifmmode\let\next=\bold@\else
   \def\next{\errmessage{Use \string\bold\space only in math
      mode}}\fi\next}\def\bold@#1{{\fam\bffam #1}}
\mathchardef\Gamma="0000
\mathchardef\Delta="0001
\mathchardef\Theta="0002
\mathchardef\Lambda="0003
\mathchardef\Xi="0004
\mathchardef\Pi="0005
\mathchardef\Sigma="0006
\mathchardef\Upsilon="0007
\mathchardef\Phi="0008
\mathchardef\Psi="0009
\mathchardef\Omega="000A
\mathchardef\varGamma="0100
\mathchardef\varDelta="0101
\mathchardef\varTheta="0102
\mathchardef\varLambda="0103
\mathchardef\varXi="0104
\mathchardef\varPi="0105
\mathchardef\varSigma="0106
\mathchardef\varUpsilon="0107
\mathchardef\varPhi="0108
\mathchardef\varPsi="0109
\mathchardef\varOmega="010A
\font\dummyft@=dummy
\fontdimen1 \dummyft@=\z@
\fontdimen2 \dummyft@=\z@
\fontdimen3 \dummyft@=\z@
\fontdimen4 \dummyft@=\z@
\fontdimen5 \dummyft@=\z@
\fontdimen6 \dummyft@=\z@
\fontdimen7 \dummyft@=\z@
\fontdimen8 \dummyft@=\z@
\fontdimen9 \dummyft@=\z@
\fontdimen10 \dummyft@=\z@
\fontdimen11 \dummyft@=\z@
\fontdimen12 \dummyft@=\z@
\fontdimen13 \dummyft@=\z@
\fontdimen14 \dummyft@=\z@
\fontdimen15 \dummyft@=\z@
\fontdimen16 \dummyft@=\z@
\fontdimen17 \dummyft@=\z@
\fontdimen18 \dummyft@=\z@
\fontdimen19 \dummyft@=\z@
\fontdimen20 \dummyft@=\z@
\fontdimen21 \dummyft@=\z@
\fontdimen22 \dummyft@=\z@
\def\fontlist@{\\{\tenrm}\\{\sevenrm}\\{\fiverm}\\{\teni}\\{\seveni}%
 \\{\fivei}\\{\tensy}\\{\sevensy}\\{\fivesy}\\{\tenex}\\{\tenbf}\\{\sevenbf}%
 \\{\fivebf}\\{\tensl}\\{\tenit}\\{\tensmc}}
\def\dodummy@{{\def\\##1{\global\let##1=\dummyft@}\fontlist@}}
\newif\ifsyntax@
\newcount\countxviii@
\def\newtoks@{\alloc@5\toks\toksdef\@cclvi}
\def\nopages@{\output={\setbox\z@=\box\@cclv \deadcycles=\z@}\newtoks@\output}
\def\syntax{\syntax@true\dodummy@\countxviii@=\count18
\loop \ifnum\countxviii@ > \z@ \textfont\countxviii@=\dummyft@
   \scriptfont\countxviii@=\dummyft@ \scriptscriptfont\countxviii@=\dummyft@
     \advance\countxviii@ by-\@ne\repeat
\dummyft@\tracinglostchars=\z@
  \nopages@\frenchspacing\hbadness=\@M}
\def\magstep#1{\ifcase#1 1000\or
 1200\or 1440\or 1728\or 2074\or 2488\or
 \errmessage{\string\magstep\space only works up to 5}\fi\relax}
{\lccode`\2=`\p \lccode`\3=`\t
 \lowercase{\gdef\tru@#123{#1truept}}}

\def\scaletype#1{\mag=#1\relax
 \hsize=\expandafter\tru@\the\hsize
 \vsize=\expandafter\tru@\the\vsize
 \dimen\footins=\expandafter\tru@\the\dimen\footins}

\def\scalefont#1#2\andcallit#3{\edef\font@{\the\font}#1\font#3=
  \fontname\font\space scaled #2\relax\font@}
\def\Mag@#1#2{\ifdim#1<1pt\multiply#1 #2\relax\divide#1 1000 \else
  \ifdim#1<10pt\divide#1 10 \multiply#1 #2\relax\divide#1 100\else
  \divide#1 100 \multiply#1 #2\relax\divide#1 10 \fi\fi}
\def\scalelinespacing#1{\Mag@\baselineskip{#1}\Mag@\lineskip{#1}%
  \Mag@\lineskiplimit{#1}}
\def\wlog#1{\immediate\write-1{#1}}
\catcode`\@=\active


\TagsOnRight
\title
On Lunn-Senior's Mathematical Model \\
of  Isomerism in Organic Chemistry. Part I
\endtitle

\author
Valentin Vankov Iliev
\endauthor

\centerline {\it Section of Algebra, Institute of Mathematics and Informatics}
\centerline {\it Bulgarian Academy of Sciences, 1113 Sofia, Bulgaria}
\centerline {\it E-mail: viliev\@math.bas.bg}

\heading
1. Introduction
\endheading

In this introduction we summarize both the Lunn-Senior's mathematical model
from [3], and the content of the present paper.

1.1.  Let $A\!R$ be a set of atoms and radicals we are interesting in.  The
structural (connectivity) formula of a given chemical molecule
is usually drawn as
an $A\!R$-labelled graph $\Gamma$, where the labels of the vertices of
$\Gamma$
represent atoms or radicals from $A\!R$, and its (possibly multiple) edges
represent valences, or, equivalently, the connectivity data.  We note that
repetitions of labels are allowed.  In the sequel, we identify the graph
$\Gamma$ with the corresponding structural formula.  Following [3, I, p. 1030],
we use the terms \lq\lq structure" and \lq\lq connectivity"
as synonyms in order
to underline their independence of the $3$-dimensional space's limitations.

The mathematical model of Lunn and Senior, which is considered in [3], is based
on fixing a certain subset $U(\Gamma)$ of the set $v(\Gamma)$ of vertices of
$\Gamma$, which has the property that each vertex in $U(\Gamma)$ is
an endpoint of
exactly one edge of $\Gamma$.  The labels of the vertices from $U(\Gamma)$ are
called {\it univalent substituents} of $\Gamma$.  The subgraph $\Sigma(\Gamma)$
of $\Gamma$, with set of vertices $v(\Gamma)\backslash U(\Gamma)$ and all edges
that connect these vertices, is said to be the {\it skeleton}
of $\Gamma$.

Obviously, the division of a structural formula into skeleton and univalent
substituents is not unique, but once fixed, this division produces certain
properties of the molecule, which, after Lunn and Senior (see [3, I, p. 1031]),
are called {\it type properties}.

Given the skeleton $\Sigma=\Sigma(\Gamma)$, the \lq\lq degrees of freedom" of
the system are constituted by the various ways of distributing the univalent
substituents among the unsatisfied valences of the skeleton.  Let $d$ be the
number of univalent substituents of $\Gamma$.  We assign to each vertex of the
skeleton with unsatisfied valence a number from $1, 2,\hdots, d$, so that
different vertices have different numbers, and denote the set of these numbers
by $[1,d]$.  There are as many different $A\!R$-labelled graphs $\Gamma$ with a
fixed skeleton $\Sigma$, as maps $i\colon [1,d]\to A\!R$, $k\mapsto i_k$.
Thus, the Cartesian product $(A\!R)^d$ classifies the variety of all structural
formulae $\Gamma$ with a given skeleton $\Sigma$.  The combinatorial analysis
of these $\Gamma$'s is governed by the representation theory of the symmetric
group $S_d$ of the set $[1,d]$.

The fact that the univalent substituents consist of \lq\lq groups of like
individuals", and that \lq\lq ...the differences between them become
qualitative, like the differences between red, blue, and yellow geometrical
points" (see [3, I, p. 1031]), can be encoded in the mathematical model via
dissecting the set $[1,d]$ into several disjoint subsets $A_k$:
$[1,d]=\cup_kA_k$.  The group $S_d$ acts naturally on the set $\Delta_d$ of
all ordered dissections
$A=(A_1, A_2,\hdots, A_d)$ of the set $[1,d]$
by virtue of the rule
$$
\zeta A= (\zeta(A_1),
\zeta(A_2),\hdots, \zeta(A_d)).\tag 1.1.1
$$
Thus, we establish a monomial representation of the
symmetric group $S_d$.
We consider the subset $T_d$ of $\Delta_d$, consisting of all ordered
dissections $A$ whose components are ordered from largest to smallest.
Clearly the elements of the latter can be identified as tabloids with $d$
nodes (see [2, Ch. 2, 2.2.16]). Since $S_d$ is $d$-transitive on the set
$[1,d]$, there
exists a one-one correspondence between the orbit space $S_d\backslash T_d$ and
the set $P_d$ of all partitions
$\lambda=(\lambda_1,\lambda_2,\hdots,\lambda_d)$ of the positive integer $d$.
This correspondence can be obtained by factoring out the surjective map
$\varphi\colon T_d\to P_d$, $(A_1, A_2,\hdots, A_d)\mapsto
(\lambda_1,\lambda_2,\hdots,\lambda_d)$, where $\lambda_k$ is the
cardinality of the set $A_k$.
The $S_d$-orbit $T_\lambda$
corresponding to the partition $\lambda\in P_d$ consists of all tabloids of
shape $\lambda$.

Once a skeleton $\Sigma$ with $d$ unsatisfied valences is fixed, any tabloid
$A\in T_d$
can be considered as structural substituens'
{\it pre-formula} of the $d$ univalent substituents.  In other words, $A$ is a
pattern of maps which assigns to each number in the component $A_1$ of that
tabloid
$\lambda_1$ identical univalent substituents $x_1$
of type $1$, to each number in the component $A_2$ --- $\lambda_2$
identical univalent substituents $x_2$ of type $2$, etc., regardless of the
nature of these substituents. Moreover, there is a one-one correspondence
between the structural substituents' pre-formulae and the structural
pre-formulae obtained after joining the skeleton.  Then the monomial
$$
x_1^{\lambda_1}x_2^{\lambda_2}\hdots x_d^{\lambda_d},\tag 1.1.2
$$
where $\lambda$ is a partition of $d$,
represents the empirical
substituents'
pre-formula common to all structural substituents' pre-formulae from the
the set $T_\lambda$.

Throughout the rest of the paper, in any particular consideration the skeleton
will be fixed, so we shall use the expression \lq\lq structural (respectively,
empirical) pre-formula" for structural (respectively, empirical) substituents'
pre-formula, and shall identify this structural pre-formula with the
corresponding tabloid.

Introducing tabloids,
we avoid their equivalent but complicated set theoretic
interpretations used in [3, II--]. In particular, our approach allows us to
generalize for any partition $\lambda$ of $d$ the adjacency relations from [3,
VI], explicitly defined by Lunn and Senior only for the case
$\lambda_1+\lambda_2=d$.

\vskip 13pt

1.2.  A simple substitution reaction
$$
x_1^{\mu_1}\hdots x_i^{\mu_i}\hdots x_j^{\mu_j}\hdots\longrightarrow
x_1^{\lambda_1}\hdots x_i^{\lambda_i}\hdots x_j^{\lambda_j}\hdots,\tag 1.2.1
$$
where $\lambda, \mu\in P_d$, and
$\mu_1=\lambda_1,\hdots$, $\mu_i=\lambda_i+1,\hdots$,
$\mu_j=\lambda_j-1,\hdots$, $\mu_d=\lambda_d$,
is reflected by the mathematical model via introducing on the sets $P_d$ and
$T_d$ the so called simple raising operators $\rho_{i,j}$, and
$R_{i,s}$, respectively (see Sections 2, 3). The operator $R_{i,s}$
acts on a particular structural pre-formula
$A=(A_1, A_2,\hdots, A_d)\in T_\lambda$ by transferring the element $s\in A_j$
to $A_i$. This operator  mimics the inverse of the operation indicated in
the chemical equation (1.2.1): The replacement of one of the univalent
substituents $x_i$ of type $i$ in $A$ with an
univalent substituent $x_j$ of type $j$. The structural
pre-formula
$$
B=(A_1,\hdots, A_i\cup\{s\},\hdots, A_j\backslash\{s\},\hdots, A_d)=R_{i,s}A
$$
thus obtained is a tabloid of shape $\mu$, and $\lambda$ and $\mu$ are
connected via
the simple raising operator $\rho_{i,j}$ (see Section 2):
$\mu=\rho_{i,j}\lambda$.

A finite product $R$ (respectively, $\rho$) of simple raising operators
$R_{i,s}$ (respectively, $\rho_{i,j}$) is said to be a raising operator.
By means of these raising operators, we introduce
partial orders on the sets $T_d$ and $P_d$:
$$
A\leq B \hbox{\ if and only if there is a raising operator\ } R \hbox{\ with\ }
B=RA,\tag 1.2.2
$$
$$
\lambda\leq\mu \hbox{\ if and only if there is a raising operator\ } \rho
\hbox{\ with\ } \mu=\rho\lambda.
$$
The latter order is the famous
dominance order which plays an important role in the representation theory of
the symmetric group (see [2]). We note that in Sections 1 and 3 we state
equivalent definitions of the partial orders $A\leq B$, and $\lambda\leq\mu$,
respectively, which allow a direct check (in particular, by a computer).

\vskip 13pt

1.3.  Now, we turn our attention to the structural pre-formulae as arranged in
equivalence classes by certain isomeric relation.  In [3], Lunn and Senior
consider three isomeric relations:

(a) Univalent substitution isomerism;

(a$^\prime$) stereoisomerism;

(a$^{\prime\prime}$) structural isomerism.

The basic assumption of Lunn and
Senior in [3, III] is that for a fixed isomeric relation among
(a) -- (a$^{\prime\prime}$),
and
for a fixed skeleton $\Sigma$, there exists a permutation group $W\leq S_d$,
such that the corresponding isomeric classes can be identified with some
$W$-orbits in $T_d$.
More precisely, the group $W$ acts on the set $T_d$ via the rule (1.1.1), and
the isomers with skeleton $\Sigma$ and with $d$ univalent substituents are
identified with the elements of the set $T_{d;W}=W\backslash T_d$ of $W$-orbits
in $T_d$.

The authors emphasize that this group $W$ can be chosen from the large
selection of subgroups of $S_d$, using considerations which have nothing in
common with the $3$-dimensional space configuration of the respective
molecule.

The set $T_\lambda$ of tabloids of shape $\lambda$ is a disjoint union of
several $W$-orbits,
and if we denote the set of these $W$-orbits by $T_{\lambda;W}$, we have
$T_{d;W}=\cup_{\lambda\in P_d}T_{\lambda;W}$. It should be mentioned that in
the set
$T_{\lambda;W}$ are gathered all isomers with empirical pre-formula (1.1.2).
Let
$n_{\lambda;W}$
be the number of elements of the set $T_{\lambda;W}$.

\vskip 13pt
1.4. Let us consider the partial order on $T_{d;W}$ obtained by
factoring-out the partial order (1.2.2) in $T_d$:
For $a$, $b\in T_{d;W}$, we write
$$
a\leq b \hbox{\sl\ if and only if\ } A\leq B \hbox{\sl\ for some\ } A\in a
\hbox{\sl\ and\ } B\in b.
$$

This partial order on $T_{d;W}$ is a natural generalization of the adjacency
relations considered in [3, VI], so it is a mathematical model of the
genetic relations among isomers in organic chemistry.

For any couple $\lambda$, and $\mu$ of adjacent partitions with
$\lambda<\mu$, and $\mu=\rho_{i,j}\lambda$,
we consider the subset $R_{\lambda,\mu;W}\subset
T_{\lambda;W}\times T_{\mu;W}$, consisting of all ordered pairs $(a,b)$ such
that $a<b$, and set
$t_{\lambda,\mu;W}=|R_{\lambda,\mu;W}|$.

\vskip 13pt

1.5. Now, we shall enunciate the main statements of
Lunn and Senior from [3], summarized in the following

\proclaim{1.5.1. Lunn-Senior's thesis} Let $\Sigma$ be a skeleton with $d$
unsatisfied single valences.
One considers  molecules with skeleton $\Sigma$ and substitution's structural
pre-formulae which have empirical formula (1.1.2).
Then

1. There exist three permutation groups $G$, $G^\prime$, $G^{\prime\prime}\leq
S_d$, such that:

(1a) Any univalent substitution isomer can be identified with  a $G$-orbit in
$T_d$;

(1a$^\prime$) any stereoisomer can be identified with  a
$G^\prime$-orbit in $T_d$;

(1a$^{\prime\prime}$) any structural isomer can be identified with
a $G^{\prime\prime}$-orbit in $T_d$.

2. The groups
$G$, $G^\prime$, and $G^{\prime\prime}\leq
S_d$, are connected in the following way:

(2a) $G=G^\prime$, in case there are no chiral pairs among the
univalent substitution isomers, and $G\leq G^\prime$ with $|G^\prime:G|=2$, in
case there are such pairs. In the first case,
the $G$- and $G^\prime$-orbits coincide and some of them inventory the
diastereomers.
In the last case, each $G^\prime$-orbit contains either

(2ae) two $G$-orbits, and the members of any
chiral pair are represented by such a couple of $G$-orbits,

or, coincide with

(2ad) one $G$-orbit, and any diastereomer is represented by such a $G$-orbit.

(2b) Any $G^{\prime\prime}$-orbit is a disjoint union of $G^\prime$-orbits.

3. Each simple substitution reaction $b\to a$ of the type (1.2.1) can be
identified with the element $(a,b)\in R_{\lambda,\mu;G}$.

4.  The terms and relations involved in the statements 1 -- 3 do not depend on
the nature of the univalent substituents, so they represent type properties of
the molecules under consideration.

\endproclaim

\proclaim{Remark 1.5.2} {\rm The chemical discourse which has resource to the
experiment, and Lunn-Senior's mathematical model, create two languages showing
some discrepancy. Below, we state explicitly the chemical definitions of the
different types of isomerism described
by the mathematical model, in terms of
the model itself. Any two compounds in a particular definition are supposed to
have the same
empirical formula, that is, the corresponding tabloids have the same shape.

Two chemical compounds are said to be {\it structural isomers} if the
$G^{\prime\prime}$-orbits of their structural formulae are different.

Two chemical compounds are called {\it stereoisomers} if the
$G^\prime$-orbits of their structural pre-formulae are different, but are
contained in the same $G^{\prime\prime}$-orbit (that is, they have the same
connectivity data).

Two chemical compounds are said to be {\it univalent substitution isomers} if
the $G$-orbits of their structural formulae are different.

Two chemical compounds are said to form an {\it chiral pair}
if the $G$-orbits of their structural formulae are different, but are contained
in, and cover the same $G^\prime$-orbit (in particular, they represent the same
stereoisomer).

Two chemical compounds are said to be {\it diastereomers}
if: (a) the $G$-orbits $O_1$ and $O_2$ of their structural formulae are
different; (b)
each of $O_1$ and $O_2$ coincide with the corresponding
$G^\prime$-orbit; (c)
both $O_1$ and $O_2$  are contained in the same
$G^{\prime\prime}$-orbit ((a) - (c) yield that $O_1$ and $O_2$ are
stereoisomers).}

\endproclaim

Let
$N_{\lambda;\Sigma}$ (respectively, $N_{\lambda;\Sigma}^\prime$,
$N_{\lambda;\Sigma}^{\prime\prime}$)
be the number of univalent substitution isomers
(respectively, stereoisomers, structural isomers) with fixed skeleton
$\Sigma$, which have empirical
pre-formula (1.1.2).
Let
$T_{\lambda,\mu;\Sigma}$
be the number of different simple
substitution reactions of the type (1.2.1) among
the univalent substitution isomers
with that skeleton $\Sigma$.

According to Lunn-Senior's
thesis we have as consequences the following
inequalities:
$$
N_{\lambda;\Sigma}\leq n_{\lambda;G}, \hbox{\ }
T_{\lambda,\mu;\Sigma}\leq t_{\lambda,\mu;G}, \hbox{\ }\lambda\in P_d, \tag
1.5.3
$$
$$
N_{\lambda;\Sigma}^\prime\leq n_{\lambda;G^\prime},\hbox{\ }\lambda\in P_d,\tag
1.5.4
$$
and
$$
N_{\lambda;\Sigma}^{\prime\prime}\leq n_{\lambda;G^{\prime\prime}},
\hbox{\ }\lambda\in P_d. \tag 1.5.5
$$

The above inequalities can be used to find the group which corresponds to the
particular type of isomerism,
as Lunn-Senior's thesis asserts:  If one of the
inequalities from a row is false for a particular subgroup of the symmetric
group $S_d$,
then this subgroup has to be rejected (see  [3, IV]).  On the other hand,
Theorem 5.2.5 shows that the family $(n_{\lambda;W})_\lambda$ of non-negative
integers defines both the permutation group $W\leq S_d$ up to combinatorial
equivalence, and the corresponding induced monomial representation
$Ind_W^{S_d}(1_W)$ of the symmetric group $S_d$ --- up to isomorphism
(here $W$ is one of the groups $G$, $G^\prime$, or $G^{\prime\prime}$).

\vskip 13pt

1.6. A disadvantage of Lunn-Senior's mathematical model is that there are no
enough tools immanent to it, in order for two $W$-orbits to be distinguished.

The aim of this article is to present a mathematical
formalism which includes Lunn-Senior's model as a particular case and makes use
of the one-dimensional characters of the group $W$, and the one-dimensional
characters of the group
$S_\lambda=S_{\lambda_1}\times S_{\lambda_2}\times\cdots\leq S_d$,
for picking out of some special $W$-orbits.
A point of departure is the
following observation. Let us suppose that there are chiral pairs
among the stereoisomers of a given molecule with empirical formula (1.1.2).
Then, according to
Lunn-Senior's thesis 1.5.1, the group $G$ is a (normal) subgroup of $G^\prime$
with
$|G^\prime:G|=2$. Let $\chi_e\colon G^\prime\to \{1,-1\}$ be the homomorphism
of groups, which assigns  $1$ to each element of $G$, and
$-1$ to each element of the complement $G^\prime\backslash G$ of $G$.
Each
$G^\prime$-orbit (which, at least potentially, represents a stereoisomer)
either coincide
with the corresponding $G$-orbit, (and potentially represents a
diastereomer) or splits into two $G$-orbits (thus potentially representing an
chiral pair). The $G^\prime$-orbits $O$ which consist of two
$G$-orbits
can be distinguished from the other $G^\prime$-orbits in the following way.
Suppose that $A\in O$ is a tabloid,
and let
$G_A^\prime$ be the stabilizer of $A$ in $G^\prime$.
We can
consider $\chi_e$ as a one-dimensional character
$\chi_e\colon G^\prime\to K$, where $K$ is the field of complex numbers. Then
$O$
splits into two $G$-orbits if and only if the character $\chi_e$ is
identically
$1$ on the subgroup $G_A^\prime$. We can count the number of those
$G^\prime$-orbits (let us call them {\it $\chi_e$-orbits}), using the machinery
developed in Section 5. Thus, the one-dimensional character $\chi_e$ of
the group $G^\prime$ produces a type property of the molecule in question.

On the other hand, it is well known that
there is a one-one
correspondence between the set $T_{\lambda;W}$ of all $W$-orbits in
$T_\lambda$, and the set of all double cosets of $S_d$ modulo $(W,S_\lambda)$.
Let $\theta$ be a one-dimensional character of the group
$S_\lambda$, and let $\chi$ be a one-dimensional character of $W$. We
consider the subset
$T_{\lambda;\chi,\theta}$
of the set
$T_{\lambda;W}$,  consisting of all $W$-orbits which satisfy
property (5.1.3), (call them $(\chi,\theta)$-{\it orbits}), and set
$n_{\lambda;\chi,\theta}=
|T_{\lambda;\chi,\theta}|$.

The hypothesis that for any pair $(\chi,\theta)$, where $W$ is
a group among $G$, $G^\prime$, and $G^{\prime\prime}$, the property (5.1.3)
is a type property of the corresponding molecule, recognizable by an
experiment, yields the following

\proclaim{1.6.1. Extended Lunn-Senior's thesis} Let $\Sigma$ be a skeleton with
$d$ unsatisfied single valences.
One considers  molecules with skeleton $\Sigma$ and substitution's structural
pre-formulae which have empirical formula (1.1.2).
Then

1. There exist three permutation groups $G$, $G^\prime$, $G^{\prime\prime}\leq
S_d$, such that:

(1a) Any univalent substitution isomer can be identified with  a $G$-orbit in
$T_d$;

(1a$^\prime$) any stereoisomer can be identified with a
$G^\prime$-orbit in $T_d$;

(1a$^{\prime\prime}$) any structural isomer can be identified with
a $G^{\prime\prime}$-orbit in $T_d$.

2. The groups
$G$, $G^\prime$, and $G^{\prime\prime}\leq
S_d$, are connected in the following way:

(2a) $G=G^\prime$, in case there are no chiral pairs among the
univalent substitution isomers, and $G\leq G^\prime$ with $|G^\prime:G|=2$, in
case there are such pairs. In the first case,
the $G$- and $G^\prime$-orbits coincide and some of them inventory the
diastereomers.
In the last case, each $G^\prime$-orbit contains either

(2ae) two $G$-orbits, and the members of any
chiral pair are represented by such a couple of $G$-orbits,

or, coincide with

(2ad) one $G$-orbit, and any diastereomer is represented by such a $G$-orbit.

The $\chi_e$-orbits are those $G^\prime$-orbits which represent the
chiral pairs.

(2b) Any $G^{\prime\prime}$-orbit is a disjoint union of $G^\prime$-orbits.

3. For each sequence $b\to\cdots\to a$ of simple substitution reactions
one has $a<b$ and the reaction $b\to a$ can be identified with the inequality
$a<b$ in $T_{d;G}$.

4.  The terms and relations involved in the statements 1 -- 3 do not depend on
the nature of the univalent substituents, so they represent type properties of
the molecules under consideration.

5. If $\theta$ is a one-dimensional character of the group $S_\lambda$,
and $\chi$ is a one-dimensional character of the group $W$, where $W$ is
one of $G$, $G^\prime$, or $G^{\prime\prime}$, then the set of
all $(\chi,\theta)$-orbits
of $W$ in $T_\lambda$ represents a type property of the molecule.

\endproclaim

The isomers which correspond to the hypothetical type property from 1.6.1, item
5, are called
$(\chi,\theta)$-{\it isomers.}
Let
$N_{\lambda;\chi,\theta;\Sigma}$
be the number of all
$(\chi,\theta)$-isomers
with fixed skeleton
$\Sigma$.

As far as Extended Lunn-Senior's
thesis is valid, we have the
inequalities
$$
N_{\lambda;\chi,\theta;\Sigma}\leq
n_{\lambda;\chi,\theta}.
$$

\vskip 13pt

1.7. In Section 2 we consider the dominance order on the set $M_d$
consisting of all $d$-tuples $m=(m_1,\hdots,m_d)$ of
non-negative integers
whose sum is $d$,
(see [2, Ch. 1, 1.4.6])
and gather the necessary
information concerning neighbourhood in $M_d$ and in its subset $P_d$ of
all partitions of $d$.

In Section 3 we introduce tabloids
and raising
operators which act on their set $T_d$ by analogy with
the
raising operators from Section 2. Inasmuch as possible, we work in the wider
set $\Delta_d$,
consisting of all ordered dissections $A=(A_1,\hdots,A_d)$ of
the set $[1,d]$.
We provide the set $\Delta_d$ with a partial order (also called dominance
order) such that if we consider the dominance order on the set $M_d$,
then the map $\varphi\colon \Delta_d\to M_d$ from (3.1.1)
is a homomorphism of partially ordered sets.
The main objective in Section 3 is the study of the equation $\varphi(X)=n$,
where $n\in M_d$ (respectively, $n\in P_d$), and the unknown $X$ varies in an
interval
$[A,B]$ in $\Delta_d$ (respectively, in $T_d$). Theorem 3.4.3 allows us
to establish
Theorem 3.5.1 which is a criterion for two ordered dissections (tabloids) $A$
and $B$ to be neighbours with respect to the corresponding partial order.
This is done by a systematical use of raising operators.

In Section 4 we factor out the constructions from Section 3 with respect to the
action of a permutation group $W\leq S_d$, and produce the sets
$\Delta_{d;W}$ and $T_{d;W}$,
the last one being the sphere of action
of the generalized Lunn-Senior's
mathematical model of isomerism. Note especially Theorem
4.2.1 which
gives necessary and sufficient conditions for two elements $a$ and $b$ to be
adjacent in $\Delta_{d;W}$
(respectively, in $T_{d;W}$), as well as Theorem 4.2.3
which is a criterion for $a$ and $b$ to be neighbours there.

Section 5 is devoted to finding explicit expressions for the maximum number of
isomers under consideration, according to Lunn-Senior's thesis 1.5.1 and its
extension 1.6.1. Here Theorem 5.2.7 is the central result.
In Corollary 5.2.10 we give another proof of
Ruch's formula which establishes an explicit expression for the numbers
$n_{\lambda;W}$ (see [6]).
We have to point out Lemma 5.4.3 which shows that when
$\theta$
is the unit character of the group $S_\lambda$, the abstract condition (5.1.3)
on the stabilizer $W_A$ of an ordered dissection $A\in a$ is equivalent to
to the following maximum property of the $W$-orbit $a$:
$$
\hbox{\sl \lq\lq the\ }
W-\hbox{\sl orbit\ } a
\hbox{\sl\ consists of\ } |W:W_\chi|
\hbox{\sl\ in number\ } W_\chi-\hbox{\sl orbits"},
$$
where $W_\chi\leq W$ is the kernel of the one-dimensional character $\chi\colon
W\to K$.

Theorem 5.3.1 is a generalization of an important result of E.  Ruch which
connects
the dominance order on the set $P_d$ and the existence of chiral
pairs, as it is shown in Subsection 6.2, Theorem 6.2.1.  The rest of
Section 6 contains
illustrations of our approach applied on well known examples:  A proof of
Kauffmann formulae for the derivatives of naphthalene, and inferences of the
genetic relations of ethene and benzene.

\heading{2. Partitions}

\endheading

2.1.
Let $N_d$ be the set of all $d$-tuples
$m=(m_1,\hdots ,m_d)$
of integers $m_j$ with $\sum_jm_j=d$.
Let $M_d$ be the subset of $N_d$ consisting of all $d$-tuples $m$ with
non-negative components. We denote by $P_d$ the
subset of $M_d$ whose elements are all $\lambda=(\lambda_1,\hdots
,\lambda_d)\in M_d$ with
$\lambda_1\geq\hdots\geq\lambda_d$. The elements of $P_d$ are called {\it
partitions of $d$}.
The partition $\lambda$ can be visualized by the
corresponding Young diagram:

$$
\lambda=
\left.
\matrix \format \c & \quad \c & \quad \c & \quad \c & \quad \c & \quad \c &
\quad \c & \quad \c \\
\times & \times & \cdots & \cdots & \cdots & \times &\hbox{\ \ } &
\lambda_1\hbox{\ nodes\ }   \cr
\times & \times & \cdots & \cdots & \times & \hbox{ \ \ } & \hbox{ \ \
} & \lambda_2\hbox{\ nodes\ }   \cr
\hbox{\ \ } & \hbox{\ \ } & \vdots & \hbox{\ \ } & \hbox{\ \ } & \hbox{\ \ }
& \hbox{\ \ } & \vdots   \cr
\times & \cdots & \hbox{\ \ } & \hbox{\ \ } & \hbox{\ \ } & \hbox{\ \ }  &
\hbox{\ \ } & \lambda_t\hbox{\ nodes\ }   \cr
\endmatrix
\right.\tag 2.1.1
$$
where $\lambda_t$ is the last nonzero component of $\lambda$.
Let
$l=(l_1,\hdots ,l_d)$ and
$m=(m_1,\hdots ,m_d)$ be two elements of the set $M_d$. In case $l_1=m_1$,
we denote by $q(l,m)$ the maximum number $q\in [1,d]$ such that
$l_1=m_1,\hdots, l_q=m_q$.
Otherwise, we set $q(l,m)=0$.

Let $\leq$ be the {\it dominance order} on $N_d$ (see [4, Ch. I, Sec. 1]).
We remind that
$l\leq m$ if and
only if $\sum_{k=1}^il_k\leq\sum_{k=1}^i m_k$
for any $1\leq i\leq d$. In this case we say that $m$ {\it dominates} $l$.
It is clear that $\leq$ is a partial order on $N_d$ (see Appendix A) which
induces
partial orders on $M_d$ and in $P_d$, the last two being denoted by the same
sign and
also named dominance order. Below, the dominance order on $P_6$ is graphically
portrayed.

$$
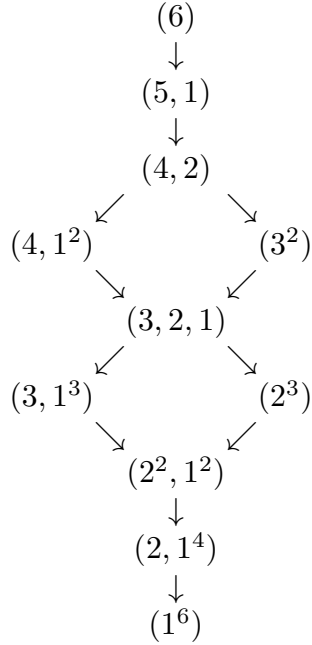
\left.
\matrix \format \r & \r & \r & \c & \l & \l & \l & \l & \l & \l & \l & \l & \l
& \l &  \l & \l & \l \\
\hbox{\ \ } & \hbox{\ \ } & \hbox{\ \ } &
(6) & \hbox{\ \ } & \hbox{\ \ } & \hbox{\ \ } \cr
\hbox{\ \ } & \hbox{\ \ } & \hbox{\ \ } & \downarrow & \hbox{\ \
} & \hbox{\ \ } & \hbox{\ \ } \cr
\hbox{\ \ } & \hbox{\ \ } & \hbox{\ \ } &
(5,1) & \hbox{\ \ } & \hbox{\ \ } & \hbox{\ \ } \cr
\hbox{\ \ } & \hbox{\ \ } & \hbox{\ \ } & \downarrow & \hbox{\ \
} &\hbox{\ \ } & \hbox{\ \ } \cr
\hbox{\ \ } & \hbox{\ \ } & \hbox{\ \ } &
(4,2) & \hbox{\ \ } & \hbox{\ \ } & \hbox{\ \ }  \cr
\hbox{\ \ } & \hbox{\ \ } & \swarrow & \hbox{\ \ }
& \searrow & \hbox{\ \ } & \hbox{\ \ } \cr
\hbox{\ \ } & (4,1^2) &  \hbox{\ \ } & \hbox{\ \ } & \hbox{\ \ } &
(3^2) &\hbox{\ \ }   \cr
\hbox{\ \ } & \hbox{\ \ } & \searrow & \hbox{\ \ } & \swarrow
& \hbox{\ \ } & \hbox{\ \ } \cr
\hbox{\ \ } & \hbox{\ \ } & \hbox{\ \ } & (3,2,1) & \hbox{\ \ } &
\hbox{\ \ } & \hbox{\ \ } \cr
\hbox{\ \ } & \hbox{\ \ } & \swarrow & \hbox{\ \ } & \searrow
& \hbox{\ \ } & \hbox{\ \ } \cr
\hbox{\ \ } & (3,1^3) &  \hbox{\ \ } & \hbox{\ \ }
& \hbox{\ \ } & (2^3) & \hbox{\ \ }
\cr
\hbox{\ \ } & \hbox{\ \ } & \searrow & \hbox{\ \ } & \swarrow
& \hbox{\ \ } & \hbox{\ \ } \cr
\hbox{\ \ } & \hbox{\ \ } & \hbox{\ \ } & (2^2,1^2) & \hbox{\ \ } & \hbox{\ \ }
& \hbox{\ \ } \cr
\hbox{\ \ } & \hbox{\ \ } & \hbox{\ \ } & \downarrow & \hbox{\ \ }
& \hbox{\ \ } & \hbox{\ \ } \cr
\hbox{\ \ } & \hbox{\ \ } & \hbox{\ \ } & (2,1^4) & \hbox{\ \ }
& \hbox{\ \ } & \hbox{\ \ }  \cr
\hbox{\ \ } & \hbox{\ \ } & \hbox{\ \ } & \downarrow & \hbox{\ \ }
& \hbox{\ \ } & \hbox{\ \ } \cr
\hbox{\ \ } & \hbox{\ \ } & \hbox{\ \ } & (1^6) & \hbox{\ \ } & \hbox{\ \ }
& \hbox{\ \ } \cr
\endmatrix
\right.
$$

Given $i$, $j\in [1,d]$, we
define an operator $\rho_{i,j}\colon N_d\to N_d$
by the formulae

$$
\rho_{i,j}(l)=
\left\{
\matrix \format \c & \c & \c & \c \\
(l_1,\hdots ,l_i+1,\hdots,l_j-1,\hdots, l_d) & \hbox{\sl\ \ if\ \ } &
1\leq i<j\leq d\cr
l & \hbox{\sl \ \ if\ \ } & 1\leq j\leq i\leq d.\cr
\endmatrix
\right.
$$

The operators $\rho_{i,j}$ are called {\it simple raising operators}
in $N_d$.
Obviously, any two simple
raising operators $\rho_{i_1,j_1}$ and $\rho_{i_2,j_2}$ in $N_d$ commute.
Any product
$\rho=\rho_{i_1,j_1}\rho_{i_2,j_2}\hdots$
of simple raising operators
is called {\it raising operator} in the set $N_d$. A raising operator is said
to be {\it non-trivial} if it moves at least one element in $N_d$.
Otherwise, it is called {\it trivial}.

\proclaim{Remark 2.1.2} {\rm We note that the subsets $M_d$, and $P_d$ of
$N_d$ are not closed with respect to the action of a non-trivial raising
operator: Given a $d$-tuple $l\in M_d$, and a simple raising
operator $\rho_{i,j}$ with $1\leq i<j\leq d$, one has $\rho_{i,j}(l)\in M_d$
if and only if $l_j\geq 1$. }

\endproclaim

If $1\leq i<j\leq d$, and $\mu=\rho_{i,j}\lambda\in P_d$, then the picture
below illustrates the action of the simple raising operator $\rho_{i,j}$
on the Young diagram representing the partition $\lambda$:

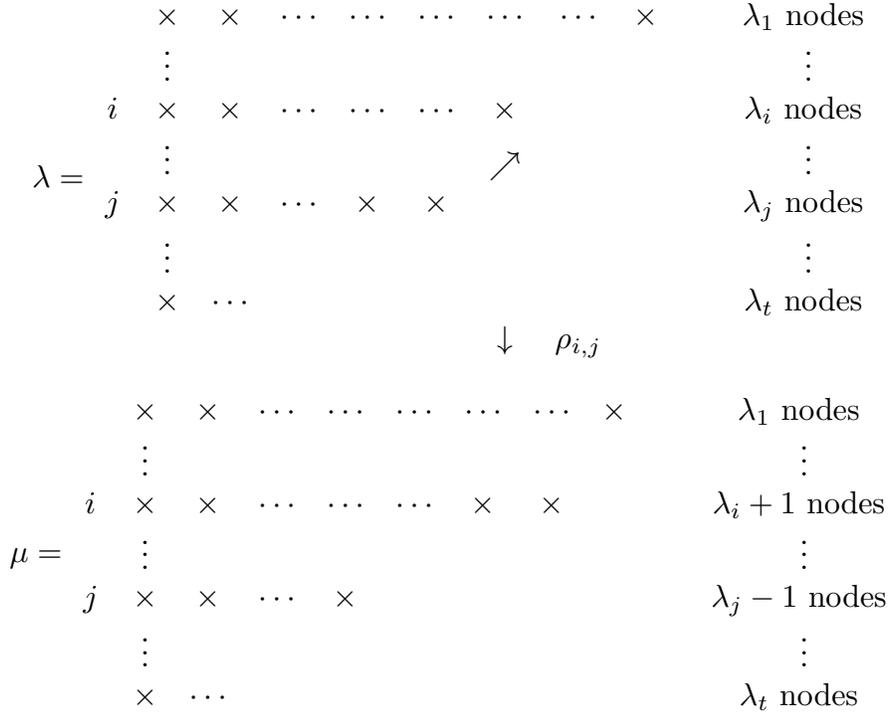
$$
\lambda=
\left.
\matrix \format \c & \quad \c & \quad \c & \quad \c & \quad \c & \quad \c &
\quad \c & \quad \c & \quad \c & \quad \c & \quad \c \\
\hbox{\ \ } & \times & \times & \cdots & \cdots & \cdots & \cdots & \cdots
& \times & \hbox{\ \ } & \lambda_1 \hbox{\ nodes\ }   \cr
\hbox{\ \ } & \vdots & \hbox{\ \ } & \hbox{\ \ } & \hbox{\ \ } &
\hbox{\ \ } & \hbox{ \ \ } & \hbox{ \ \ } & \hbox{\ \ } & \hbox{\ \ } &
\vdots   \cr
i & \times & \times & \cdots & \cdots & \cdots & \times & \hbox{\ \ }
& \hbox{\ \ } & \hbox{\ \ } & \lambda_i \hbox{\ nodes\ }   \cr
\hbox{\ \ } & \vdots & \hbox{\ \ } & \hbox{\ \ } & \hbox{\ \ } &
\hbox{\ \ } & \nearrow & \hbox{ \ \ } & \hbox{\ \ } & \hbox{\ \ } &
\vdots   \cr
j & \times & \times & \cdots & \times & \times & \hbox{\ \ } & \hbox{\ \ }
& \hbox{\ \ } & \hbox{\ \ } & \lambda_j \hbox{\ nodes\ }   \cr
\hbox{\ \ } & \vdots & \hbox{\ \ } & \hbox{\ \ } & \hbox{\ \ } &
\hbox{\ \ } & \hbox{ \ \ } & \hbox{ \ \ } & \hbox{\ \ } & \hbox{\ \ } &
\vdots   \cr
\hbox{\ \ } & \times & \cdots & \hbox{\ \ } & \hbox{\ \ } & \hbox{\ \ } &
\hbox{\ \
} & \hbox{\ \ } & \hbox{\ \ } & \hbox{\ \ } & \lambda_t \hbox{\ nodes\ }   \cr
\hbox{\ \ } & \hbox{\ \ } & \hbox{\ \ } & \hbox{\ \ } & \hbox{\ \ } &
\hbox{\ \ } & \downarrow & \rho_{i,j} & \hbox{\ \ } & \hbox{\ \ }
& \hbox{\ \ } \cr
\endmatrix
\right.
$$
$$
\mu=
\left.
\matrix \format \c & \quad \c & \quad \c & \quad \c & \quad \c & \quad \c &
\quad \c & \quad \c & \quad \c & \quad \c & \quad \c \\
\hbox{\ \ } & \times & \times & \cdots & \cdots & \cdots & \cdots & \cdots
& \times & \hbox{\ \ } & \lambda_1 \hbox{\ nodes\ }   \cr
\hbox{\ \ } & \vdots & \hbox{\ \ } & \hbox{\ \ } & \hbox{\ \ } &
\hbox{\ \ } & \hbox{ \ \ } & \hbox{ \ \ } & \hbox{\ \ } & \hbox{\ \ } &
\vdots   \cr
i & \times & \times & \cdots & \cdots & \cdots & \times & \times
& \hbox{\ \ } & \hbox{\ \ } & \lambda_i+1 \hbox{\ nodes\ }   \cr
\hbox{\ \ } & \vdots & \hbox{\ \ } & \hbox{\ \ } & \hbox{\ \ } &
\hbox{\ \ } & \hbox{\ \ } & \hbox{ \ \ } & \hbox{\ \ } & \hbox{\ \ } &
\vdots   \cr
j & \times & \times & \cdots & \times & \hbox{\ \ } & \hbox{\ \ } & \hbox{\ \ }
& \hbox{\ \ } & \hbox{\ \ } & \lambda_j-1 \hbox{\ nodes\ }   \cr
\hbox{\ \ } & \vdots & \hbox{\ \ } & \hbox{\ \ } & \hbox{\ \ } &
\hbox{\ \ } & \hbox{ \ \ } & \hbox{ \ \ } & \hbox{\ \ } & \hbox{\ \ } &
\vdots   \cr
\hbox{\ \ } & \times & \cdots & \hbox{\ \ } & \hbox{\ \ } & \hbox{\ \ } &
\hbox{\ \
} & \hbox{\ \ } & \hbox{\ \ } & \hbox{\ \ } & \lambda_t \hbox{\ nodes\ }   \cr
\endmatrix
\right.
$$

\vskip 13pt

2.2.  The main aim of the rest of this section is to discuss the conditions
under which two elements of $M_d$ (respectively, of $P_d$) are neighbours (see
Apendix A).

For any ordered pair $(l,m)$ of elements of $N_d$ we define a sequence of
integers
$$
r_k=r_k(l,m)=
\sum_{i=1}^k(m_i-l_i),\hbox{ \ } k=1,\hdots ,d-1,
$$
and set $r=r(l,m)=\sum_{k=1}^{d-1}r_k$.
It is evident that $l\leq m$ if and only if $r_k(l,m)\geq 0$ for all
$k=1,\hdots ,d-1$.

We borrow part (i) of the next lemma from [4, Ch. I, Sec. 1], and modify it in
part (ii).

\proclaim{Lemma 2.2.1} (i)
If $l$, $m\in N_d$,
then one has $l\leq m$ if and
only if there exists a raising operator $\rho$ with $m=\rho(l)$;

(ii) if $l,m\in M_d$, and if $l\leq m$, then
$m=\rho(l)$ for
a raising
operator $\rho$ having the following property: There exists a sequence of
$r=r(l,m)$ non-trivial simple raising operators $\rho_1,\rho_2,\hdots, \rho_r$
of the type $\rho_{i,i+1}$, $1\leq i\leq d-1$, such that:

(a) $\rho=\rho_r\cdots\rho_2\rho_1$;

(b) $\rho_1(l)\in M_d$,
$\rho_2\rho_1(l)\in M_d,\hdots,$
$\rho_{r-1}\hdots\rho_2\rho_1(l)\in M_d$;

(c) $l<\rho_1(l)<
\rho_2\rho_1(l)<\hdots <
\rho_{r-1}\hdots\rho_2\rho_1(l)<m$.

\endproclaim

\demo{Proof} (i) Suppose that
there exists a raising operator $\rho$ with $m=\rho(l)$. We can
assume that $\rho=\rho_{i,j}$, and in this case apparently
$l\leq m$. Conversely, let $l\leq m$. Then $m=\rho(l)$ for
$$
\rho=\prod_{k=1}^{d-1}\rho_{k,k+1}^{r_k},
$$
where
$r_k=r_k(l,m)$;

(ii) We shall prove this statement by induction with respect to $r=r(l,m)\geq
0$. If $r=0$, then $l=m$ and the trivial operator $\rho$ which is a product of
$0$ in number simple raising operators $\rho_k$ works. If $r=1$, then $l<m$,
and there exists an index $k=i$ such that $r_i=1$, and $r_k=0$ for
all $k\neq i$, $k=1,\hdots, d-1$. This implies $m=\rho_{i,i+1}(l)$, and we can
set $\rho=\rho_{i,i+1}$. Suppose that part (ii) is true for all $l, m\in
M_d$ with $l\leq m$, and $r\leq k$,
and let $r=k+1$. We set $q=q(l,m)$. Thus, $q<d-1$, and $r_1=\cdots=r_q=0$, and
$r_{q+1}\geq 1$. Let $\kappa\geq 2$ be the smallest integer with
$r_{q+\kappa}=0$ (integers $\kappa$ with the property $r_{q+\kappa}=0$ exist:
For instance, $\kappa=d-q$). We have
$$
l_{q+2}+\cdots+l_{q+\kappa}\geq
l_{q+2}+\cdots+l_{q+\kappa}-
m_{q+2}+\cdots+m_{q+\kappa}=r_{q+1}\geq 1,
$$
and hence there exists an index $j$, $q+2\leq j\leq q+\kappa$, with $l_j\geq
1$. We set $i=q+1$, and $l^\prime=\rho_{i,j}(l)$. Then $\l^\prime\in M_d$,
$l<l^\prime$, and we have $r_k(l^\prime,m)=r_k-1$ when $k=i,\hdots, j-1$, and
$r_k(l^\prime,m)=r_k$ otherwise.
Since $r_k\geq 1$ for all $k=i,\hdots
,q+\kappa-1$, then
$r_k(l^\prime,m)\geq 0$ for all $k$, so $l^\prime\leq m$. Moreover,
$r(l^\prime,m)=r(l,m)-(j-i)=k+1-(j-i)\leq k$, and the inductive assumption
yields that there exist $r^\prime=r(l^\prime,m)$ simple raising operators
$\rho_1^\prime,\hdots, \rho_{r^\prime}^\prime$ of the
desired type, such that
$m=\rho^\prime(l^\prime)$ for
$\rho^\prime=\rho_{r^\prime}^\prime\cdots\rho_1^\prime$, and conditions (b) and
(c) are satisfied. Taking into account that $\rho_{i,j}=\rho_{i,i+1}\cdots
\rho_{j-1,j}$, and that $r=r^\prime+(j-i)$, we get our statement.

\enddemo

\proclaim{Theorem 2.2.2}
The $d$-tuples $l$, $m$ are neighbours in $M_d$ with $l< m$ if and only
if there exists $i\in [1,d]$ such that $m=\rho_{i,i+1}(l)$.

\endproclaim

\demo{Proof}
Let $m=\rho_{i,i+1}(l)$ and $l\leq n\leq m$. We have
$n_k=l_k=m_k$ for $1\leq k\leq i-1$. Then $l_i\leq n_i\leq m_i=l_i+1$, so
either $n_i=l_i$, or $n_i=m_i$. Further,
$l_i+l_{i+1}\leq n_i+n_{i+1}\leq
l_i+1+l_{i+1}-1$, hence $l_i+l_{i+1}=n_i+n_{i+1}=m_i+m_{i+1}$. The two cases
$n_i=l_i$, or
$n_i=m_i$, imply $n=l$, or $n=m$, respectively.  Therefore $l$ and $m$ are
neighbours with $l< m$.

Now, suppose that the $d$-tuples $l$, $m$ are neighbours in $M_d$ with $l< m$.
According to Lemma 2.2.1, (ii), we have $m=\rho(l)$, where the raising operator
$\rho$ satisfies all conditions (a) - (c). This yields $r=1$, and hence
$m=\rho_{i,i+1}(l)$.

\enddemo

\vskip 13pt

2.3. Here we state
[2, Ch. 1, Theorem 1.4.10] which
gives necessary and sufficient conditions for
two partitions $\lambda$, $\mu\in\ P_d$ to be neighbours in $P_d$, and refer
to the corresponding proof there. It reads as follows:

\proclaim{Theorem 2.3.1}
The partitions $\lambda$, $\mu$ are neighbours in $P_d$ with $\lambda<\mu$
if and only if there exist a pair of integers $(i,j)$ with $1\leq i<j\leq d$,
and such that the following two conditions hold:

(i) One has $\mu=\rho_{i,j}(\lambda)$;

(ii) one has $j=i+1$, or $\lambda_i=\lambda_j$.

\endproclaim

In terms of Young diagrams we move the node from the end of $j$-th row of
$\lambda$ to the end of its $i$-th row and this move is minimal with the
property that we do not leave the subset $P_d\subset \Delta_d$. The last
minimum property is equivalent to (ii).

\heading{3. Dominance among ordered dissections and tabloids}

\endheading

3.1.
By an {\it ordered dissection} of the integer-valued interval
$[1,d]=\{1,2,\hdots ,d\}$
we mean a
$d$-tuple $A=(A_1,\hdots ,A_d)$
of disjoint subsets $A_i\subset [1,d]$
with $\cup_{i=1}^dA_i=[1,d]$.  Sometimes, we shall think of an ordered
dissection $A$ as an infinite sequence
$(A_1,\hdots ,A_d, A_{d+1},\hdots)$, where $A_k=\emptyset$ for $k>d$.
We denote by $\Delta_d$ the
set of all ordered dissections of
$[1,d]$,
and define
the surjective map
$$
\varphi\colon \Delta_d\to M_d,\tag 3.1.1
$$
$$
(A_1,\hdots ,A_d)\to (|A_1|,\hdots ,|A_d|).
$$

\vskip 13pt

3.2. Each ordered dissection
$A=(A_1,\hdots ,A_d)$ of $[1,d]$ with $|A_1|\geq\hdots\geq |A_d|$ is called
{\it tabloid}.
Let $T_d$ be the subset of $\Delta_d$ consisting of all tabloids.
Obviously, $T_d=\varphi^{-1}(P_d)$.
The tabloid $A$ can be visualized by
placing the elements of $A_k$ in the $k$-th row of the  Young diagram (2.1.1)
corresponding to the partition $\lambda=\varphi(A)$ without taking into
account their
order, for $k=1,\hdots,t$. The next figure illustrates both the tabloid $A$
and the map $\varphi$:

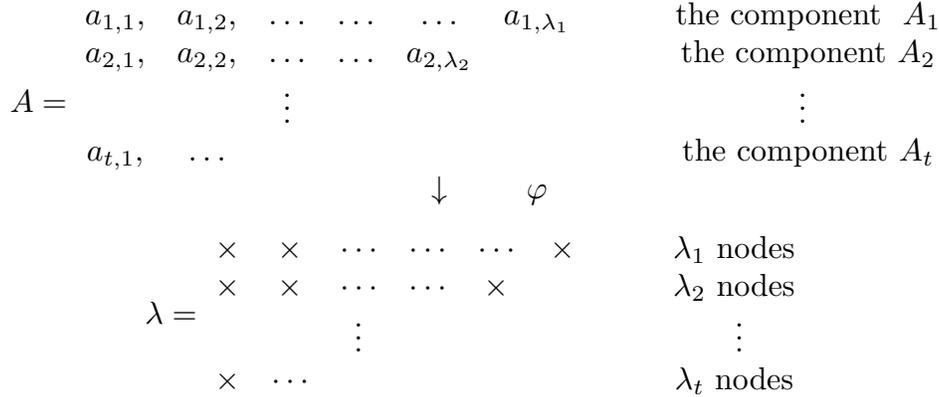
$$
A=
\left.
\matrix \format \c & \quad \c & \quad \c & \quad \c & \quad \c & \quad \c &
\quad \c & \quad \c \\
a_{1,1}, & a_{1,2}, & \hdots & \hdots & \hdots & a_{1,\lambda_1} &\hbox{\ \ } &
\hbox{\ the component \ } A_1  \cr
a_{2,1}, & a_{2,2}, & \hdots & \hdots & a_{2,\lambda_2} & \hbox{ \ \ } & \hbox{
\ \ } & \hbox{\ the component\ } A_2  \cr
\hbox{\ \ } & \hbox{\ \ } & \vdots & \hbox{\ \ } & \hbox{\ \ } & \hbox{\ \ }
& \hbox{\ \ } & \vdots   \cr
a_{t,1}, & \hdots & \hbox{\ \ } & \hbox{\ \ } & \hbox{\ \ } & \hbox{\ \ }  &
\hbox{\ \ } & \hbox{\ the component\ } A_t  \cr
\hbox{\ \ } & \hbox{\ \ } & \hbox{\ \ } & \hbox{\ \ } & \downarrow & \varphi &
\hbox{\ \ } & \hbox{\ \ }  \cr
\endmatrix
\right.
$$
$$
\lambda=
\left.
\matrix \format \c & \quad \c & \quad \c & \quad \c & \quad \c & \quad \c &
\quad \c & \quad \c \\
\times & \times & \cdots & \cdots & \cdots & \times &\hbox{\ \ } &
\lambda_1\hbox{\ nodes\ }   \cr
\times & \times & \cdots & \cdots & \times & \hbox{ \ \ } & \hbox{ \ \
} & \lambda_2\hbox{\ nodes\ }   \cr
\hbox{\ \ } & \hbox{\ \ } & \vdots & \hbox{\ \ } & \hbox{\ \ } & \hbox{\ \ }
& \hbox{\ \ } & \vdots   \cr
\times & \cdots & \hbox{\ \ } & \hbox{\ \ } & \hbox{\ \ } & \hbox{\ \ }  &
\hbox{\ \ } & \lambda_t\hbox{\ nodes\ }   \cr
\endmatrix
\right.
$$

We define a partial order on $\Delta_d$ via the rule
$$
A\leq B \hbox{\sl\ if and only if\ }
\cup_{k=1}^iA_k\subset\cup_{k=1}^i B_k,
\hbox{\sl\ for any\ } 1\leq i\leq d,
$$
and call it {\it dominance order.} In case $A\leq B$ we say that $B$ {\it
dominates} $A$.

For each $s\in [1,d]$ and each $A\in \Delta_d$ there exists
a unique $j\in [1,d]$,
such that $s\in A_j$. We set $\varepsilon_A(s)=j$.
Thus, any $A\in \Delta_d$ produces a map
$\varepsilon_A\colon [1,d]\to [1,d]$.

We introduce a partial order on the set of all maps $[1,d]\to
[1,d]$ by virtue of the rule: $\alpha\leq\beta$ if and only if
$\alpha(s)\leq\beta(s)$ for all $s\in [1,d]$.

For any two integers $1\leq i,\hbox{\ } s\leq d$, we
define an operator $R_{i,s}\colon \Delta_d\to \Delta_d$
by the formulae

$$
R_{i,s}(A)=
\left\{
\matrix \format \c & \c & \c & \c \\
(A_1,\hdots ,A_i\cup\{s\},\hdots
,A_{\varepsilon_A(s)}\backslash\{s\},\hdots, A_d) & \hbox{\sl\ \ if\ \ } &
\varepsilon_A(s)>i \cr
A & \hbox{\sl\ \ if\ \ } &
\varepsilon_A(s)\leq i.\cr
\endmatrix
\right.
$$

The operators $R_{i,s}$ are said to be {\it simple raising operators}
in $\Delta_d$.
Any product
$R=R_{i_1,s_1}R_{i_2,s_2}\hdots $ of simple raising operators
is called {\it raising operator} on the set $\Delta_d$.

The action of the simple raising operator $R_{i,s}$ on the tabloid $A$ with
$B=R_{i,s}A\in T_d$ can be illustrated by the picture below:

$$
A=
\left.
\matrix \format \c & \quad \c & \quad \c & \quad \c & \quad \c & \quad \c &
\quad \c & \quad \c & \quad \c & \quad \c & \quad \c \\
a_{1,1}, & \ast & \hdots & \ast & \hdots & \ast & \ast
 & a_{1,\lambda_1} & \hbox{\ \ } & \hbox{\ the component\ } A_1  \cr
\vdots & \hbox{\ \ } & \hbox{\ \ } & \hbox{\ \ } &
\hbox{\ \ } & \hbox{ \ \ } & \hbox{ \ \ } & \hbox{\ \ } & \hbox{\ \ } &
\vdots   \cr
a_{i,1}, & \ast & \hdots & \ast & \hdots & a_{i,\lambda_i} &
\hbox{\ \ } & \hbox{\ \ } & \hbox{\ \ } & \hbox{\ the component\ } A_i  \cr
\vdots & \hbox{\ \ } & \hbox{\ \ } & \hbox{\ \ } &
\hbox{\ \ } & \nearrow & \hbox{ \ \ } & \hbox{\ \ } & \hbox{\ \ } &
\vdots   \cr
a_{j,1}, & \ast & \hdots  & \ast & s & \hbox{\ \ } &
\hbox{\ \ } & \hbox{\ \ } & \hbox{\ \ } & \hbox{\ the component\ } A_j   \cr
\vdots & \hbox{\ \ } & \hbox{\ \ } & \hbox{\ \ } &
\hbox{\ \ } & \hbox{ \ \ } & \hbox{ \ \ } & \hbox{\ \ } & \hbox{\ \ } &
\vdots   \cr
a_{t,1}, & \hdots & \hbox{\ \ } & \hbox{\ \ } & \hbox{\ \ } &
\hbox{\ \
} & \hbox{\ \ } & \hbox{\ \ } & \hbox{\ \ } & \hbox{\ the component\ } A_t  \cr
& \hbox{\ \ } & \hbox{\ \ } & \hbox{\ \ } & \hbox{\ \ } &
\hbox{\ \ } & \downarrow & R_{i,s} & \hbox{\ \ } & \hbox{\ \ } &
\hbox{\ \ } \cr
\endmatrix
\right.
$$
$$
B=
\left.
\matrix \format \c & \quad \c & \quad \c & \quad \c & \quad \c & \quad \c &
\quad \c & \quad \c & \quad \c & \quad \c & \quad \c \\
a_{1,1}, & \ast & \hdots & \ast & \hdots & \ast & \ast
 & a_{1,\lambda_1} & \hbox{\ \ } & \hbox{\ the component\ } B_1  \cr
\vdots & \hbox{\ \ } & \hbox{\ \ } & \hbox{\ \ } &
\hbox{\ \ } & \hbox{ \ \ } & \hbox{ \ \ } & \hbox{\ \ } & \hbox{\ \ } &
\vdots   \cr
a_{i,1}, & \ast & \hdots & \ast & \hdots  & a_{i,\lambda_i},  & s
& \hbox{\ \ } & \hbox{\ \ } & \hbox{\ the component\ } B_i  \cr
\vdots & \hbox{\ \ } & \hbox{\ \ } & \hbox{\ \ } &
\hbox{\ \ } & \hbox{\ \ } & \hbox{ \ \ } & \hbox{\ \ } & \hbox{\ \ } &
\vdots   \cr
a_{j,1}, & \ast & \hdots  & \ast & \hbox{\ \ } & \hbox{\ \ }
& \hbox{\ \ } & \hbox{\ \ } & \hbox{\ \ } & \hbox{\ the component \ } B_j   \cr
\vdots & \hbox{\ \ } & \hbox{\ \ } & \hbox{\ \ } &
\hbox{\ \ } & \hbox{ \ \ } & \hbox{ \ \ } & \hbox{\ \ } & \hbox{\ \ } &
\vdots   \cr
a_{t,1}, & \hdots & \hbox{\ \ } & \hbox{\ \ } & \hbox{\ \ } &
\hbox{\ \
} & \hbox{\ \ } & \hbox{\ \ } & \hbox{\ \ } & \hbox{\ the component\ } B_t  \cr
\endmatrix
\right.
$$

It is easy to see that any two simple raising operators
commute.
Thus, for any
$i\in
[1,d]$,
and for any subset $X\subset [1,d]$ we can define without ambiguity
$R_{i,X}=\prod_{x\in X}R_{i,x}$.

For any
$i\in
[1,d]$, and for any finite family $J=(j_x)_{x\in X}$ of elements of  $[1,d]$,
we define a raising operator in $N_d$ by
$\rho_{i,J}=\prod_{x\in X}\rho_{i,j_x}$.

\proclaim{Lemma 3.2.1}
(i) For any $A\in \Delta_d$ and any raising operator
$R=R_{i_1,s_1}R_{i_2,s_2}\hdots$,
one has the inequality
$\varepsilon_{R\left(A\right)}\leq\varepsilon_A$.
If there exists a pair $i_k,s_k$ with $\varepsilon_A(s_k)>i_k$, then
$\varepsilon_{R\left(A\right)}<\varepsilon_A$;

(ii)
for any subset $X\subset [1,d]$,
one has
$\varphi(R_{i,X}A)=\rho_{i,\varepsilon_A\left(X\right)}\varphi(A)$;

(iii) the map $\varphi\colon \Delta_d\to M_d$ is a homomorphism of partially
ordered sets: $\varphi(A)\leq\varphi(B)$ for $A\leq B$;
if $A\leq B$ and $\varphi(A)=\varphi(B)$, then
$A=B$.

\endproclaim

\demo{Proof}
(i) It is enough to prove the first statement for $R=R_{i,s}$. When
$\varepsilon_A(s)\leq i$, it is obvious. Now, let
$\varepsilon_A(s)>i$; Since $i=\varepsilon_{R\left(A\right)}(s)$ and since
$\varepsilon_A(t)=\varepsilon_{R\left(A\right)}(t)$ for $t\neq s$, then
$\varepsilon_{R\left(A\right)}<\varepsilon_A$, and we have proved both the
first statement and the second statement for $R=R_{i,s}$.

For the second statement, we write
$R=R^\prime R_{i_k,s_k}$.
Then $R(A)=
R^\prime R_{i_k,s_k}(A)$ and
$\varepsilon_{R\left(A\right)}\leq
\varepsilon_{R_{i_k,s_k}\left(A\right)}<\varepsilon_A$.

(ii) We shall use induction with respect to the number of elements in the
set $X$. When $|X|=1$, this is trivial. Suppose $|X|\geq 2$, and set
$X^\prime=X\backslash\{s\}$, where $s\in X$,
$B=R_{i,X^\prime}A$, and
$j=\varepsilon_B(s)$.
We have
$$
\varphi(R_{i,X}A)=
\varphi(R_{i,s}R_{i,X^\prime}A)=
\varphi(R_{i,s}B)=
\rho_{i,j}\varphi(B)=
$$
$$
\rho_{i,j}
\varphi(R_{i,X^\prime}A)=
\rho_{i,j}
\rho_{i,\varepsilon_A\left(X^\prime\right)}\varphi(A).
$$
Since $s\notin X^\prime$, then
$j=\varepsilon_B(s)=
\varepsilon_A(s)$, so part (ii) is proved.

(iii) This is a direct consequence of the definitions of the partial orders
on $\Delta_d$ and $M_d$.

\enddemo

\proclaim{Lemma 3.2.2} Let $A\in \Delta_d$.
If
$R=R_{i_1,s_1}R_{i_2,s_2}\hdots$
is a raising operator,
then
$A\leq R(A)$. In particular, if
there exists a pair $i_k,s_k$ with $\varepsilon_A(s_k)>i_k$,
then $A<R(A)$.

\endproclaim

\demo{Proof} Let $B=R(A)$. We can suppose that $R=R_{i,s}$ and in this case
the inequality $A\leq R(A)$ is obvious.
Now, Lemma 3.2.1, (i), yields the statement.

\enddemo

\vskip 13pt

3.3.
Let $A, B\in \Delta_d$ with $A\leq B$, and let $l=\varphi(A)$ and
$m=\varphi(B)$.
According to Lemma 3.2.1, (iii), the map $\varphi$, defined via (3.1.1),
is a homomorphism of partially ordered sets. In particular, $\varphi$ maps the
interval $[A,B]$ into the interval $[l,m]$. In the next two lemmas we begin
the study of the equation $\varphi(X)=n$, where $X\in [A,B]$, for various $n\in
[l,m]$.

\proclaim{Lemma 3.3.1}
Let $A, B\in \Delta_d$ with $A\leq B$, and let $l=\varphi(A)$ and
$m=\varphi(B)$.
Suppose $l\leq n\leq m$, where $n\in M_d$.
If for some
$i$, $1\leq i\leq d$, one has
$$
i-1=q(l,n),
$$
then there exists a raising operator $R_{i,X}$ with $X\subset
A_{i+1}\cup\hdots\cup A_d$, such that $A^\prime=R_{i,X}(A)$ and
$l^\prime =\varphi(A^\prime)$
satisfy the conditions
$A<A^\prime\leq B$,
and
$l<l^\prime\leq n$, and
$$
i\leq q(l^\prime,n).
$$

\demo{Proof} If $i=d$, then $n=l$, and we choose $X$ to be the empty set.
Now, let $i<d$. The equality
$i-1=q(l,n)$ implies $l_1=n_1,\hdots, l_{i-1}=n_{i-1}$ and $l_i<n_i$. Hence,
$$
l_1+\cdots +l_i<
n_1+\cdots +n_i\leq
m_1+\cdots +m_i.
$$
We choose a subset
$X\subset B_1\cup\hdots\cup B_i\backslash A_1\cup\hdots\cup A_i$
consisting of $n_i-l_i$ elements. Obviously, $X\subset A_{i+1}\cup\hdots\cup
A_d$. We set
$A^\prime=R_{i,X}(A)$. Then $l^\prime
=\rho_{i,\varepsilon_A\left(X\right)}(l)$, and the conditions of the lemma are
satisfied.

\enddemo

\proclaim{Lemma 3.3.2}
Let $A, B\in \Delta_d$ with $A<B$, and let $l=\varphi(A)\in M_d$ and
$m=\varphi(B)$. Suppose
that $m=\rho_{i,j}l$, where $1\leq i<j\leq d$, and that
there exist an integer $r\geq 1$, and two
sequences $(i_k)_{k=1}^r$
and $(s_\kappa)_{\kappa=1}^r$ in the interval $[1,d]$, such that
$$
i=i_1<i_2<\hdots <i_r<j, \hbox{\sl\ and\ }
\varepsilon_A(s_\kappa)=i_{\kappa+1},\hbox{\sl\ for all\ }1\leq \kappa\leq r-1,
\hbox{\sl\ and\ } \varepsilon_A(s_r)>i_r,
$$
and that the components of the ordered dissections
$R_{i_1,s_1}\hdots R_{i_r,s_r}A$ and $B$ coincide for all indices in the closed
interval $[1,i_r]$. Then there exist two integers $i_{r+1}$, $s_{r+1}$ in
$[1,d]$, such that $i_r<i_{r+1}\leq j$, and $\varepsilon_A(s_{r+1})=i_{r+1}$,
and in case $i_{r+1}<j$ the components of the ordered dissections
$R_{i_1,s_1}\hdots R_{i_{r+1},s_{r+1}}A$ and $B$ coincide for all
indices in the closed interval $[1,i_{r+1}]$, or one has
$B=R_{i_1,s_1}\hdots R_{i_r,s_r}A$
in case $i_{r+1}=j$.

\endproclaim

\demo{Proof} It is obvious that the elements $s_1,\hdots ,s_r\in [1,d]$ are
pairwise different.
The condition yields
$B_{i_k}=(A_{i_k}\backslash\{s_{k-1}\})\cup\{s_k\}$ for all $2\leq k\leq r$,
and
$B_{i_1}=A_{i_1}\cup\{s_1\}$,
and
$B_k=A_k$  for all $1\leq k\leq i_r$ with
$k\notin\{i_1,\hdots ,i_r\}$.

We shall prove the following

\proclaim{Sublemma}
(i) One has $A_k=B_k$ for all $k\in
[i_r+1,\min\{\varepsilon_A(s_r),j\}-1]$;

(ii) one has $\varepsilon_A(s_r)\leq j$.

\endproclaim

\demo{Proof}
When
$\min\{\varepsilon_A(s_r),j\}=i_r+1$,
that is, the interval
$[i_r+1,\min\{\varepsilon_A(s_r),j\}-1]$
is empty, the
statement is trivial. Let $\min\{\varepsilon_A(s_r),j\}>i_r+1$. We have
$$
A_{i_1}\cup\hdots\cup A_{i_r}\cup A_{i_r+1}\subset
B_{i_1}\cup\hdots\cup B_{i_r}\cup B_{i_r+1}.
$$
Since
$$
B_{i_1}\cup\hdots\cup B_{i_r}=
A_{i_1}\cup\hdots\cup A_{i_r}\cup\{s_r\},
$$
and since $s_r\notin A_{i_r+1}$, we obtain
$A_{i_r+1}\subset B_{i_r+1}$.
Then
$l_{i_r+1}=m_{i_r+1}$ implies
$A_{i_r+1}=B_{i_r+1}$.
Suppose that
$$
A_{i_r+1}=B_{i_r+1},\hdots, A_{k-1}=B_{k-1},
$$
for $i_r+1<k\leq \min\{\varepsilon_A(s_r),j\}-1$. Then we get
$A_k\subset B_k\cup\{s_r\}$,
and because of $s_r\notin A_k$, we obtain
$A_k\subset B_k$.
Then
$l_k=m_k$ implies
$A_k=B_k$. Thus, part (i) is proved by induction.

(ii) Suppose the opposite, that is, $\varepsilon_A(s_r)>j$.
Then, according to part (i), we have
$A_j\subset B_j\cup\{s_r\}$.
Again $s_r\notin A_j$ yields
$A_j\subset B_j$. On the other hand,
$l_j-1=m_j$,
which is a contradiction.

\enddemo

We set $i_{r+1}=\varepsilon_A(s_r)$. According to the above Sublemma,
$i_r<i_{r+1}\leq j$ and
$A_k=B_k$ for all $i_r<k<i_{r+1}$. Thus, we have
$A_{i_{r+1}}\subset B_{i_{r+1}}\cup\{s_r\}$,
so
$A_{i_{r+1}}\backslash\{s_r\}\subset
B_{i_{r+1}}$.

Case 1. $i_{r+1}<j$.

Since $l_{i_{r+1}}=m_{i_{r+1}}$,
there exists an element $s_{r+1}\in
B_{i_{r+1}}$ such that
$s_{r+1}\notin A_{i_{r+1}}\backslash\{s_r\}$, and
$B_{i_{r+1}}=(A_{i_{r+1}}\backslash\{s_r\})\cup\{s_{r+1}\}$.
Since $s_k\in B_{i_k}$,
we have $s_{r+1}\neq s_k$
for $1\leq k\leq r$.
This
implies $s_{r+1}\notin A_{i_{r+1}}$; hence $\varepsilon_A(s_{r+1})>i_{r+1}$.
Having this information, it is not hard to check that
the components of the ordered dissections
$B$ and $R_{i_1,s_1}\hdots R_{i_r,s_r}R_{i_{r+1},s_{r+1}}A$ coincide for all
indices in the closed interval $[1,i_{r+1}]$.

Case 2. $i_{r+1}=j$.

Since $l_j-1=m_j$, then
$A_j\backslash\{s_r\}=B_j$,
so
the components of the ordered dissections
$R_{i_1,s_1}\hdots R_{i_r,s_r}A$ and $B$ coincide for all
indices in the closed interval $[1,j]$.
Now, we shall prove that $B_k=A_k$ for all $j+1\leq k\leq d$. We have
$\cup_{k=1}^j A_k=
\cup_{k=1}^j B_k$,
so
$A_{j+1}\subset B_{j+1}$.
Therefore the equality
$l_{j+1}=m_{j+1}$ gives
$A_{j+1}=B_{j+1}$.
Obvious induction finishes the proof.

\enddemo

\vskip 13pt

3.4. We say that $l\in M_d$ and $m\in M_d$ are {\it adjacent}
with $l<m$ if $m=\rho_{i,j}l$ for some pair of integers $(i,j)$ with $1\leq
i<j\leq
d$. Given $A$, $B\in \Delta_d$, we set $l=\varphi(A)$ and $m=\varphi(B)$.
The ordered dissections $A$ and $B$ are called {\it adjacent}
with $A<B$ if $A<B$, and $l$ and $m$ are adjacent (with $l<m$).
The ordered dissections $A$ and $B$ are said to be {\it strongly
adjacent} with $A<B$
if $B=R_{i,s}l$ for some pair $i, s\in [1,d]$ such that
$\varepsilon_A(s)>i$.
According to Lemma 3.2.1, (ii), if $A$ and $B$ are strongly adjacent in
$\Delta_d$ with $A<B$, then $A$ and $B$ are adjacent.  The converse statement
is not true.  The situation is clarified in the next theorem.

\proclaim{Theorem 3.4.1}
Let $A, B\in \Delta_d$ be adjacent with $A<B$,
and let $l=\varphi(A)\in M_d$ and
$m=\varphi(B)$. Suppose
that $m=\rho_{i,j}l$, where $1\leq i<j\leq d$. Then
there exist an integer $r\geq 1$, and two
sequences $(i_k)_{k=1}^{r+1}$
and $(s_\kappa)_{\kappa=1}^r$ in the interval $[1,d]$, such that
$$
i=i_1<i_2<\hdots <i_{r+1}=j, \hbox{\sl\ and\ }
\varepsilon_A(s_\kappa)=i_{\kappa+1}\hbox{\sl\ for all\ }1\leq \kappa\leq r,
$$
and that
$B=R_{i_1,s_1}\hdots R_{i_r,s_r}A$.

\endproclaim

\demo{Proof} We apply several times Lemma 3.3.2. In order to begin, we note
that $q(l,m)=i-1$, and use Lemma 3.3.1 in case $n=m$, thereby producing the
first pair
$(i_1,s_1)$ with $i_1=i$, and $\varepsilon_A(s_1)>i$. It is obvious that the
components
of the ordered dissections $B$ and $R_{i_1,s_1}A$ coincide for all indices in
the interval $[1,i_1]$.

\enddemo

\proclaim{Theorem 3.4.2}
Let $A, B\in \Delta_d$ with $A\leq B$, and let $l=\varphi(A)\in M_d$ and
$m=\varphi(B)\in M_d$. For any $n\in M_d$ with $l\leq n\leq m$, and
$q(l,n)=q$, there exists a raising operator of the
type $R=R_{d,X_d}\hdots R_{q+1,X_{q+1}}$ with
$X_k\subset A_{k+1}\cup\hdots\cup A_d$,
such that $A^\prime=R(A)$
satisfies the conditions
$A\leq A^\prime\leq B$,
and
$\varphi(A^\prime)=n$.

\endproclaim

\demo{Proof} We shall use induction with respect to $q=q(l,n)$.
If $q=d$, then $l=n$ and the ordered dissection
$A^\prime=R(A)=A$ for the trivial operator $R=R_{d,X_d}R_{d+1,X_{d+1}}$,
$X_d=X_{d+1}=\emptyset$, works. Suppose that if $i\leq q\leq d$, then
there exists a raising operator of the type
$R=R_{d,X_d}\hdots R_{q+1,X_{q+1}}$,
such that $A^\prime=R(A)$
satisfies the conditions
$A\leq A^\prime\leq B$,
and
$\varphi(A^\prime)=n$.
If $q=i-1$, then Lemma 3.3.1 yields the
existence of a raising operator of the type $R^{\prime\prime}=R_{i,X_i}$ with
$X_i\subset A_{i+1}\cup\hdots\cup A_d$, such that
$A^{\prime\prime} =R^{\prime\prime}(A)$, and
$l^{\prime\prime} =\rho_{i,\varepsilon_A\left(X_i\right)}(l)$,
satisfy the conditions
$l<l^{\prime\prime}\leq n$, and
$i\leq q(l^{\prime\prime},n)\leq d$,
and
$A<A^{\prime\prime}\leq B$, and $\varphi(A^{\prime\prime})=l^{\prime\prime}$.
Hence,
there exists a raising operator
$$
R^\prime=R_{d,X_d}\hdots R_{i+1,X_{i+1}},
$$
such that $A^\prime=R^\prime(A^{\prime\prime})$
satisfies the conditions
$A^{\prime\prime}<A^\prime\leq B$,
and
$\varphi(A^\prime)=n$.
Since
$A^\prime=R(A)$ for
$$
R=R^\prime R^{\prime\prime}=
R_{d,X_d}\hdots R_{i+1,X_{i+1}}R_{i,X_i}=
R_{d,X_d}\hdots R_{q+1,X_{q+1}},
$$
the induction is done.

\enddemo

\proclaim{Theorem 3.4.3}
(i) Let $A, B\in \Delta_d$ with $A<B$. Let
$l=\varphi(A)$ and $m=\varphi(B)$. Then
the
restriction $\varphi_1$ of the map $\varphi$ on the interval $[A,B]$
in $\Delta_d$
is a surjection
$$
\varphi_1\colon [A,B]\to [l,m],
$$
and
one has $\varphi_1^{-1}((l,m))=(A,B)$;

(ii) let $A, B\in T_d$ with $A\leq B$. Let
$\lambda=\varphi(A)$ and $\mu=\varphi(B)$. Then
the
restriction $\varphi_2$ of the map $\varphi$ on the interval $[A,B]$ in $T_d$
is a surjection
$$
\varphi_2\colon [A,B]\to [\lambda,\mu],
$$
and
one has $\varphi_2^{-1}((\lambda,\mu))=(A,B)$;

\endproclaim

\demo{Proof} (i) The surjectivity of $\varphi_1$ is a consequence of
Theorem 3.4.2.
The inclusion
$\varphi_1^{-1}((l,m))\subset (A,B)$ is obvious. Suppose that $A<C<B$.
Then
the assumption that $\varphi(C)=l$, or $\varphi(C)=m$ leads to a contradiction
with Lemma 3.2.1, (iii).

(ii)
If $C\in M_d$ with
$\varphi(C)\in P_d$, then $C\in T_d$, so part (i) assures
that the map
$\varphi_2$ is surjective. The rest of the proof is identical to
that of part (i).

\enddemo

\proclaim{Theorem 3.4.4} If $A, B\in \Delta_d$
then
$A\leq B$ if and only if there exists a raising operator $R$ such that
$B=R(A)$.

\endproclaim

\demo{Proof} The \lq\lq if" part follows from Lemma 3.2.2. Now, let
$A, B\in \Delta_d$, $A\leq B$, with
$l=\varphi(A)$ and $m=\varphi(B)$. In case $A=B$ we choose $R$ to be the
trivial operator. Now, let $A<B$. We apply Theorem 3.4.2 in the particular
case $n=m$ to produce a raising operator $R$ such that the ordered dissection
$A^\prime=R(A)$ satisfies $A\leq A^\prime\leq B$,
and
$\varphi(A^\prime)=m=\varphi(B)$. Then Lemma 3.2.1, (iii), yields that
$B=A^\prime=R(A)$.

\enddemo

\vskip 13pt

3.5. Here we find necessary and sufficient conditions for two ordered
dissections, or for two tabloids to be neighbours with respect to the partial
orders on $\Delta_d$ and
on $T_d$, respectively (see Appendix A).

\proclaim{Theorem 3.5.1}
(i) The ordered dissections $A, B\in \Delta_d$ are neighbours in $\Delta_d$ with
$A<B$, if and only if
there exist $i\in [1,d]$ and $s\in [1,d]$, such that
$\varepsilon_A(s)=i+1$ and
$B=R_{i,s}(A)$;

(ii)
the tabloids  $A, B\in T_d$ are neighbours in $T_d$ with
$A<B$,
if and only if
there exist a pair of integers $(i,j)$ with $1\leq i<j\leq d$,
an integer $r\geq 1$,
and two sequences $(i_k)_{k=1}^{r+1}$
and $(s_\kappa)_{\kappa=1}^r$ in the interval $[1,d]$, such that:
$$
j=i+1 \hbox{\ or\ } |A_i|=|A_j|,\tag 3.5.2
$$
and
$$
i=i_1<i_2<\hdots <i_{r+1}=j, \hbox{\sl\ and\ }
\varepsilon_A(s_\kappa)=
i_{\kappa+1},\hbox{\sl\ for all\ }1\leq \kappa\leq r,\tag 3.5.3
$$
and that
$$
B=R_{i_1,s_1}\hdots R_{i_r,s_r}A.\tag 3.5.4
$$

\endproclaim

\demo{Proof}
(i)
We set $l=\varphi(A)$, and $m=\varphi(B)$.
Suppose that the pair $A, B\in \Delta_d$ is such that  $B=R_{i,s}(A)$ with
$\varepsilon_A(s)=i+1$.
Then
Lemma 3.2.1, (ii), yields $m=\rho_{i,i+1}l$. Hence, according to Theorem 2.2.2
we have that $l$ and $m$ are neighbours with $l<m$, and now
Theorem 3.4.3, (i), yields that $A$ and $B$ are neighbours  in
$\Delta_d$ with $A<B$.

Assume that $A, B\in \Delta_d$ are neighbours in $\Delta_d$ with
$A<B$. Theorem 3.4.3, (i), implies that
$l$ and $m$ are neighbours in $M_d$ with $l<m$.
Then, due to Theorem 2.2.2 there exist an integer $1\leq i<d$, such
that $m=\rho_{i,i+1}l$,
and Theorem 3.4.1 yields the existence of an element $s\in [1,d]$ with
$\varepsilon_A(s)=i+1$ and $B=R_{i,s}A$.

(ii) Suppose
that $A, B\in T_d$ are neighbours in $T_d$ with
$A<B$.
Denote $\lambda=\varphi(A)$ and $\mu=\varphi(B)$.
Theorem 3.4.3, (ii), implies that the partitions $\lambda$ and $\mu$ are
neighbours in $P_d$ with
$\lambda <\mu$. Due to Theorem 2.3.1, there is a pair of
integers $(i,j)$ with $1\leq i<j\leq d$, and
such that $\mu=\rho_{i,j}\lambda$. Therefore, according to Theorem
3.4.1, there exist an integer $r\geq 1$, and two
sequences $(i_k)_{k=1}^{r+1}$
and $(s_\kappa)_{\kappa=1}^r$ in the interval $[1,d]$, such that (3.5.3) and
(3.5.4) hold.
Moreover,
Theorem 2.3.1 yields
(3.5.2).

Conversely, suppose that the conditions (3.5.2) -- (3.5.4) are satisfied.
Applying
the map $\varphi$ on the equality (3.5.4),
we obtain
$$
\mu=\varphi(B)=\rho_{i_1,i_2}\hdots\rho_{i_r,i_{r+1}}\varphi(A)=\rho_{i,j}
\lambda.
$$
Therefore Theorem 2.3.1 assures that
the partitions $\lambda$ and $\mu$ are
neighbours in $P_d$ with
$\lambda <\mu$. Now, according to
Theorem 3.4.3, (ii), the tabloids
$A$ and $B$ are neighbours in
$T_d$ with $A<B$.

\enddemo

The next picture illustrates Theorem 3.5.1, (ii), case $j>i+1$, when  there
exists a sequence of \lq\lq virtual substitutions"
which starting with $A$ produces $B$. Here \lq\lq virtual" means that during
the intermediate steps we leave the set $T_d$ of tabloids.

$$
A=
\left.
\matrix \format \c & \quad \c & \quad \c & \quad \c & \quad \c & \quad \c &
\quad \c & \quad \c & \quad \c & \quad \c & \quad \c & \quad \c & \quad \c \\
\hbox{\ \ } & a_{1,1}, & a_{1,2}, & \ast & \ast & \hdots & \ast
& \ast & \ast & a_{1,\lambda_1} & \hbox{\ \ }
A_1  \cr
\hbox{\ \ } & \vdots & \hbox{\ \ } & \hbox{\ \ } & \hbox{\ \ } & \hbox{\ \ } &
\hbox{ \ \ } & \hbox{ \ \ } & \hbox{\ \ } & \hbox{\ \ } &
\vdots   \cr
i=i_1 & a_{i,1}, & a_{i,2}, & \ast & \ast & \hdots & \ast &
a_{i,\lambda_i} & \hbox{\ \ } & \hbox{\ \ } & \hbox{\ \
} A_i  \cr
\hbox{\ \ } & \vdots & \hbox{\ \ } & \hbox{\ \ } & \hbox{\ \ } & \hbox{\ \ } &
 \hbox{ \ \ } & \hbox{ \ \ } & \hbox{\ \ } & \hbox{\ \ } &
\vdots   \cr
i_r & \ast & \ast & \ast & \ast & \hdots & \ast &
\ast & \hbox{\ \ } & \hbox{\ \ } & \hbox{\ \
} A_{i_r}  \cr
\hbox{\ \ } & \vdots & \hbox{\ \ } & \nwarrow & \hbox{\ \ } & \hbox{\ \
} & \hbox{ \ \ } & \hbox{ \ \ } & \hbox{\ \ } & \hbox{\ \ } &
\hbox{\ \ }   \cr
j=i_{r+1} & a_{j,1}, & \ast & \ast & s_r, & \hdots & \ast &
a_{j,\lambda_j} & \hbox{\ \ }& \hbox{\ \ } &
 \hbox{\ \ } A_j   \cr
\hbox{\ \ } & \vdots & \hbox{\ \ } & \hbox{\ \ } & \hbox{\ \ } &
\hbox{\ \ } & \hbox{ \ \ } & \hbox{ \ \ } & \hbox{\ \ } & \hbox{\
\ } & \vdots   \cr
\hbox{\ \ } & a_{t,1}, & \hdots & \hbox{\ \ } & \hbox{\ \ } & \hbox{\ \ } &
\hbox{\ \
} & \hbox{\ \ } & \hbox{\ \ } & \hbox{\ \ } & \hbox{\ \ } A_t  \cr
\hbox{\ \ } & \hbox{\ \ } & \hbox{\ \ } & \hbox{\ \ } & \hbox{\ \ } &
\hbox{\ \ } & \downarrow R_{i_r,s_r} & \hbox{\ \ } & \hbox{\ \ }
& \hbox{\ \ } & \hbox{\ \ } \cr
\endmatrix
\right.
$$
$$
\left.
\matrix \format \c & \quad \c & \quad \c & \quad \c & \quad \c & \quad \c &
\quad \c & \quad \c & \quad \c & \quad \c & \quad \c & \quad \c & \quad \c \\
\hbox{\ \ } & a_{1,1}, & a_{1,2}, & \ast & \ast & \hdots & \ast
& \ast & \ast & a_{1,\lambda_1} & \hbox{\ \ }
A_1  \cr
\hbox{\ \ } & \vdots & \hbox{\ \ } & \hbox{\ \ } & \hbox{\ \ } & \hbox{\ \ } &
\hbox{ \ \ } & \hbox{ \ \ } & \hbox{\ \ } & \hbox{\ \ } &
\vdots   \cr
i=i_1 & a_{i,1}, & a_{i,2}, & \ast & \ast & \hdots & \ast &
a_{i,\lambda_i} & \hbox{\ \ } & \hbox{\ \ } & \hbox{\ \
} A_i  \cr
\hbox{\ \ } & \vdots & \hbox{\ \ } & \hbox{\ \ } & \hbox{\ \ } & \hbox{\ \ } &
\hbox{ \ \ } & \hbox{ \ \ } & \hbox{\ \ } & \hbox{\ \ } &
\vdots   \cr
i_{r-1} & \ast & \ast  & \ast & \ast & \hdots & \ast &
\ast & \hbox{\ \ } & \hbox{\ \ } & \hbox{\ \
} A_{i_{r-1}}  \cr
\hbox{\ \ } & \vdots & \hbox{\ \ } & \hbox{\ \ } & \nearrow &  \hbox{\ \ } &
 \hbox{ \ \ } & \hbox{ \ \ } & \hbox{\ \ } & \hbox{\ \ } &
\vdots   \cr
i_r & \ast & s_r,  & s_{r-1}, & \ast & \hdots & \ast &
\ast & \ast & \hbox{\ \ } & \hbox{\ \
} A_{i_r}\cup\{s_r\}  \cr
\hbox{\ \ } & \vdots & \hbox{\ \ } & \hbox{\ \ } & \hbox{\ \ } & \hbox{\ \
} & \hbox{ \ \ } & \hbox{ \ \ } & \hbox{\ \ } & \hbox{\ \ } &
\hbox{\ \ }   \cr
j=i_{r+1} & a_{j,1}, & \ast & \ast & \hat{s_r}, & \hdots  &
a_{j,\lambda_j} & \hbox{\ \ }& \hbox{\ \ } & \hbox{\ \ } &
 \hbox{\ \ } A_j\backslash\{s_r\}   \cr
\hbox{\ \ } & \vdots & \hbox{\ \ } & \hbox{\ \ } & \hbox{\ \ } &
\hbox{\ \ } & \hbox{ \ \ } & \hbox{ \ \ } & \hbox{\ \ } & \hbox{\
\ } & \vdots   \cr
\hbox{\ \ } & a_{t,1}, & \hdots & \hbox{\ \ } & \hbox{\ \ } & \hbox{\ \ } &
\hbox{\ \
} & \hbox{\ \ } & \hbox{\ \ } & \hbox{\ \ } & \hbox{\ \ } A_t  \cr
\hbox{\ \ } & \hbox{\ \ } & \hbox{\ \ } & \hbox{\ \ } & \hbox{\ \ } &
\hbox{\ \ } & \downarrow & R_{i_{r-1},s_{r-1}} & \hbox{\ \ }
& \hbox{\ \ } & \hbox{\ \ } \cr
\hbox{\ \ } & \hbox{\ \ } & \hbox{\ \ } & \hbox{\ \ } & \hbox{\ \ } & \hbox{\ \
} & \vdots & \hbox{ \ \ } & \hbox{\ \ } & \hbox{\ \ } &
\hbox{\ \ }   \cr
\hbox{\ \ } & \hbox{\ \ } & \hbox{\ \ } & \hbox{\ \ } & \hbox{\ \ } &
\hbox{\ \ } & \downarrow & \hbox{\ \ } & \hbox{\ \ }
& \hbox{\ \ } & \hbox{\ \ } \cr
\endmatrix
\right.
$$
$$
B=
\left.
\matrix \format \c & \quad \c & \quad \c & \quad \c & \quad \c & \quad \c &
\quad \c & \quad \c & \quad \c & \quad \c & \quad \c & \quad \c & \quad \c \\
\hbox{\ \ } & a_{1,1}, & a_{1,2}, & \ast & \ast & \hdots & \ast
& \ast & \ast & a_{1,\lambda_1} & \hbox{\ \ }
A_1  \cr
\hbox{\ \ } & \vdots & \hbox{\ \ } & \hbox{\ \ } & \hbox{\ \ } & \hbox{\ \ } &
\hbox{ \ \ } & \hbox{ \ \ } & \hbox{\ \ } & \hbox{\ \ } &
\vdots   \cr
i=i_1 & a_{i,1}, & a_{i,2}, & \ast & \ast & \hdots & \ast &
a_{i,\lambda_i}, & s_1 & \hbox{\ \ } & \hbox{\ \
} A_i\cup\{s_1\}  \cr
\hbox{\ \ } & \vdots & \hbox{\ \ } & \hbox{\ \ } & \hbox{\ \ } & \hbox{\ \ } &
 \hbox{ \ \ } & \hbox{ \ \ } & \hbox{\ \ } & \hbox{\ \ } &
\vdots   \cr
j=i_{r+1} & a_{j,1}, & \ast & \ast & \hat{s_r}, & \hdots  &
a_{j,\lambda_j} & \hbox{\ \ } & \hbox{\ \ }& \hbox{\ \ } &
 \hbox{\ \ } A_j\backslash\{s_r\}   \cr
\hbox{\ \ } & \vdots & \hbox{\ \ } & \hbox{\ \ } & \hbox{\ \ } &
\hbox{\ \ } & \hbox{ \ \ } & \hbox{ \ \ } & \hbox{\ \ } & \hbox{\
\ } & \vdots   \cr
\hbox{\ \ } & a_{t,1}, & \hdots & \hbox{\ \ } & \hbox{\ \ } & \hbox{\ \ } &
\hbox{\ \
} & \hbox{\ \ } & \hbox{\ \ } & \hbox{\ \ } & \hbox{\ \ } A_t  \cr
\endmatrix
\right.
$$
(The hat over a number stands for absence of that number.)

\heading{4. The model}

\endheading

4.1. The symmetric group $S_d$ acts on the set $\Delta_d$ of all ordered
dissections of $[1,d]$ by the rule (1.1.1).
Let $W\leq S_d$ be a subgroup of the symmetric group $S_d$. Then the group
$W$ acts on the set $\Delta_d$ via the same rule. We denote by
$\Delta_{d;W}$ the
factor-set $W\backslash \Delta_d$ and by
$T_{d;W}$ --- the factor-set $W\backslash
T_d$. Let $\psi_W\colon \Delta_d\to \Delta_{d;W}$ be the natural
surjection. For any $A\in \Delta_d$
we denote by $O_W(A)$ its $W$-orbit in $\Delta_d$, so
$\psi_W(A)=O_W(A)$. Since $\varphi(\sigma A)=\varphi(A)$
for any $A\in \Delta_d$ and for any $\sigma\in W$,
the map $\varphi$ factors out to a map $\varphi_W\colon \Delta_{d;W}\to M_d$.

\enddemo

\proclaim{Lemma 4.1.1} If $A$, $B\in \Delta_d$ are neighbours in $\Delta_d$
with $A<B$ and
if $\zeta A\leq B$, then
$\zeta A<B$, and the ordered dissections
$\zeta A$ and $B$ are neighbours in $\Delta_d$.

\endproclaim

\demo{Proof} The equalities
$\zeta A=B$, $\varphi(\zeta A)=\varphi(B)$, together with Lemma 3.2.1, (iii),
yield $A=B$ which is a contradiction. Hence
$\zeta A<B$. According to Theorem 3.5.1, (i), the fact that $A$ and $B$ are
neighbours implies $B=R_{i,s}A$ for some $i\in [1,d]$ and $s\in [1,d]$ with
$\varepsilon_A(s)=i+1$. Then, using Lemma 3.2.1, (ii), we obtain
$\varphi(B)=\varphi(R_{i,s}A)=\rho_{i,i+1}(\varphi(A))
=\rho_{i,i+1}(\varphi(\zeta A))$. Now, we apply Lemma 3.3.2 for the pair
$\zeta A$ and $B$,
and get the existence of an integer $s_1\in [1,d]$ with
$\varepsilon_A(s_1)=i+1$, such that
$B=R_{i,s_1}(\zeta A)$. The neighbourhood of
$\zeta A$ and $B$ follows from Theorem 3.5.1, (i).

\enddemo

Let $a$, $b\in \Delta_{d;W}$, and $A\in a$, $B\in b$.
We define a partial order $\leq$ on the factor-set $\Delta_{d;W}$ via the rule:
$$
a\leq b\hbox{\sl\ if and only if there exists a\ }\sigma\in W,
\hbox{\sl\ such that\ } \sigma A\leq B.
$$

\proclaim{Theorem 4.1.2}
(i)
Let $a$, $b\in \Delta_{d;W}$, and $A\in a$, $B\in b$ with $A\leq B$.
Then the
restriction $\psi_1$ of the map $\psi_W$ on the union of the intervals
$[\sigma A,B]$, $\sigma\in W$, in $\Delta_d$, is a surjection
$$
\psi_1\colon \cup_{\sigma\in W}[\sigma A,B]\to [a,b]
$$
onto the interval $[a,b]$ in $\Delta_{d;W}$,
and
one has $\psi_1^{-1}((a,b))=\cup_{\sigma\in W}(\sigma A,B)$;

(ii)
let $a$, $b\in T_{d;W}$, and $A\in a$, $B\in b$ with $A\leq B$.
Then the
restriction $\psi_2$ of the map $\psi_W$ on the union of the intervals
$[\sigma A,B]$, $\sigma\in W$, in $T_d$, is a surjection
$$
\psi_2\colon \cup_{\sigma\in W}[\sigma A,B]\to [a,b]
$$
onto the interval $[a,b]$ in $T_{d;W}$,
and
one has $\psi_2^{-1}((a,b))=\cup_{\sigma\in W}(\sigma A,B)$;

\endproclaim

\demo{Proof}
(i) By definition $a\leq b$.
Suppose that
$a\leq c\leq b$, where $c\in \Delta_{d;W}$ and let $C\in c$. There exist
$\sigma$, $\tau\in W$, such that $\sigma A\leq C$ and $\tau C\leq B$.
Then $\tau\sigma A\leq\tau C\leq B$ and $\psi_1(\tau C)=c$, so the surjectivity
of the map $\psi_1$ is proved. Assume that $\sigma A<C<B$, for some $\sigma\in
W$, and some $C\in \Delta_d$.
By definition, $a\leq c\leq b$, where $c=\psi_1(C)$. If
$a=c$, or $c=b$, then $\tau C<C$, or $\tau B<B$, respectively, for an
appropriate $\tau\in W$, which contradicts to Lemma 3.2.1, (iii). Therefore
$\cup_{\sigma\in W}(\sigma A,B)\subset\psi_1^{-1}((a,b))$, and
part (i) holds.

(ii) We note that $c\in T_{d;W}$, and $C\in c$, where $C\in \Delta_d$, yield
$C\in T_d$. Thus, the proof of part (i) holds in this case, too.

\enddemo

\vskip 13pt

4.2. It is said that $a, b\in \Delta_{d;W}$ are {\it adjacent} with $a<b$ if
there exist $A\in a$, and $B\in b$, which are adjacent with $A<B$. In other
words, there exists a pair of integers $(i,j)$ with $1\leq i<j\leq d$, and
such that $\varphi_W(b)=\rho_{i,j}(\varphi_W(a))$.

\proclaim{Theorem 4.2.1}
The elements $a, b\in \Delta_{d;W}$
are adjacent with
$a<b$,
if and only if
there exist $A\in a$, and $B\in b$, with $A<B$, and
there exist a pair of integers $(i,j)$ with $1\leq i<j\leq d$,
an integer $r\geq 1$,
and two sequences $(i_k)_{k=1}^{r+1}$
and $(s_\kappa)_{\kappa=1}^r$ in the interval $[1,d]$, such that:
$$
i=i_1<i_2<\hdots <i_{r+1}=j, \hbox{\sl\ and\ }
\varepsilon_A(s_\kappa)=
i_{\kappa+1},\hbox{\sl\ for all\ }1\leq \kappa\leq r,
$$
and that
$$
B=R_{i_1,s_1}\hdots R_{i_r,s_r}A. \tag 4.2.2
$$

\endproclaim

\demo{Proof}
The necessity holds because of Theorem 3.4.1. For the converse statement
we apply the map $\varphi$ on the equality (4.2.2) and obtain
$\varphi(B)=\rho_{i,j}(\varphi(A))$. Hence $a$ and $b$ are adjacent with
$a<b$.

\enddemo

\proclaim{Theorem 4.2.3}
(i) The elements $a, b\in \Delta_{d;W}$ are neighbours in $\Delta_{d;W}$ with
$a<b$, if and only if
there exist $A\in a$, and $B\in b$, with $A<B$, and
there exist $i\in [1,d]$ and $s\in [1,d]$, such that
$\varepsilon_A(s)=i+1$ and
$B=R_{i,s}(A)$;

(ii) the elements $a, b\in T_{d;W}$
are neighbours in $T_{d;W}$ with
$a<b$,
if and only if
there exist $A\in a$, and $B\in b$, with $A<B$, and
there exist a pair of integers $(i,j)$ with $1\leq i<j\leq d$,
an integer $r\geq 1$,
and two sequences $(i_k)_{k=1}^{r+1}$
and $(s_\kappa)_{\kappa=1}^r$ in the interval $[1,d]$, such that:
$$
j=i+1 \hbox{\ or\ } |A_i|=|A_j|,
$$
and
$$
i=i_1<i_2<\hdots <i_{r+1}=j, \hbox{\sl\ and\ }
\varepsilon_A(s_\kappa)=
i_{\kappa+1},\hbox{\sl\ for all\ }1\leq \kappa\leq r,
$$
and that
$B=R_{i_1,s_1}\hdots R_{i_r,s_r}A$.

\endproclaim

\demo{Proof} Using Lemma 4.1.1, and Theorem 4.1.2, (i) (respectively
(ii)), we get that
$a$ and $b$ are neighbours in $\Delta_{d;W}$
(respectively, in $T_{d;W}$) with $a<b$
if and only
if $A$ and $B$ are neighbours in $\Delta_d$ (respectively, in $T_d$) with
$A<B$. Then Theorem 3.5.1, (i) (respectively (ii)), finishes the proof of
part (i) (respectively, of part (ii)).

\enddemo
\heading{5. Counting of isomers}

\endheading

5.1. The set $T_d$ can be stratified
using the fibres $T_\lambda=\varphi^{-1}(\lambda)$ of the
map $\varphi\colon T_d\to P_d$,
where $\lambda$ runs through the set $P_d$. Clearly, $T_\lambda$ is the
set of all tabloids of shape $\lambda$. Since the symmetric group $S_d$ is
$d$-transitive on $[1,d]$,
the set of fibres $T_\lambda$, $\lambda\in P_d$, coincides with
the set $S_d\backslash T_d$ of $S_d$-orbits in $T_d$.

The orbit $T_\lambda$ contains the tabloid
$I$ with components
$I_1=[1,\lambda_1], I_2=[\lambda_1+1,\lambda_1+\lambda_2],\hdots$, and its
stabilizer
is the subgroup $S_\lambda=S_{\lambda_1}\times\cdots\times S_{\lambda_d}\leq
S_d$. Thus
$$
S_d/S_\lambda\simeq
T_\lambda,\tag 5.1.1
$$
$$
\upsilon S_\lambda\mapsto\upsilon I,
$$
is an isomorphism of $S_d$-sets.

Let us fix a $S_d$-orbit
$T_\lambda$ and consider
the action of the permutation group $W\leq S_d$ on
$T_\lambda$,
which is induced by the action
(1.1.1) of $S_d$. Let us denote by
$T_{\lambda;W}$ the orbit space
$W\backslash T_\lambda$.
Then the isomorphism (5.1.1) of $S_d$-sets can also be considered as an
isomorphism of $W$-sets, and moreover, it factors out to a bijection
$$
W\backslash S_d/S_\lambda\simeq
T_{\lambda;W},
$$
$$
W\upsilon S_\lambda\mapsto\upsilon I,
$$
between the set
of double cosets of $S_d$
modulo
$(W,S_\lambda)$, and
the set
of $W$-orbits
in $T_\lambda$.

Let $A=\upsilon I\in T_\lambda$. The stabilizer $W_A$ of $A$ in the
group $W$ consists of all $\sigma\in W$ such that
$\upsilon^{-1}\sigma\upsilon\in S_\lambda$,
or, equivalently,
$\sigma\in\upsilon S_\lambda\upsilon^{-1}$.
Hence
$W_A=W\cap\upsilon S_\lambda\upsilon^{-1}$.

We fix a one-dimensional character $\chi\colon W\to K$ and a one-dimensional
character $\theta\colon S_\lambda\to K$. For a
given $\upsilon\in S_d$, and $A=\upsilon I$, the rule
$$
\beta_\upsilon\colon W_A\to K,\tag 5.1.2
$$
$$
\beta_\upsilon(\sigma)=\chi(\sigma)\theta(\upsilon^{-1}\sigma\upsilon),
$$
defines a one-dimensional character of the stabilizer $W_A$.

If $B=\tau A$ for some $\tau\in W$, then $B=\tau\upsilon I$ and
$W_B=\tau W_A\tau^{-1}$. For the corresponding
one-dimensional character
$\beta_{\tau\upsilon}\colon W_B\to K$,
we have
$$
\beta_{\tau\upsilon}(\tau\sigma\tau^{-1})
=\chi(\tau\sigma\tau^{-1})\theta(\upsilon^{-1}
\tau^{-1}\tau\sigma\tau^{-1}\tau\upsilon)=
\beta_\upsilon(\sigma),
$$
where $\sigma\in W_A$.
Therefore, given a $W$-orbit $a\in
T_{\lambda;W}$,
the statements
$$
\lq\lq\beta_\upsilon(\sigma)=1 \hbox{\sl\ for any\ } \sigma\in W_A"\tag 5.1.3
$$
are simultaneously true or false regardless of the representative
$A=\upsilon I\in a$.
We denote by
$T_{\lambda;\chi,\theta}$
the subset of
$T_{\lambda;W}$ consisting of those $W$-orbits $a$ for which the
statement (5.1.3) is true for some representative
$A=\upsilon I\in a$, and call them {\it $(\chi,\theta)$-orbits} of the group
$W$.
In particular, $T_{\lambda;1_W,1_{S_\lambda}}=T_{\lambda;W}$.

In case
$\theta=1_{S_\lambda}$ for all $\lambda\in P_d$,
the $(\chi,\theta)$-orbits of the group
shall be called simply {\it
$\chi$-orbits} of the
group $W$. Thus the $\chi$-orbits are those $W$-orbits $a\in T_d$ for
which there exists a tabloid $A\in a$ such that the character $\chi$ is
identically $1$ on its stabilizer $W_A$ (see (5.1.3)). Then the last condition
holds for
all tabloids $A\in a$. We set
$T_{\lambda;\chi}=T_{\lambda;\chi,1_{S_\lambda}}$, and
$T_{d;\chi}=\cup_{\lambda\in P_d}T_{\lambda;\chi}$.

We introduce the following
families of non-negative integers:
$n_{\lambda;\chi,\theta}=
|T_{\lambda;\chi,\theta}|$,
$n_{\lambda;\chi}=|T_{\lambda;\chi}|$, and
$n_{\lambda;W}=
|T_{\lambda;W}|$,
where $\lambda\in P_d$.
Note that $n_{\lambda;\chi}=n_{\lambda;\chi,1_{S_\lambda}}$, and
$n_{\lambda;W}=
n_{\lambda;1_W}=
n_{\lambda;1_W,1_{S_\lambda}}$.

\vskip 13pt

5.2. Now, our aim is to find an explicit formula for the number
$n_{\lambda;\chi,\theta}$
 of $(\chi,\theta)$-orbits of the group $W$ in the
set $T_\lambda$, where
$\lambda\in P_d$.
We shall use terminology,
notation and results from [4] and [7].

For any finite set $X$ we denote by $|X|$ the number
of its elements. For a partition $\lambda\in P_d$ we shall use also the
notation
$(1^{m_1},2^{m_2},\hdots, d^{m_d})$,
where $m_k$ is the number of
the parts of $\lambda$, which are equal to $k$, $1\leq k\leq d$. Given
a permutation $\zeta\in S_d$, we denote by $\varrho(\zeta)$, and also, by
$(1^{c_1\left(\zeta\right)},
2^{c_2\left(\zeta\right)},\hdots,
d^{c_d\left(\zeta\right)})$
the corresponding partition of
the number $d$. We set
$$
C(W,S_\lambda)=\{(\sigma,\eta)\in W\times
S_\lambda\mid\varrho(\sigma)=\varrho(\eta)\}.
$$
Let $t$ be the length of the partition $\lambda$. Then
$S_\lambda=S_{\lambda_1}\times\cdots\times S_{\lambda_t}$, so any $\eta\in
S_\lambda$ has the form $\eta=\eta_1\hdots\eta_t$, where $\eta_k\in
S_{\lambda_k}$. Thus
$\varrho(\eta)=\varrho(\eta_1)\cup\hdots\cup\varrho(\eta_t)$, where
$\varrho(\eta_k)\in P_{\lambda_k}$.

The one-dimensional character $\theta$ has a unique decomposition
$\theta=\theta_1\hdots\theta_t$, where $\theta_k$ is either the signature or
the unit character of $S_{\lambda_k}$.
We set
$$
L_\lambda=\{(\alpha,\alpha^{\left(1\right)},\hdots,\alpha^{\left(t\right)})\in
P_d\times P_{\lambda_1}\times\cdots\times P_{\lambda_t}\mid
\alpha=\alpha^{\left(1\right)}\cup\hdots\cup\alpha^{\left(t\right)}\},
$$
and define a map
$$
\gamma_1\colon W\times S_{\lambda_1}\times\cdots\times S_{\lambda_t}\to
P_d\times P_{\lambda_1}\times\cdots\times P_{\lambda_t},
$$
$$
(\sigma,\eta_1\hdots,\eta_t)\mapsto
(\varrho(\sigma),\varrho(\eta_1),\hdots,\varrho(\eta_t)).
$$
Then
$C(W,S_\lambda)=\gamma_1^{-1}(L_\lambda)$. Let $L^\prime(W,S_\lambda)\subset
L_\lambda$ be the image of $C(W,S_\lambda)$ via the map $\gamma_1$. The
restriction of $\gamma_1$ on $C(W,S_\lambda)$ is a surjective map
$$
\gamma\colon C(W,S_\lambda)\to L^\prime(W,S_\lambda).
$$
If
$(\alpha,\alpha^{\left(1\right)},\hdots,\alpha^{\left(t\right)})\in
L^\prime(W,S_\lambda)$, then
$$
\gamma^{-1}(\alpha,\alpha^{\left(1\right)},\hdots,\alpha^{\left(t\right)})=
W_\alpha\times K_{\alpha^{\left(1\right)}}\times\cdots\times
K_{\alpha^{\left(t\right)}},\tag 5.2.1
$$
where $W_\alpha$ is the subset of the group $W$,
consisting of all permutations of cyclic type $\alpha$, and
$K_{\alpha^{\left(k\right)}}$ is the conjugacy class in $S_{\lambda_k}$,
corresponding to the partition
$\alpha^{\left(k\right)}\in P_{\lambda_k}$. The set $W_\alpha$ is a union of
conjugacy classes of the group $W$:
$$
W_\alpha= C_1^{\left(\alpha\right)}\cup\hdots\cup
C_{i_\alpha}^{\left(\alpha\right)}.\tag
5.2.2
$$

We set
$$
L(W,S_\lambda)=L^\prime(W,S_\lambda)\backslash\{((1), (1),\hdots, (1))\}.\tag
5.2.3
$$

Let $h_\lambda=h_{\lambda_1}h_{\lambda_2}\hdots$, where $h_{\lambda_k}$ is the
$\lambda_k$-th complete
symmetric function (see [4, Ch. I, Sec. 2]).

\proclaim{Lemma 5.2.4}
Let $W\leq S_d$ be a
permutation group and $\chi\colon W\to K$ be a one-dimensional character.

(i) The characteristic of the induced
monomial representation $Ind_W^{S_d}(\chi )$ is equal to the generalized cyclic
index
$$
Z(\chi;p_1,\hdots ,p_d)
=\frac{1}{|W|}\sum_{\sigma\in W}\chi (\sigma)
p_1^{c_1\left(\sigma\right )}\hdots p_d^{c_d\left (\sigma\right )},
$$
where $p_s=x_0^s+x_1^s+\cdots$ are the power sums;

(ii) one has
$n_{\lambda;\chi,\theta}=
\langle Ind_W^{S_d}(\chi),
Ind_{S_\lambda}^{S_d}(\theta)\rangle_{S_d}$;

(iii) one has
$n_{\lambda;\chi}=
\langle Z(\chi;p_1,\hdots ,p_d),
h_\lambda\rangle$;

(iv) if $Z(\chi;p_1,\hdots ,p_d)= \sum_{\lambda\in P_d}a_\lambda m_\lambda$,
where $m_\lambda$ are the monomial symmetric functions, then $a_\lambda=
n_{\lambda;\chi}$.

\endproclaim

\demo{Proof} (i) Let $\psi$ be the map which assigns to each substitution
$\zeta\in S_d$ the symmetric function $p_{\varrho\left(\zeta\right)}$
(see [4, Ch. I, Sec. 7]).
According to Frobenius reciprocity law, we have
$$
ch(Ind_W^{S_d}(\chi))=\langle Ind_W^{S_d}(\chi),\psi\rangle_{S_d}=
$$
$$
\langle\chi,Res_W^{S_d}(\psi)\rangle_W=
Z(\chi;p_1,\hdots ,p_d).
$$

(ii) Using Frobenius reciprocity law and [7, Ch. II, 7.4, Proposition 15], we
have $$
\langle Ind_W^{S_d}(\chi),
Ind_{S_\lambda}^{S_d}(\theta)\rangle_{S_d}=
\langle\chi,Res_W^{S_d}Ind_{S_\lambda}^{S_d}(\theta)\rangle_{S_d}=
$$
$$
\langle\chi,\sum_{\upsilon\in\Upsilon}
Ind_{W_\upsilon}^W(\theta_\upsilon)\rangle_W,
$$
where $\Upsilon$ is a system of representatives of the double cosets of $S_d$
modulo
$(W,S_\lambda)$, and $\theta_\upsilon$ is the one-dimensional character of the
group $W\cap\upsilon S_\lambda\upsilon^{-1}$ given by the formula
$\theta_\upsilon(x)=\theta(\upsilon^{-1}x\upsilon)$. Further,
$$
\langle Ind_W^{S_d}(\chi),
Ind_{S_\lambda}^{S_d}(\theta)\rangle_{S_d}=
\sum_{\upsilon\in\Upsilon}
\langle\chi,Ind_{W_\upsilon}^W(\theta_\upsilon)\rangle_W=
$$
$$
\sum_{\upsilon\in\Upsilon}
\langle Res_{W_\upsilon}^W(\chi),\theta_\upsilon\rangle_{W_\upsilon}=
|\{\upsilon\in\Upsilon\mid Res_{W_\upsilon}^W(\chi)=\theta_\upsilon\}|.
$$
Since $\theta^2=1$, then
$$
\{\upsilon\in\Upsilon\mid Res_{W_\upsilon}^W(\chi)=\theta_\upsilon\}=
\{\upsilon\in\Upsilon\mid
\beta_\upsilon(\sigma)=1 \hbox{\sl\ for any\ } \sigma\in W_\upsilon\}.
$$
Therefore, using the isomorphism (5.1.1) of $W$-sets we get
$$
\langle Ind_W^{S_d}(\chi),
Ind_{S_\lambda}^{S_d}(\theta)\rangle_{S_d}=
n_{\lambda;\chi,\theta}.
$$

(iii) Indeed,
the characteristic map $ch$ is an isometric isomorphism of rings
(see [4, Ch. I, Sec. 7, 7.3]),
so in particular,
$$
\langle Ind_W^{S_d}(\chi),
Ind_{S_\lambda}^{S_d}(1_{S_\lambda})\rangle_{S_d}=
\langle ch(Ind_W^{S_d}(\chi)),
ch(Ind_{S_\lambda}^{S_d}(1_{S_\lambda}))\rangle .
$$
Evidently, $ch(Ind_{S_\lambda}^{S_d}(1_{S_\lambda}))=h_\lambda$.
According to Lemma 5.2.4, (i), we have
$$
ch(Ind_W^{S_d}(\chi))= Z(\chi;p_1,\hdots ,p_d).
$$
Therefore
$$
n_{\lambda;\chi}=
n_{\lambda;\chi,1_{S_\lambda}}=
\langle Ind_W^{S_d}(\chi),
Ind_{S_\lambda}^{S_d}(1_{S_\lambda})\rangle_{S_d}=
\langle
Z(\chi;p_1,\hdots ,p_d),
h_\lambda\rangle .
$$

(iv) Using part (iii) we obtain
$$
n_{\lambda;\chi}=
\langle\sum_{\alpha\in P_d}a_\alpha m_\alpha,
h_\lambda\rangle =
\sum_{\lambda\in P_d}a_\alpha\langle m_\alpha,
h_\lambda\rangle = a_\lambda.
$$
The last equality holds because of
[4, Ch. I, Sec. 4, 4.5].

\enddemo

\proclaim{Theorem 5.2.5}
Let $W$, $W^\prime\leq S_d$ be two
permutation groups. Then the following four statements are equivalent:

(i) One has $n_{W;\lambda}=n_{W^\prime;\lambda}$ for all $\lambda\in P_d$;

(ii) one has
$$
Z(1_W;p_1,\hdots ,p_d)=
Z(1_{W^\prime};p_1,\hdots ,p_d);
$$

(iii) the induced monomial representations
$Ind_W^{S_d}(1_W)$, and
$Ind_{W^\prime}^{S_d}(1_{W^\prime})$,
of the symmetric group $S_d$, are isomorphic;

(iv) there exists a one-one correspondence between the groups $W$ and
$W^\prime$,
such that the corresponding permutations have the same type of cycle
decomposition.

\endproclaim

\demo{Proof} Lemma 5.2.4, (i) and (iv), applied for $1_W$, and
$1_{W^\prime}$, and
[4, Ch. I, Sec. 7, 7.3]) yield the  equivalence of (i), (ii), and (iii).
It is easily seen that the equality of cyclic indices in (ii) is equivalent to
(iv).

\enddemo

\proclaim{Remark 5.2.6} {\rm According to [3, IV], two permutation groups
$W$, $W^\prime\leq S_d$ which satisfy (iv) are said to be {\it literally
conformal}. In [5, Ch. I, Sec. 25], it is shown that each of (ii) and (iv) is
equivalent to the so called {\it combinatorial equivalence} of $W$ and
$W^\prime$.}

\endproclaim

For any $\lambda\in P_d$,
$\lambda=(1^{m_1},2^{m_2},\hdots, d^{m_d})$, we set
$z_\lambda=1^{m_1}m_1!2^{m_2}m_2!\hdots d^{m_d}m_d!$.

\proclaim{Theorem 5.2.7} One has
$$
n_{\lambda;\chi,\theta}=
\frac{d!}{|W|\lambda_1!\hdots\lambda_d!}+
$$
$$
\frac{1}{|W|}
\sum_{\left(\alpha,\alpha^{\left(1\right)},
\hdots,\alpha^{\left(t\right)}\right)\in
L\left(W,S_\lambda\right)}
(\sum_{i=1}^{i_\alpha}|C_i^{\left(\alpha\right)}|\chi
(C_i^{\left(\alpha\right)}))
\frac{z_\alpha}{z_{\alpha^{\left(1\right)}}
\hdots z_{\alpha^{\left(t\right)}}}
\theta_1(K_{\alpha^{\left(1\right)}})\hdots
\theta_t(K_{\alpha^{\left(t\right)}}).
$$

\endproclaim

\demo{Proof}
The characteristic map $ch$ is an isometry,
and Lemma 5.2.4, (i), holds, so
$$
\langle Ind_W^{S_d}(\chi),
Ind_{S_\lambda}^{S_d}(\theta)\rangle_{S_d}=
\langle ch(Ind_W^{S_d}(\chi)),
ch(Ind_{S_\lambda}^{S_d}(\theta))\rangle =
$$
$$
\langle  Z(\chi;p_1,\hdots ,p_d), Z(\theta;p_1,\hdots ,p_d)\rangle=
\langle \frac{1}{|W|}\sum_{\sigma\in
W}\chi(\sigma)p_{\varrho\left(\sigma\right)},
\frac{1}{|S_\lambda|}\sum_{\eta\in
S_\lambda}\theta(\eta)p_{\varrho\left(\eta\right)}\rangle=
$$
$$
\frac{1}{|W||S_\lambda|}\sum_{\sigma\in W}\sum_{\eta\in S_\lambda}
\chi(\sigma)\theta(\eta)\langle
p_{\varrho\left(\sigma\right)},p_{\varrho\left(\eta\right)}\rangle.
$$
According to [4, Ch. I, Sec. 4, 4.7], we obtain
$$
\langle Ind_W^{S_d}(\chi),
Ind_{S_\lambda}^{S_d}(\theta)\rangle_{S_d}=
\frac{1}{|W||S_\lambda|}\sum_{\left(\sigma,\eta\right)\in
C\left(W, S_\lambda\right)}
\chi(\sigma)\theta(\eta)z_{\varrho\left(\sigma\right)}.\tag 5.2.8
$$
Further, we use the partition of the set $C(W,S_\lambda)$ into the
fibres of the surjective map $\gamma$, as well as their representation (5.2.1).
Thus, we have
$$
\langle Ind_W^{S_d}(\chi),
Ind_{S_\lambda}^{S_d}(\theta)\rangle_{S_d}=
$$
$$
\frac{1}{|W||S_\lambda|}\sum_{
\left(\alpha,\alpha^{\left(1\right)},\hdots,\alpha^{\left(t\right)}\right)\in
L^\prime\left(W,S_\lambda\right)}
\sum_{\left(\sigma,\eta_1,\hdots,\eta_t\right)\in
W_\alpha\times K_{\alpha^{\left(1\right)}}\times\cdots\times
K_{\alpha^{\left(t\right)}}}
\chi(\sigma)\theta_1(\eta_1)\hdots\theta_t(\eta_t)z_\alpha=
$$
$$
\frac{1}{|W|\lambda_1!\hdots\lambda_t!}\sum_{
\left(\alpha,\alpha^{\left(1\right)},\hdots,\alpha^{\left(t\right)}\right)\in
L^\prime\left(W,S_\lambda\right)}
(\sum_{i=1}^{i_\alpha}|C_i^{\left(\alpha\right)}|\chi
(C_i^{\left(\alpha\right)}))
\prod_{k=1}^l\frac{\lambda_k!}{z_{\alpha^{\left(k\right)}}}
\theta_k(K_{\alpha^{\left(k\right)}})
z_\alpha=
$$
$$
\frac{1}{|W|}
\sum_{\left(\alpha,\alpha^{\left(1\right)},
\hdots,\alpha^{\left(t\right)}\right)\in
L^\prime\left(W,S_\lambda\right)}
(\sum_{i=1}^{i_\alpha}|C_i^{\left(\alpha\right)}|\chi
(C_i^{\left(\alpha\right)}))
\frac{z_\alpha}{z_{\alpha^{\left(1\right)}}
\hdots z_{\alpha^{\left(t\right)}}}
\theta_1(K_{\alpha^{\left(1\right)}})\hdots
\theta_t(K_{\alpha^{\left(t\right)}})=
$$
$$
\frac{d!}{|W|\lambda_1!\hdots\lambda_d!}+
$$
$$
\frac{1}{|W|}
\sum_{\left(\alpha,\alpha^{\left(1\right)},
\hdots,\alpha^{\left(t\right)}\right)\in
L\left(W,S_\lambda\right)}
(\sum_{i=1}^{i_\alpha}|C_i^{\left(\alpha\right)}|\chi
(C_i^{\left(\alpha\right)}))
\frac{z_\alpha}{z_{\alpha^{\left(1\right)}}
\hdots z_{\alpha^{\left(t\right)}}}
\theta_1(K_{\alpha^{\left(1\right)}})\hdots
\theta_t(K_{\alpha^{\left(t\right)}}).
$$
In the last two equalities we make use of (5.2.2) and (5.2.3). Now,
Lemma 5.2.4, (ii), yields the result.

\enddemo

The specialization
$\chi=1_W$, and $\theta=1_{S_\lambda}$, in Theorem 5.2.7 entails

\proclaim{Corollary 5.2.9} One has
$$
n_{\lambda;W}=
\frac{d!}{|W|\lambda_1!\hdots\lambda_d!}+
$$
$$
\frac{1}{|W|}
\sum_{\left(\alpha,\alpha^{\left(1\right)},
\hdots,\alpha^{\left(t\right)}\right)\in
L\left(W,S_\lambda\right)}|W_\alpha|
\frac{z_\alpha}{z_{\alpha^{\left(1\right)}}
\hdots z_{\alpha^{\left(t\right)}}}.
$$

\endproclaim

\proclaim{Corollary 5.2.10 (Ruch's formula)} One has
$$
n_{\lambda;W}=
\frac{n!}{|W||S_\lambda|}\sum_{\alpha\in P_d}
\frac{|W_\alpha||(S_\lambda)_\alpha|}
{|K_\alpha|}.
$$

\endproclaim

\demo{Proof} Using the equality (5.2.8) for
$\chi=1_W$, and $\theta=1_{S_\lambda}$, we have
$$
n_{\lambda;W}=
\frac{1}{|W||S_\lambda|}\sum_{\left(\sigma,\eta\right)\in
C\left(W, S_\lambda\right)}
z_{\varrho\left(\sigma\right)}=
$$
$$
\frac{1}{|W||S_\lambda|}\sum_{\alpha\in P_d}
\sum_{\left(\sigma,\eta\right)\in
W\times
S_\lambda ,\varrho\left(\sigma\right)=\varrho\left(\eta\right)=\alpha}
z_\alpha=
$$
$$
\frac{1}{|W||S_\lambda|}\sum_{\alpha\in P_d}|W_\alpha||(S_\lambda)_\alpha|
z_\alpha=
$$
$$
\frac{n!}{|W||S_\lambda|}\sum_{\alpha\in P_d}
\frac{|W_\alpha||(S_\lambda)_\alpha|}
{|K_\alpha|}.
$$
\enddemo

\proclaim{Remark 5.2.11} {\rm Let $\Gamma$ and $\Delta$ be two graphs with $d$
vertices, and with
automorphism groups $W$ and $S_\lambda$, respectively. The number calculated
in Theorem 5.2.7 coincides with the number of superpositions of $\Gamma$ and
$\Delta$, such that the one-dimensional character (5.1.2) is identically $1$ on
their stabilizers (that is, their automorphism groups). The last number can
also be obtained by
an appropriate generalization of Redfield's superposition theorem (see [1]).}

\endproclaim

\vskip 13pt

5.3. Here we shall consider
the family of non-negative integers
$n_{\lambda;\chi}=|T_{\lambda;\chi}|$, $\lambda\in P_d$.

\proclaim{Theorem 5.3.1}
Let  $\chi$
be a one-dimensional character $\chi$ of the group $W\leq S_d$,
and let $\lambda, \mu\in P_d$. If $\lambda\leq\mu$,
then
$n_{\lambda;\chi}\geq n_{\mu;\chi}$.

\endproclaim

\demo{Proof} According to Lemma 5.2.4, (iii), we have
$n_{\lambda;\chi}
-n_{\mu;\chi}=
\langle Z(\chi;p_1,\hdots ,p_d),
h_\lambda-h_\mu\rangle$. Then [4, Ch. I, Sec. 7, Example 9 (b)] implies that
the difference
$h_\lambda-h_\mu$ is a non-negative integral linear combination of the
Schur functions $s_\nu$, $\nu\in P_d$. On the other hand,
$Z(\chi;p_1,\hdots ,p_d)$ is the characteristic of the induced monomial
representation
$Ind_W^{S_d}(\chi )$, so it also is a linear combination of
$s_\nu$ with non-negative integral coefficients. Therefore the above scalar
product is non-negative.

\enddemo

The specialization $\chi=1_W$ yields

\proclaim{Corollary 5.3.2}
 If $\lambda, \mu\in P_d$, and $\lambda\leq\mu$,
then
$n_{\lambda;W}\geq n_{\mu;W}$.

\endproclaim

\vskip 13pt

5.4.
Below,
Lemma 5.4.3 for $I=T_d$ gives a combinatorial interpretation of the set
$T_{d;\chi}$ of all $\chi$-orbits.
We shall work in a more general setup.

Let $W$ be a finite group which acts on a set $I.$ For each element $i\in
I$ we denote by $W_i$ its stabilizer in $W$. Let
$\chi$ be a one-dimensional character of the group $W$ with kernel $H\leq W$.

\proclaim {Lemma 5.4.1} The following statements hold:

(i) The inclusion
$W_i\leq H$, and the equality $|W_i:H_i|=1$ are equivalent for any $i\in I$;

(ii) if $O$ is a $W$-orbit in $I$,
then all $H$-orbits in $O$ have the same
number of elements, and their number is a divisor of the index $|W:H|$.

\endproclaim

\demo {Proof}
(i) It is enough to note that $H_i=H\cap W_i$.

(ii) Let $i\in O$. Since $H$ is a normal subgroup of $W$,
then $\sigma H_i\sigma^{-1}\leq H$ for all $\sigma\in W$. Therefore
$|H:H_{\sigma
i}|=|H:\sigma H_i\sigma^{-1}|=|H:H_i|$, that is, each $H$-orbit in $O$ has the
same number of elements. Then using
the equality
$$
|W:H||H:H_i|=|W:W_i||W_i:H_i|, \tag 5.4.2
$$
where $i\in I$, we obtain immediately that
the number of $H$-orbits in $O$  is a divisor of $|W:H|$.

\enddemo

If $O$ is a $W$-orbit in $I$, and if one has
$W_i\leq H$ for some $i\in O$ (and, hence, for all $i\in O$),
then $O$ is said to be a {\it $\chi$-orbit}.

\proclaim {Lemma 5.4.3} The following two statements are equivalent:

(i) The $W$-orbit $O$ is a $\chi$-orbit;

(ii) the $W$-orbit $O$ contains exactly
$|W:H|$ in number $H$-orbits;

(iii) the $W$-orbit $O$ contains maximum
number $H$-orbits.

\endproclaim

\demo {Proof} Let $O$ be a $\chi$-orbit.
The equality (5.4.2)
and Lemma 5.4.1, (i),
yield $|W:H||H:H_i|=|O|$ for $i\in O$. Because of Lemma  5.4.1, (ii), the
indices $|H:H_i|$ do not depend on $i\in O$ and all are equal to the number of
elements of any $H$-orbit in $O$.  Therefore (ii) holds. Conversely, suppose
that the
$W$-orbit $O$ contains exactly
$|W:H|$ in number $H$-orbits, and
let $i\in O$. Lemma 5.4.1, (ii), implies
$|W:W_i|=|O|=|W:H||H:H_i|$. Comparing with (5.4.2), we obtain
$|W_i:H_i|=1$. Due to Lemma 5.4.1, (i), $O$ is a $\chi$-orbit.
Finally, Lemma 5.4.1, (ii), yields that part (iii) is equivalent to part
(ii).

\enddemo

\heading
6. First applications
\endheading

6.1. Now, we apply Corollary 5.2.9 to obtain Kauffmann formulae for the number
of the derivatives of naphthalene, $C_{10}H_8$.

The group of substitution isomerism of
naphthalene is the subgroup $G$ of $S_8$, consisting of the elements
$$
(1), (12)(34)(56)(78), (13)(24)(57)(68), (14)(23)(58)(67)
$$
(see [3, IX, D]).
The unit $(1)$ of $G$ produces the term
$\frac{1}{4}\frac{8!}{\lambda_1!\hdots\lambda_8!}$.
The other $3$ elements of $G$
have cyclic structure $(2^4)$, so $|G_{\left(2^4\right)}|=3$.

Suppose that the
set $L\left(G,S_\lambda\right)$ contains an element
$(\alpha,\alpha^{\left(1\right)},\hdots,\alpha^{\left(t\right)})$ with
$\alpha=(2^4)$.
Then $\alpha$ is to be the cyclic type of an element of
$S_\lambda=S_{\lambda_1}\times\cdots\times S_{\lambda_t}\leq S_8$, where $t$ is
the length of the partition $\lambda$ of $8$.

In case at least one of the
components $\lambda_k$ is odd, we establish a contradiction, so in this case
$$
n_{\lambda;G}=\frac{1}{4}\frac{8!}{\lambda_1!\hdots\lambda_8!}.
$$

In the rest of the cases, all
components $\lambda_k$ have to be even, so
$\lambda=(2\mu_1, 2\mu_2,\hdots,2\mu_t)$, where
$\mu=(\mu_1,\mu_2,\hdots,\mu_t)$ is a partition of $4$.
Now,
$$
L\left(W,S_\lambda\right)=\{((2^4), (2^{\mu_1}),
(2^{\mu_2}),\hdots,(2^{\mu_t}))\}
$$
and
$$
n_{\lambda;G}=
\frac{1}{4}\frac{8!}{\lambda_1!\hdots\lambda_8!}+
\frac{1}{4}
|W_{\left(2^4\right)}|
\frac{z_{\left(2^4\right)}}{z_{\left(2^{\mu_1}\right)}
z_{\left(2^{\mu_2}\right)}\hdots}=
$$
$$
\frac{1}{4}\frac{8!}{\lambda_1!\hdots\lambda_8!}+
\frac{3}{4}
\frac{2^44!}{2^{\mu_1}\mu_1!
2^{\mu_2}\mu_2!\hdots}=
$$
$$
\frac{1}{4}\frac{8!}{\lambda_1!\hdots\lambda_8!}+
\frac{3}{4}
\frac{4!}{\left(\frac{\lambda_1}{2}\right)!
\left(\frac{\lambda_2}{2}\right)!\hdots}.
$$
Thus, we have obtained Kauffmann formulae.

\vskip 13pt

6.2. This subsection is devoted to chiral pairs.
Let $\Sigma$ be a skeleton with $d$
unsatisfied single valences. Suppose that among the substitution
derivatives of a given parent substance with skeleton $\Sigma$ there is an
chiral pair. Then according to Lunn-Senior Thesis 1.5.1, (2a), the
group
$G^\prime\leq S_d$ of stereoisomerism contains the group $G$ of substitution
isomerism as a subgroup of index $2$. In particular, $G$ is a normal subgroup
of $G^\prime$. Let $\chi_e\colon G^\prime\to K$ be the one-dimensional complex
valued character with kernel $G$. We have $\chi_e(\sigma)=1$ for $\sigma\in G$
and $\chi_e(\sigma)=-1$ for $\sigma\in G^\prime\backslash G$. Lemma 5.4.1,
(ii), for $I=T_d$ and $W=G^\prime$ implies that each $G^\prime$-orbit
$O$ contains either two or one $G$-orbit. Lunn-Senior Thesis 1.5.1, part (2ae),
and part (2ad), makes the corresponding identifications with the
chiral pairs, and with the diastereomers, respectively.
Lemma 5.4.3 applied
 for $I=T_d$ and $W=G^\prime$ shows that the set $T_{d;\chi_e}$ of
$\chi_e$-orbits contains the chiral pairs. In particular, part (2a) of
the Extended Lunn-Senior Thesis 1.6.1 is justified.

Now, as a direct consequence of Theorem 5.3.1 for $W=G^\prime$ and
$\chi=\chi_e$, we obtain a result of E. Ruch of special beauty.

\proclaim {Theorem 6.2.1 (Ruch)} If a distribution of ligands according to the
partition $\mu$
amounts to a chiral molecule, and $\lambda$ is dominated by $\mu$ , then also
a distribution according to $\lambda$  yields a chiral molecule.

\endproclaim

\vskip 13pt

6.3.  Now, using our approach, we shall present the K\"orner's relations
between
the di-, and tri-substitution derivatives of benzene, $C_6H_6$.  The exposition
below follows that of Lunn and Senior.

Let the skeleton $\Sigma$ be the six carbon atom
ring of benzene.
According to [3, VI], the group $G$ of substitution isomerism of benzene
has the following elements:
$$
(1), (123456), (135)(246), (14)(25)(36), (153)(264), (165432),
$$
$$
(13)(46),(12)(36)(45), (26)(35), (16)(25)(34), (15)(24), (14)(23)(56).
$$
Clearly, $G$ coincides with the dihedral group $D_6=\langle r,s\rangle$, where
$r=(123456)$, and $s=(13)(46)$.

Case 1. $\lambda=(4,2)$.

There are three isomeric forms of the di-substitution products of benzene,
called {\it para}, {\it ortho}, and {\it meta} derivatives. Therefore
$N_{\left(4,2\right);\Sigma}=3$, which is in
agreement with the equality $n_{\left(4,2\right);G}=3$.

We have $T_{\left(4,2\right);G}=\{a_{\left(4,2\right)},b_{\left(4,2\right)},
c_{\left(4,2\right)}\}$, where:

$a_{\left(4,2\right)}$ is the $G$-orbit
$$
\{(\{2,3,5,6\},\{1,4\},\emptyset ,\emptyset,\emptyset,\emptyset),
(\{1,3,4,6\},\{2,5\},\emptyset ,\emptyset,\emptyset,\emptyset),
(\{1,2,4,5\},\{3,6\},\emptyset ,\emptyset,\emptyset,\emptyset)\}
$$
of the tabloid
$A^{\left(4,2\right)}=
(\{2,3,5,6\},\{1,4\},\emptyset ,\emptyset,\emptyset,\emptyset)$;

$b_{\left(4,2\right)}$ is the $G$-orbit
$$
\{(\{1,2,3,4\},\{5,6\},\emptyset ,\emptyset,\emptyset,\emptyset),
(\{2,3,4,5\},\{1,6\},\emptyset ,\emptyset,\emptyset,\emptyset),
(\{3,4,5,6\},\{1,2\},\emptyset ,\emptyset,\emptyset,\emptyset),
$$
$$
(\{1,4,5,6\},\{2,3\},\emptyset ,\emptyset,\emptyset,\emptyset),
(\{1,2,5,6\},\{3,4\},\emptyset ,\emptyset,\emptyset,\emptyset),
(\{1,2,3,6\},\{4,5\},\emptyset ,\emptyset,\emptyset,\emptyset)\}
$$
of the tabloid
$B^{\left(4,2\right)}=
(\{1,2,3,4\},\{5,6\},\emptyset ,\emptyset,\emptyset,\emptyset)$;

$c_{\left(4,2\right)}$ is the $G$-orbit
$$
\{(\{2,4,5,6\},\{1,3\},\emptyset ,\emptyset,\emptyset,\emptyset),
(\{1,3,5,6\},\{2,4\},\emptyset ,\emptyset,\emptyset,\emptyset),
(\{1,2,4,6\},\{3,5\},\emptyset ,\emptyset,\emptyset,\emptyset),
$$
$$
(\{1,2,3,5\},\{4,6\},\emptyset ,\emptyset,\emptyset,\emptyset),
(\{2,3,4,6\},\{1,5\},\emptyset ,\emptyset,\emptyset,\emptyset),
(\{1,3,4,5\},\{2,6\},\emptyset ,\emptyset,\emptyset,\emptyset)\}
$$
of the tabloid
$C^{\left(4,2\right)}=
(\{2,4,5,6\},\{1,3\},\emptyset ,\emptyset,\emptyset,\emptyset)$.

Case 2. $\lambda=(3^2)$.

The tri-substitution products of benzene exist in three isomeric forms if all
the substituents are the same. They are known as {\it asymmetrical}, {\it
vicinal}, and {\it symmetrical} derivatives. Thus
$N_{\left(3^2\right);\Sigma}=3$, which agrees with $n_{\left(3^2\right);G}=3$.

We have $T_{\left(3^2\right);G}=\{a_{\left(3^2\right)},b_{\left(3^2\right)},
c_{\left(3^2\right)}\}$, where:

$a_{\left(3^2\right)}$ is the $G$-orbit
$$
\{(\{1,2,4\},\{3,5,6\},\emptyset ,\emptyset,\emptyset,\emptyset),
(\{2,3,5\},\{1,4,6\},\emptyset ,\emptyset,\emptyset,\emptyset),
(\{3,4,6\},\{1,2,5\},\emptyset ,\emptyset,\emptyset,\emptyset),
$$
$$
(\{1,4,5\},\{2,3,6\},\emptyset ,\emptyset,\emptyset,\emptyset),
(\{2,5,6\},\{1,3,4\},\emptyset ,\emptyset,\emptyset,\emptyset),
(\{1,3,6\},\{2,4,5\},\emptyset ,\emptyset,\emptyset,\emptyset),
$$
$$
(\{2,3,6\},\{1,4,5\},\emptyset ,\emptyset,\emptyset,\emptyset),
(\{1,2,5\},\{3,4,6\},\emptyset ,\emptyset,\emptyset,\emptyset),
(\{1,4,6\},\{2,3,5\},\emptyset ,\emptyset,\emptyset,\emptyset),
$$
$$
(\{3,5,6\},\{1,2,4\},\emptyset ,\emptyset,\emptyset,\emptyset),
(\{2,4,5\},\{1,3,6\},\emptyset ,\emptyset,\emptyset,\emptyset),
(\{1,3,4\},\{2,5,6\},\emptyset ,\emptyset,\emptyset,\emptyset)\}
$$
of the tabloid
$A^{\left(3^2\right)}=
(\{1,2,4\},\{3,5,6\},\emptyset,\emptyset,\emptyset,\emptyset)$;

$b_{\left(3^2\right)}$ is the $G$-orbit
$$
\{(\{1,2,3\},\{4,5,6\},\emptyset ,\emptyset,\emptyset,\emptyset),
(\{2,3,4\},\{1,5,6\},\emptyset ,\emptyset,\emptyset,\emptyset),
(\{3,4,5\},\{1,2,6\},\emptyset ,\emptyset,\emptyset,\emptyset),
$$
$$
(\{4,5,6\},\{1,2,3\},\emptyset ,\emptyset,\emptyset,\emptyset),
(\{1,5,6\},\{2,3,4\},\emptyset ,\emptyset,\emptyset,\emptyset),
(\{1,2,6\},\{3,4,5\},\emptyset ,\emptyset,\emptyset,\emptyset)\}
$$
of the tabloid
$B^{\left(3^2\right)}=
(\{1,2,3\},\{4,5,6\},\emptyset,\emptyset,\emptyset,\emptyset)$;

$c_{\left(3^2\right)}$ is the $G$-orbit
$$
\{(\{1,3,5\},\{2,4,6\},\emptyset ,\emptyset,\emptyset,\emptyset),
(\{2,4,6\},\{1,3,5\},\emptyset ,\emptyset,\emptyset,\emptyset)\}
$$
of the tabloid
$C^{\left(3^2\right)}=
(\{1,3,5\},\{2,4,6\},\emptyset,\emptyset,\emptyset,\emptyset)$.

Since
$A^{\left(3^2\right)}<(135)(246)A^{\left(4,2\right)}$,
$A^{\left(3^2\right)}<B^{\left(4,2\right)}$,
$A^{\left(3^2\right)}<(135)(246)C^{\left(4,2\right)}$,
$B^{\left(3^2\right)}<B^{\left(4,2\right)}$,
$B^{\left(3^2\right)}<(14)(25)(36)C^{\left(4,2\right)}$, and
$C^{\left(3^2\right)}<(123456)C^{\left(4,2\right)}$,
we have
$$
a_{\left(3^2\right)} < a_{\left(4,2\right)},\hbox{\ }
a_{\left(3^2\right)} < b_{\left(4,2\right)},\hbox{\ }
a_{\left(3^2\right)} < c_{\left(4,2\right)},
$$
$$
b_{\left(3^2\right)} < b_{\left(4,2\right)},\hbox{\ }
b_{\left(3^2\right)} < c_{\left(4,2\right)},\hbox{\ }
c_{\left(3^2\right)} < c_{\left(4,2\right)}.
$$

The above inequalities coincide with the classical K\"orner relations between
di- and tri-substitution products of benzene, and serve for complete
identification of these six derivatives:

$$
\left.
\matrix \format \c & \c & \c & \c & \c & \c & \c & \c & \c & \c & \c & \c & \c
& \c &  \c  & \c & \c & \c & \c \\
a_{\left(4,2\right)}  & \hbox{\ \ } & \hbox{\ \ } & \hbox{\ \ } & \hbox{\ \ }
& \hbox{\ \ } & \hbox{\ \ } & \hbox{\ \ } & b_{\left(4,2\right)}  & \hbox{\ \
} & \hbox{\ \ } & \hbox{\ \ } & \hbox{\ \ } & \hbox{\ \ } & \hbox{\ \ } &
\hbox{\ \ } & c_{\left(4,2\right)} & \hbox{\ \ } & \hbox{\ \ }  \cr
\downarrow & \hbox{\ \ } & \hbox{\ \ } & \hbox{\ \ } & \hbox{\ \ } & \hbox{\ \
} & \hbox{\ \ } & \swarrow & \downarrow & \hbox{\ \ } & \hbox{\ \ } &
\hbox{\ \ } & \hbox{\ \ } & \hbox{\ \ } & \hbox{\ \ } & \swarrow & \downarrow &
\searrow & \hbox{\ \ } \cr
a_{\left(3^2\right)} & \hbox{\ \ } &\hbox{\ \ } & \hbox{\ \ } & \hbox{\ \ } &
\hbox{\ \ } & a_{\left(3^2\right)} & \hbox{\ \ } & b_{\left(3^2\right)} &
\hbox{\ \ } & \hbox{\ \ } & \hbox{\ \ } & \hbox{\ \ } & \hbox{\ \ } &
a_{\left(3^2\right)} & \hbox{\ \ } & b_{\left(3^2\right)} & \hbox{\ \ } &
c_{\left(3^2\right)} \cr
\endmatrix \right.
$$
Here the arrow
$a\rightarrow b$ means that the isomers $a$ and $b$ are neighbours with $a>b$
and $b$ can be obtained from $a$ via a simple substitution reaction. The
K\"orner's diagrams yield that $a_{\left(4,2\right)}$ represents the para
compound,
$b_{\left(4,2\right)}$ represents the ortho compound,
$c_{\left(4,2\right)}$ represents the meta compound,
$a_{\left(3^2\right)}$ represents the asymmetrical compound,
$b_{\left(3^2\right)}$ represents the vicinal compound, and
$c_{\left(3^2\right)}$ represents the symmetrical compound.

\vskip 13pt

6.4. Here we shall discuss the derivatives of ethene, $C_2H_4$,
and their genetic relations, taking into account
the exposition from [3, VI].
The group $G$ of substitution isomerism of
ethene is the Klein subgroup of $S_4$:
$$
G=\{(1), (12)(34), (13)(24), (14)(23)\}.
$$
Since there are no chiral pairs, $G^\prime=G$. For the
group $G^{\prime\prime}$ we can choose any one of the three conjugated Sylow
$2$-subgroups of $S_4$, for instance
$$
G^{\prime\prime}=\{(1), (12)(34), (13)(24), (14)(23), (13), (24), (1234),
(1432)\}.
$$

The group $G^{\prime\prime}$ coincides with the dihedral group $D_4=\langle
r,s\rangle$, where $r=(1234)$, and $s=(13)$. Thus $r^2=
(13)(24)$, $r^3=(1432)$, $sr=(12)(34)$, $sr^2=(24)$, $sr^3=
(14)(23)$.

These groups are defined in [3, VI] by using the inequalities
(1.5.3) -- (1.5.5).

The Abelian group $G$ has
four one-dimensional characters:
The unit character,
the character $\chi_1$ with kernel $\langle (13)(24)\rangle$,
the character $\chi_2$
with kernel $\langle (12)(34)\rangle$, and
the character $\chi_3$ with kernel $\langle (14)(23)\rangle$.

Let $\Sigma$ be the two carbon atom skeleton of ethene.

Case 1. $\lambda=(4)$.

Then
$N_{\left(4\right);\Sigma}=n_{\left(4\right);G}=1$, and
$N_{\left(4\right);\Sigma}^{\prime\prime}=n_{\left(4\right);G^{\prime\prime}}=1
$.

We have
$T_{\left(4\right);G}=T_{\left(4\right),G^{\prime\prime}}=\{a_{\left(4\right)}
\}$ , where
$a_{\left(4\right)}$ is the only $G$-, and $G^{\prime\prime}$-orbit of the
tabloid
$A^{\left(4\right)}=(\{1,2,3,4\},\emptyset,\emptyset,\emptyset)$. The only
$G$-orbit $a_{\left(4\right)}$ represents the parent substance of ethene.

Case 2. $\lambda=(3,1)$.

In this case, again
$N_{\left(3,1\right);\Sigma}=n_{\left(3,1\right);G}=1$, and
$N_{\left(3,1\right);\Sigma}^{\prime\prime}=
n_{\left(3,1\right);G^{\prime\prime}}=1$.

We have
$T_{\left(3,1\right);G}=T_{\left(3,1\right),G^{\prime\prime}}=
\{a_{\left(3,1\right)}\}$, where $a_{\left(3,1\right)}$
is the only $G$-, and $G^{\prime\prime}$-orbit of the
tabloid
$A^{\left(3,1\right)}=(\{1,2,3\},\{4\},\emptyset,\emptyset)$.
This is because both $G$, and $G^{\prime\prime}$ are transitive subgroups of
$S_4$.

Moreover,
$a_{\left(3,1\right)} < a_{\left(4\right)}$, since
$A^{\left(3,1\right)}<A^{\left(4\right)}$.

Case 3. $\lambda=(2^2)$.

Then
$N_{\left(2^2\right);\Sigma}=n_{\left(2^2\right);G}=3$, and
$N_{\left(2^2\right);\Sigma}^{\prime\prime}=
n_{\left(2^2\right);G^{\prime\prime}}=2$ .

We have $T_{\left(2^2\right);G}=\{a_{\left(2^2\right)},b_{\left(2^2\right)},
c_{\left(2^2\right)}\}$, where:

$a_{\left(2^2\right)}$ is the $G$-orbit
$$
\{(\{1,2\},\{3,4\},\emptyset,\emptyset),
(\{3,4\},\{1,2\},\emptyset,\emptyset)\}
$$
of the tabloid
$A^{\left(2^2\right)}=(\{1,2\},\{3,4\},\emptyset,\emptyset)$
with stabilizer
$G_{A^{\left(2^2\right)}}=\langle (12)(34)\rangle$;

$b_{\left(2^2\right)}$ is the $G$-orbit
$$
\{(\{1,4\},\{2,3\},\emptyset,\emptyset),
(\{2,3\},\{1,4\},\emptyset,\emptyset)\}
$$
of the tabloid
$B^{\left(2^2\right)}=(\{1,4\},\{2,3\},\emptyset,\emptyset)$
with stabilizer $G_{B^{\left(2^2\right)}}=\langle (14)(23)\rangle$;

$c_{\left(2^2\right)}$ is the $G$-orbit
$$
\{(\{1,3\},\{2,4\},\emptyset,\emptyset),
(\{2,4\},\{1,3\},\emptyset,\emptyset)\}
$$
of the tabloid
$C^{\left(2^2\right)}=(\{1,3\},\{2,4\},\emptyset,\emptyset)$
with stabilizer $G_{C^{\left(2^2\right)}}=\langle (13)(24)\rangle$.

For the group $G^{\prime\prime}$,
we have
$T_{\left(2^2\right);G^{\prime\prime}}=
\{u_{\left(2^2\right)},v_{\left(2^2\right)}\}$, where:

$u_{\left(2^2\right)}$ is the $G^{\prime\prime}$-orbit
$$
\{(\{1,2\},\{3,4\},\emptyset,\emptyset),
(\{3,4\},\{1,2\},\emptyset,\emptyset),
(\{1,4\},\{2,3\},\emptyset,\emptyset),
(\{2,3\},\{1,4\},\emptyset,\emptyset)\}
$$
of the tabloid
$A^{\left(2^2\right)}=(\{1,2\},\{3,4\},\emptyset,\emptyset)$;

$v_{\left(2^2\right)}$ is the $G^{\prime\prime}$-orbit
$$
\{(\{1,3\},\{2,4\},\emptyset,\emptyset),
(\{2,4\},\{1,3\},\emptyset,\emptyset)\}
$$
of the tabloid
$C^{\left(2^2\right)}=(\{1,3\},\{2,4\},\emptyset,\emptyset)$.

Evidently,
$u_{\left(2^2\right)}=a_{\left(2^2\right)}\cup
b_{\left(2^2\right)}$, and $v_{\left(2^2\right)}=
c_{\left(2^2\right)}$.

Moreover we have,
$$
a_{\left(2^2\right)} < a_{\left(3,1\right)}, \hbox{\ }
b_{\left(2^2\right)} < a_{\left(3,1\right)}, \hbox{\ }
c_{\left(2^2\right)} < a_{\left(3,1\right)},
$$
since
$A^{\left(2^2\right)}<A^{\left(3,1\right)}$,
$(13)(24)B^{\left(2^2\right)}<A^{\left(3,1\right)}$, and
$C^{\left(2^2\right)}<A^{\left(3,1\right)}$, respectively.

Case 4. $\lambda=(2,1^2)$.

Then
$N_{\left(2,1^2\right);\Sigma}=n_{\left(2,1^2\right);G}=3$, and
$N_{\left(2,1^2\right);\Sigma}^{\prime\prime}=
n_{\left(2,1^2\right);G^{\prime\prime}}=2$.

We have
$T_{\left(2,1^2\right);G}=\{a_{\left(2,1^2\right)},b_{\left(2,1^2\right)},
c_{\left(2,1^2\right)}\}$, where:

$a_{\left(2,1^2\right)}$ is the $G$-orbit
$$
\{(\{1,2\},\{3\},\{4\},\emptyset),
\{(\{1,2\},\{4\},\{3\},\emptyset),
(\{3,4\},\{1\},\{2\},\emptyset)\},
(\{3,4\},\{2\},\{1\},\emptyset)\}
$$
of the
tabloid
$A^{\left(2,1^2\right)}=(\{1,2\},\{3\},\{4\},\emptyset)$
with stabilizer $G_{A^{\left(2,1^2\right)}}=\{(1)\}$;

$b_{\left(2,1^2\right)}$ is the $G$-orbit
$$
\{(\{1,4\},\{2\},\{3\},\emptyset),
\{(\{1,4\},\{3\},\{2\},\emptyset),
(\{2,3\},\{1\},\{4\},\emptyset)\},
(\{2,3\},\{4\},\{1\},\emptyset)\}
$$
of the tabloid
$B^{\left(2,1^2\right)}=(\{1,4\},\{2\},\{3\},\emptyset)$
with stabilizer $G_{B^{\left(2,1^2\right)}}=\{(1)\}$;

$c_{\left(2,1^2\right)}$ is the $G$-orbit
$$
\{(\{1,3\},\{2\},\{4\},\emptyset),
\{(\{1,3\},\{4\},\{2\},\emptyset),
(\{2,4\},\{1\},\{3\},\emptyset)\},
(\{2,4\},\{3\},\{1\},\emptyset)\}
$$
of the
tabloid
$C^{\left(2,1^2\right)}=(\{1,3\},\{2\},\{4\},\emptyset)$
with stabilizer $G_{C^{\left(2,1^2\right)}}=\{(1)\}$.

For the group $G^{\prime\prime}$
we have
$T_{\left(2,1^2\right);G^{\prime\prime}}=
\{u_{\left(2,1^2\right)},v_{\left(2,1^2\right)}\}$, where:

$u_{\left(2,1^2\right)}$ is the $G^{\prime\prime}$-orbit
$$
\{(\{1,2\},\{3\},\{4\},\emptyset),
(\{1,2\},\{4\},\{3\},\emptyset),
(\{3,4\},\{1\},\{2\},\emptyset),
(\{3,4\},\{2\},\{1\},\emptyset),
$$
$$
(\{1,4\},\{2\},\{3\},\emptyset),
(\{1,4\},\{3\},\{2\},\emptyset),
(\{2,3\},\{1\},\{4\},\emptyset),
(\{2,3\},\{4\},\{1\},\emptyset)\}
$$
of the tabloid
$A^{\left(2,1^2\right)}=(\{1,2\},\{3\},\{4\},\emptyset)$;

$v_{\left(2,1^2\right)}$ is the $G^{\prime\prime}$-orbit
$$
\{(\{1,3\},\{2\},\{4\},\emptyset),
(\{1,3\},\{4\},\{2\},\emptyset),
(\{2,4\},\{1\},\{3\},\emptyset),
(\{2,4\},\{3\},\{1\},\emptyset)\}
$$
of the
tabloid
$C^{\left(2,1^2\right)}=(\{1,3\},\{2\},\{4\},\emptyset)$.

Clearly,
$u_{\left(2,1^2\right)}=a_{\left(2,1^2\right)}\cup
b_{\left(2,1^2\right)}$, and
$v_{\left(2,1^2\right)}=c_{\left(2,1^2\right)}$.

Moreover,
$$
a_{\left(2,1^2\right)} < a_{\left(2^2\right)},\hbox{\ }
b_{\left(2,1^2\right)} < b_{\left(2^2\right)},\hbox{\ }
c_{\left(2,1^2\right)} < c_{\left(2^2\right)},
$$
since
$A^{\left(2,1^2\right)}<A^{\left(2^2\right)}$,
$B^{\left(2,1^2\right)}<B^{\left(2^2\right)}$,
$C^{\left(2,1^2\right)}<C^{\left(2^2\right)}$,
and
$$
a_{\left(2,1^2\right)} < a_{\left(3,1\right)}, \hbox{\ }
b_{\left(2,1^2\right)} < a_{\left(3,1\right)}, \hbox{\ }
c_{\left(2,1^2\right)} < a_{\left(3,1\right)},\tag 6.4.1
$$
because
$A^{\left(2,1^2\right)}<A^{\left(3,1\right)}$,
$(12)(34)B^{\left(2,1^2\right)}<A^{\left(3,1\right)}$, and
$C^{\left(2,1^2\right)}<A^{\left(3,1\right)}$.

Case 4. $\lambda=(1^4)$.

Then
$N_{\left(1^4\right);\Sigma}=n_{\left(1^4\right);G}=6$, and
$N_{\left(1^4\right);\Sigma}^{\prime\prime}=
n_{\left(1^4\right);G^{\prime\prime}}=3$.

We have
$$
T_{\left(1^4\right);G}=
\{a_{\left(1^4\right)},
b_{\left(1^4\right)}, c_{\left(1^4\right)}, e_{\left(1^4\right)},
f_{\left(1^4\right)}, h_{\left(1^4\right)}\},
$$
where:

$a_{\left(1^4\right)}$ is the $G$-orbit
$$
\{(\{1\},\{2\},\{3\},\{4\}),
(\{2\},\{1\},\{4\},\{3\}),
(\{3\},\{4\},\{1\},\{2\}),
(\{4\},\{3\},\{2\},\{1\})\}
$$
of the
tabloid
$A^{\left(1^4\right)}=(\{1\},\{2\},\{3\},\{4\})$
(the right coset $G$ of $S_4$ modulo $G$);

$b_{\left(1^4\right)}$ is the $G$-orbit
$$
\{(\{1\},\{2\},\{4\},\{3\}),
(\{2\},\{1\},\{3\},\{4\}),
(\{3\},\{4\},\{2\},\{1\}),
(\{4\},\{3\},\{1\},\{2\})\}
$$
of the
tabloid
$B^{\left(1^4\right)}=(\{1\},\{2\},\{4\},\{3\})$
(the right coset $G(34)$ of $S_4$ modulo $G$);

$c_{\left(1^4\right)}$ is the $G$-orbit
$$
\{(\{1\},\{4\},\{2\},\{3\}),
(\{2\},\{3\},\{1\},\{4\}),
(\{3\},\{2\},\{4\},\{1\}),
(\{4\},\{1\},\{3\},\{2\})\}
$$
of the
tabloid
$C^{\left(1^4\right)}=(\{1\},\{4\},\{2\},\{3\})$
(the right coset $G(243)$ of $S_4$ modulo $G$);

$e_{\left(1^4\right)}$ is the $G$-orbit
$$
\{(\{1\},\{3\},\{2\},\{4\}),
(\{2\},\{4\},\{1\},\{3\}),
(\{3\},\{1\},\{4\},\{2\}),
(\{4\},\{2\},\{3\},\{1\})\},
$$
of the
tabloid
$E^{\left(1^4\right)}=(\{1\},\{3\},\{2\},\{4\})$
(the right coset $G(23)$ of $S_4$ modulo $G$);

$f_{\left(1^4\right)}$ is the $G$-orbit
$$
\{(\{3\},\{1\},\{2\},\{4\}),
(\{4\},\{2\},\{1\},\{3\}),
(\{1\},\{3\},\{4\},\{2\}),
(\{2\},\{4\},\{3\},\{1\})\},
$$
of the
tabloid
$F^{\left(1^4\right)}=(\{3\},\{1\},\{2\},\{4\})$
(the right coset $G(132)$ of $S_4$ modulo $G$);

$h_{\left(1^4\right)}$ is the $G$-orbit
$$
\{(\{3\},\{2\},\{1\},\{4\}),
(\{4\},\{1\},\{2\},\{3\}),
(\{1\},\{4\},\{3\},\{2\}),
(\{2\},\{3\},\{4\},\{1\})\},
$$
of the
tabloid
$H^{\left(1^4\right)}=(\{3\},\{2\},\{1\},\{4\})$
(the right coset $G(13)$ of $S_4$ modulo $G$).

For the group $G^{\prime\prime}$,
we have
$$
T_{\left(1^4\right);G^{\prime\prime}}=
\{u_{\left(1^4\right)},v_{\left(1^4\right)},
w_{\left(1^4\right)}\},
$$
where:
$u_{\left(1^4\right)}$ is the $G^{\prime\prime}$-orbit
$$
\{(\{1\},\{2\},\{3\},\{4\}),
(\{2\},\{1\},\{4\},\{3\}),
(\{3\},\{4\},\{1\},\{2\}),
(\{4\},\{3\},\{2\},\{1\}),
$$
$$
(\{3\},\{2\},\{1\},\{4\}),
(\{4\},\{1\},\{2\},\{3\}),
(\{1\},\{4\},\{3\},\{2\}),
(\{2\},\{3\},\{4\},\{1\})\},
$$
of the
tabloid
$A^{\left(1^4\right)}=(\{1\},\{2\},\{3\},\{4\})$
(the right coset $G^{\prime\prime}$ of $S_4$ modulo $G^{\prime\prime}$);

$v_{\left(1^4\right)}$ is the $G^{\prime\prime}$-orbit
$$
\{(\{1\},\{2\},\{4\},\{3\}),
(\{2\},\{1\},\{3\},\{4\}),
(\{3\},\{4\},\{2\},\{1\}),
(\{4\},\{3\},\{1\},\{2\}),
$$
$$
(\{1\},\{4\},\{2\},\{3\}),
(\{2\},\{3\},\{1\},\{4\}),
(\{3\},\{2\},\{4\},\{1\}),
(\{4\},\{1\},\{3\},\{2\})\},
$$
of the
tabloid
$B^{\left(1^4\right)}=(\{1\},\{2\},\{4\},\{3\})$
(the right coset $G^{\prime\prime}(34)$ of $S_4$ modulo $G^{\prime\prime}$);

$w_{\left(1^4\right)}$ is the $G^{\prime\prime}$-orbit
$$
\{(\{1\},\{3\},\{2\},\{4\}),
(\{2\},\{4\},\{1\},\{3\}),
(\{3\},\{1\},\{4\},\{2\}),
(\{4\},\{2\},\{3\},\{1\}),
$$
$$
(\{3\},\{1\},\{2\},\{4\}),
(\{4\},\{2\},\{1\},\{3\}),
(\{1\},\{3\},\{4\},\{2\}),
(\{2\},\{4\},\{3\},\{1\})\},
$$
of the
tabloid
$E^{\left(1^4\right)}=(\{1\},\{3\},\{2\},\{4\})$
(the right coset $G^{\prime\prime}(23)$ of $S_4$ modulo $G^{\prime\prime}$).

Clearly,
$u_{\left(1^4\right)}=a_{\left(1^4\right)}\cup
h_{\left(1^4\right)}$,
$v_{\left(1^4\right)}=b_{\left(1^4\right)}\cup
c_{\left(1^4\right)}$,
$w_{\left(1^4\right)}=e_{\left(1^4\right)}\cup
f_{\left(1^4\right)}$.

Further we have
$$
a_{\left(1^4\right)} < a_{\left(2,1^2\right)}, \hbox{\ }
b_{\left(1^4\right)} < a_{\left(2,1^2\right)}, \hbox{\ }
c_{\left(1^4\right)} < b_{\left(2,1^2\right)},
$$
$$
e_{\left(1^4\right)} < b_{\left(2,1^2\right)}, \hbox{\ }
f_{\left(1^4\right)} < c_{\left(2,1^2\right)}, \hbox{\ }
h_{\left(1^4\right)} < c_{\left(2,1^2\right)},
$$
because
$A^{\left(1^4\right)} < A^{\left(2,1^2\right)}$,
$B^{\left(1^4\right)} < A^{\left(2,1^2\right)}$,
$C^{\left(1^4\right)} < B^{\left(2,1^2\right)}$,
$(14)(23)E^{\left(1^4\right)} < B^{\left(2,1^2\right)}$,
$F^{\left(1^4\right)} < C^{\left(2,1^2\right)}$, and
$H^{\left(1^4\right)} < C^{\left(2,1^2\right)}$.

Here is the diagram which represents the derivatives of ethene.
$$
\left.
\matrix \format \c & \c & \c & \c & \c & \c & \c & \c & \c & \c & \c & \c & \c
& \c &  \c  \\
(4) & \hbox{\ \ } & \hbox{\ \ } & \hbox{\ \ } & \hbox{\ \ } & \hbox{\ \ } &
\hbox{\ \ } & \hbox{\ \ } & a_{\left(4\right)} & \hbox{\ \ } & \hbox{\ \ } &
\hbox{\ \ } & \hbox{\ \ } & \hbox{\ \ } & \hbox{\ \ } \cr
\downarrow & \hbox{\ \ } & \hbox{\ \ } & \hbox{\ \ } & \hbox{\ \ } & \hbox{\ \
} & \hbox{\ \ } & \hbox{\ \ } & \downarrow & \hbox{\ \ } & \hbox{\ \ } &
\hbox{\ \ } & \hbox{\ \ } & \hbox{\ \ } & \hbox{\ \ } \cr
(3,1) & \hbox{\ \ } & \hbox{\ \ } & \hbox{\ \ } & \hbox{\ \ } & \hbox{\ \ } &
\hbox{\ \ } & \hbox{\ \ } & a_{\left(3,1\right)} & \hbox{\ \ } & \hbox{\ \ } &
\hbox{\ \ } & \hbox{\ \ } & \hbox{\ \ } & \hbox{\ \ } \cr
\downarrow & \hbox{\ \ } & \hbox{\ \ } & \hbox{\ \ } & \hbox{\ \ } & \hbox{\ \
} & \hbox{\ \ } & \swarrow & \downarrow & \searrow & \hbox{\ \ } & \hbox{\ \ }
& \hbox{\ \ } & \hbox{\ \ } & \hbox{\ \ } \cr
(2^2) & \hbox{\ \ } & \hbox{\ \ } & \hbox{\ \ } & \hbox{\ \ } & \hbox{\ \ } &
a_{\left(2^2\right)} & \buildrel{u_{\left(2^2\right)}}\over {\leftrightarrow} &
b_{\left(2^2\right)} & \hbox{\ \ } & c_{\left(2^2\right)} & \hbox{\ \ } &
\hbox{\ \ } & \hbox{\ \ } & \hbox{\ \ } \cr
\downarrow & \hbox{\ \ } & \hbox{\ \ } & \hbox{\ \ } & \hbox{\ \ } & \swarrow
& \hbox{\ \ } & \hbox{\ \ } & \downarrow & \hbox{\ \ } & \hbox{\ \ } & \searrow
& \hbox{\ \ } & \hbox{\ \ } & \hbox{\ \ } \cr
(2,1^2) & \hbox{\ \ } & \hbox{\ \ } & \hbox{\ \ } & a_{\left(2,1^2\right)} &
\hbox{\ \ } & \buildrel{u_{\left(2,1^2\right)}}\over {\leftrightarrow} &
\hbox{\ \ } & b_{\left(2,1^2\right)} & \hbox{\ \ } & \hbox{\ \ } & \hbox{\ \ }
& c_{\left(2,1^2\right)} & \hbox{\ \ } & \hbox{\ \ } \cr
\downarrow & \hbox{\ \ } & \hbox{\ \ } & \swarrow & \hbox{\ \ } & \searrow &
\hbox{\ \ } & \hbox{\ \ } & \downarrow & \searrow & \hbox{\ \ } & \hbox{\ \ } &
\downarrow & \searrow & \hbox{\ \ } \cr
(1^4) & \hbox{\ \ } & a_{\left(1^4\right)} & \hbox{\ \ } & \hbox{\ \ } &
\hbox{\ \ } & b_{\left(1^4\right)} & \buildrel{v_{\left(1^4\right)}}\over
{\leftrightarrow} & c_{\left(1^4\right)} & \hbox{\ \ } & e_{\left(1^4\right)}
&\buildrel{w_{\left(1^4\right)}}\over {\leftrightarrow} & f_{\left(1^4\right)}
& \hbox{\ \ } & h_{\left(1^4\right)} \cr
\hbox{\ \ } & \hbox{\ \ } & \hbox{\ \ } & \nwarrow & \underline{\ \ \ } &
\underline{\ \ \ } & \underline{\ \ \ } & \underline{\ \ \ } & \underline{\ \ \
} & \underline{\ \ \ } & \underline{\ \ \ } & \underline{\ \ \ } & \underline{\
\ \ } & \nearrow & \hbox{\ \ } \cr
\hbox{\ \ } & \hbox{\ \ } & \hbox{\ \ } & \hbox{\ \ } & \hbox{\ \ } & \hbox{\ \
} & \hbox{\ \ } & \hbox{\ \ } & u_{\left(1^4\right)} & \hbox{\ \ } & \hbox{\ \
} & \hbox{\ \ } & \hbox{\ \ } & \hbox{\ \ } & \hbox{\ \ } \cr
\endmatrix \right.
$$
The arrow
$a\rightarrow b$ means that the isomers $a$ and $b$ are neighbours with $a>b$
and $b$ can be obtained from $a$ via a simple substitution reaction.
The horizontal double arrow means that the two isomers are diastereomers
and the
letter above/below it denotes the corresponding structural isomer. The above
diagram does not indicates the simple substitution reactions  from (6.4.1),
where the isomers are not neighbours.

Our extended approach confirms the conclusion of Lunn and Senior from [3, VI]
that there are no type properties which distinguish the members of the
pairs of diastereomers
$a_{\left(2^2\right)},
b_{\left(2^2\right)}$ and
$a_{\left(2,1^2\right)},
b_{\left(2,1^2\right)}$.
It is clear that the genetic relations from the above diagram fail to make any
difference between
them. On the level of one-dimensional characters of the group $G$, the members
of the second pair are indistinguishable because the stabilizers of their
elements coincide with the unit group. At first
sight each one of the characters $\chi_2$ and $\chi_3$ of the group $G$ \lq\lq
distinguishes" $a_{\left(2^2\right)}$  and
$b_{\left(2^2\right)}$: For instance $\chi_2$ is identically $1$ on the
stabilizer
$G_{A^{\left(2^2\right)}}=\langle (12)(34)\rangle$
of the tabloid $A^{\left(2^2\right)}\in
a_{\left(2^2\right)}$ and $\chi_2$ is not identically $1$ on the stabilizer
$G_{B^{\left(2^2\right)}}=\langle (14)(23)\rangle$ of the tabloid
$B^{\left(2^2\right)}\in
b_{\left(2^2\right)}$. The same is true for $\chi_3$ if we  replace
$a_{\left(2^2\right)}$ for
$b_{\left(2^2\right)}$ and vice versa. It is not hard to check that the
presence of a non-trivial one-dimensional character
$\theta$ of $S_{\left(2^2\right)}$ in the formula (5.1.2) has the same effect
as the interchange of $\chi_2$ and $\chi_3$.
Unfortunately, the characters $\chi_2$ and $\chi_3$ can
not be distinguished: Each one of them can be obtained from the other by a
special automorphism of the group $G$, that is induced by a renumbering the
unsatisfied valences of the skeleton $\Sigma$. Therefore we can only conclude
that both $a_{\left(2^2\right)}$ and
$b_{\left(2^2\right)}$ are elements of the symmetric difference
$$
(T_{\left(2^2\right);\chi_2}\backslash
T_{\left(2^2\right);\chi_3})\cup
(T_{\left(2^2\right);\chi_3}\backslash
T_{\left(2^2\right);\chi_2}),
$$
so
the type properties corresponding to
$\chi_2$ and $\chi_3$ via the Extended Lunn-Senior Thesis 1.6.1 can not be used
to make difference between the members of this pair of diastereomers.
Thus, \lq\lq ...Whenever diameric pairs of
disubstitution derivatives of ethylene have been investigated, it has been
necessary to fall back on the specific properties of the molecules in question
in order to decide which one is the cis and which one the trans isomer" (see
[3, VI]).

\heading{Appendix A}

\endheading

Let $X$ be a set with partial order
$\leq$.
This means that the binary relation
$\leq$ defined in $X$, satisfies the following two properties:

(a) $x\leq x$ for any $x\in X$;

(b) $x\leq y$ and $y\leq z$ implies $x\leq z$ for all $x$, $y$, $z\in X$.

We write $x<y$ when $x\leq y$ and $x\neq y$.

Any subset $Y\subset X$ inherits the structure of partially ordered set from
$X$.

For any $x$, $y\in X$, the set $\{z\in X|x\leq z\leq y\}$ is
denoted by $[x,y]$, and is called {\it closed interval in $X$ with endpoints
$x$ and $y$}. The set $\{z\in X|x<z<y\}$ is denoted by $(x,y)$, and
is called {\it open interval in $X$}.

When $x<y$ and $(x,y)=\emptyset$, then $x$ and $y$ are said to be {\it
neighbours} in $X$ with $x<y$. The elements $x$ and $y$
are called {\it neighbours}
in $X$ if $x$ and $y$ are neighbours in $X$ with $x<y$, or with $y<x$.

\heading{Resume of Part I}

\endheading

There are four themes in this paper, which may be of interest to a chemist:

1) The determination of the structural formula of a potentially existing isomer
with given skeleton $\Sigma$,
starting from any tabloid in the $G$-orbit which represents this
isomer according to Lunn-Senior thesis 1.5.1 (here $G$ is the symmetry group
of $\Sigma$);

2) the partial order $\leq$ on the set of all $G$-orbits of tabloids (Section 4
and the subsidiary Sections 1 - 3);

3) the hypothesis that the set of $(\chi,\theta)$-orbits determines a type
property of the molecule under consideration (Section 1, 1.6.1), and
the count of $(\chi,\theta)$-orbits (Section 5);

4) the attempt to breathe new life into the philosophy of the original
Lunn-Senior's paper [3].

The way of construction of the structural formula of an isomer with given
skeleton $\Sigma$ is explicitly
build in the representation of this isomer by a tabloid $A=(A_1, A_2,\hdots )$:
If $i\in A_k$, then we attach the univalent substituent $x_k$ to
$\Sigma$'s unsatisfied valence number $i$, for $k=1, 2,\hdots$ . Since there is
no \lq\lq canonical" numbering of the unsatisfied valences, a
main
problem of the present model is the identification of the real substances (if
any) having these structural formulae, in terms of the model itself. The
partial
order $\leq$, and the $(\chi,\theta)$-orbits can be applied for this problem to
be solved (at least partially).

The partial order may also be used in the following way:

The relation $a < b$ between the isomers $a$ and $b$ is an indication of
the existence of a finite sequence of simple substitution reactions
$b \rightarrow  c_1\rightarrow\cdots \rightarrow c_r \rightarrow a$, where the
compounds $c_1,\hdots, c_r$, are intermediate stages in a synthesis of $a$.
Such a sequence
$c_1,\hdots, c_r$ (which is far-away of being unique), can be constructed by
means of Theorem 3.4.4.

The relation $a\not\leq b$ implies that the isomer $a$ for sure can not be
obtained from the isomer $b$ via a finite sequence of simple substitution
reactions.

The partial order is tested in Section 6 for finding the genetic relations of
the substitution derivatives of ethene. It is applied also in the case of
di-, and tri-substitution derivatives of benzene and yields the classical
K\"orner's
relations. These two applications are considered also in Lunn-Senior's paper,
Part VI. It goes without saying that the adequacy of this partial
order to the chemical reality needs more experimental verifications.

A central topic in the paper is a detailed study of the notion of
"neighbourhood" with respect to the above partial order.  If two isomers $a$
and $b$ are neighbours with $a<b$, then probably there exists a chemical
reaction $b \rightarrow  a$,
but it is certain that this reaction can not be represented
as $b \rightarrow c \rightarrow  a$, where $c$ is a isomer.
The main result in this direction
is Theorem 4.2.3, (ii), which characterizes mathematically the pairs
of neighbours $a<b$,
and in this case predicts the existence of a chain
$b \rightarrow  c_1\rightarrow\cdots \rightarrow c_r \rightarrow a$,
where the
intermediate "reactions" are "virtual", that is, $c_1,\hdots ,c_r$ are not
represented by tabloids, but by ordered dissections.

Item 5 of the Extended Lunn-Senior Thesis 1.6.1 is our hypothesis.  If
$\chi$ is a one-dimensional character of the symmetric group $G$ of the
molecule, and if $\theta$ is a one-dimensional character of the group
$S_\lambda$ (this group reflects the empirical formula (1.1.2) of the
molecule), then the couple $(\chi,\theta)$ produces via condition (5.1.3) a
subset of the set of all $G$-orbits, which, we suppose, represents a type
property of this molecule.  This is true when $\chi=1_G$, and
$\theta=1_{S_\lambda}$.  In this particular case we obtain the set of all
$G$-orbits, each one of
them possibly representing an isomer due to Lun-Senior Thesis 1.5.1.  This also
is true in case $\chi=\chi_e$, and $\theta=1_{S_\lambda}$ (see  Section 6,
6.1), and we get a set which represents the chiral pairs.  Theorem 5.3.1 is a
wide generalization of a crucial result of Ruch (see  Theorem 6.2.1) which
connects the existence of chiral pairs with the dominance order among the
partitions.  This theorem holds out a hope that the Extended Lunn-Senior Thesis
1.6.1 is valid.  Which couples $(\chi,\theta)$ are within the scope of 1.6.1,
item 5, is a matter of the experiment.  We guess that there are no exceptions.
Theorem 5.2.7 gives an explicit formula for the number of the
$(\chi,\theta)$-orbits.

\heading{Acknowledgement}

\endheading

I would like to thank Prof. Adalbert Kerber --- my referee who saw no reason to
keep from me his identity. As a result we had an e-mail discussion which has
greatly benefited the present paper.

\heading {References}
\endheading

\noindent [1] V. V. Iliev, A generalization of Redfield's Master Theorem,
xxx.lanl.gov/abs/math.RT/ 9902089.

\noindent [2] G. James, A. Kerber, The Representation Theory of the Symmetric
Group, {\it in} Encyclopedia of Mathematics and its Applications, Vol. 16,
Addison-Wesley Publishing Company, 1981.

\noindent [3] A. C. Lunn, J. K. Senior, Isomerism and Configuration, J.
Phys. Chem. 33 (1929), 1027 - 1079.

\noindent [4] I. G. Macdonald,  Symmetric Functions and Hall
Polynomials, Clarendon Press, Oxford, 1995.

\noindent [5] G. P\'olya, Kombinatorische Anzahlbestimmungen f\"ur Gruppen,
Graphen und chemische Verbindungen, Acta Math. 68 (1937), 145 --
254. English translation:
G. P\'olya and R. C. Read, Combinatorial Enumeration of Groups,
Graphs and Chemical Compounds, Springer-Verlag New York Inc., 1987.

\noindent [6] E. Ruch, W. H\"asselbarth, B. Richter, Doppelnebenklassen
als Klassenbegriff und No\-menklaturprinzip f\"ur Isomere und ihre
Abz\"ahlung, Theoret. chim. Acta (Berl.), 19 (1970), 288 -- 300.

\noindent [7] J.-P. Serre, Repr\'esentations Lin\'eaires des Groupes
Finis, Hermann, Paris, 1967.

\end